\documentclass{article}
\usepackage[left=1in,top=1in,right=1in,bottom=1in,nohead,footskip=15pt]{geometry}

\usepackage[pdftex, colorlinks=true, urlcolor=MyDarkBlue, citecolor=MyDarkRed, linkcolor=MyDarkGreen]{hyperref}
\usepackage[utf8]{inputenc}
\usepackage{amsmath,amsthm,epsf,amscd,amssymb,verbatim}
\usepackage[square,sort,comma,numbers]{natbib}
\usepackage{graphicx}
\usepackage{mathrsfs}
\usepackage{color}
\usepackage{titlesec}
\usepackage{enumerate}
\usepackage[shortlabels]{enumitem}
\usepackage{subfigure}
\usepackage{datetime}
\usepackage[normalem]{ulem}
\usepackage[font=footnotesize]{caption}
\usepackage{abstract}

\usepackage{chemarr}
\usepackage{extarrows}
\usepackage[all]{xy}

\usepackage{color}
\definecolor{MyDarkRed}{rgb}{0.5,0,0.1}
\definecolor{MyDarkBlue}{rgb}{0.1,0.1,0.5}
\definecolor{MyDarkGreen}{rgb}{0.1,0.5,0.1}
\definecolor{MyRed}{rgb}{1.0,0,0}
\definecolor{MyBlue}{rgb}{0,0,1.0}
\definecolor{MyGreen}{rgb}{0,0.8,0}
\definecolor{lightgray}{rgb}{0.96,0.96,0.96}
\definecolor{darkgray}{rgb}{0.4,0.4,0.4}

\newtheorem{theorem}{Theorem}
\newtheorem{lemma}{Lemma}
\newtheorem{cor}{Corollary}
\newtheorem{prop}{Proposition}

\theoremstyle{definition}
\newtheorem{definition}{Definition}
\theoremstyle{remark}

\newtheorem{mynote}{Note}

\newenvironment{note}{\begin{quote}\vspace{-0.1in}\begin{mynote}}{\end{mynote}\end{quote}}

\newcommand{\wipe}[1]{}

\newcommand{\subsubsubsection}[1]{\vspace{0.1in}\noindent \emph{\textbf{#1}:}}
\newcommand{\imphead}[1]{\vspace{0.07in}\noindent\emph{{#1}:}}
\newcommand{\mathquote}[1]{{\text{`}{#1}\text{'}}}

\newcommand{\zA}{{\mathsf{A}}}
\newcommand{\zB}{{\mathsf{B}}}
\newcommand{\zC}{{\mathsf{C}}}
\newcommand{\zD}{{\mathsf{D}}}
\newcommand{\zE}{{\mathsf{E}}}
\newcommand{\zF}{{\mathsf{F}}}
\newcommand{\zG}{{\mathsf{G}}}

\newcommand{\zS}{{\mathsf{S}}}
\newcommand{\zO}{{\mathsf{O}}}

\newcommand{\transitionarrow}{\leadsto}
\newcommand{\ordind}[1]{{\langle {#1} \rangle}}
\newcommand{\inlinebox}[1]{{\fbox{{#1}}}}
\newcommand{\lite}[1]{{\color{darkgray}{#1}}}

\definecolor{MyRed}{rgb}{1.0,0,0}
\definecolor{MyBlue}{rgb}{0,0,1.0}
\definecolor{MyPurple}{rgb}{0.5,0,0.5}

\newcounter{pathnum}
\newcounter{condnum}[pathnum]

\newcommand{\state}[3]{ {[{#1};~ z\in ({#2}, {#3}) ]} }

\ifx\hidespecial\undefined
	\newcommand{\todo}[1]{ {\color{MyRed}\textbf{TODO:} #1} }

\else
	\newcommand{\todo}[1]{}

\fi

\ifx\hideproof\undefined
	\usepackage{framed}
	\definecolor{leftBarGray}{rgb}{0.8,0.8,0.8}
	\newenvironment{myleftbar}{%
	\MakeFramed {\advance\hsize-\width \FrameRestore}}%
	{\endMakeFramed}
	\newenvironment{quoteproof}{

	\begin{myleftbar}\noindent\textit{Proof.}\vspace{0.1in}

	\small}{\qed\end{myleftbar}}
\else
	\usepackage{environ}
	\NewEnviron{quoteproof}{}
\fi

\titleformat{\subsubsection}[runin]{\bfseries}{\thesubsubsection~}{0pt}{\bfseries}[:]

\title{A Classification of Configuration Spaces of Planar Robot Arms\\ with Application to a Continuous Inverse Kinematics Problem}

\author{
Subhrajit Bhattacharya\thanks{{Department of Mathematics, University of Pennsylvania, Philadelphia PA 19104. e-mail: \texttt{subhrabh@math.upenn.edu}. Phone: (001) 267-252-6638. (Corresponding aothor).}}
~~ and ~~
Mihail Pivtoraiko\thanks{{Department of Mechanical Engineering and Applied Mechanics, University of Pennsylvania, Philadelphia PA 19104. e-mail: \texttt{mihailp@seas.upenn.edu}}}
}
\date{}

\begin{document}

\maketitle

\setlength{\absleftindent}{0.1in}
\setlength{\absrightindent}{0.1in}
\begin{abstract}
Using results on the topology of moduli space of polygons~\cite{jaggi:thesis:paper:92,Kapovich94onthe}, it can be shown that for a planar robot arm with $n$ segments there are some values of the base-length, $z$ (\emph{i.e.}, length of line joining the base of the arm with its end-effector), at which the configuration space of the constrained arm (arm with its end effector fixed) has two disconnected components, while there are other values at which the constrained configuration space has one connected component. We first review some of these known results relating the value of $z$ with the connectivity of the constrained configuration space with end-effector fixed.

Then the main design problem addressed in this paper is the construction of pairs of continuous inverse kinematics for arbitrary robot arms, with the property that the two inverse kinematics agree (\emph{i.e.} return the same configuration) when the constrained configuration space has a single connected component, but they give distinct configurations (one in each connected component) when the configuration space of the constrained arm has two components.
This design is made possible by a fundamental theoretical contribution in this paper -- a classification of configuration spaces of robot arms such that the type of path that the system (robot arm) takes through certain critical values of the forward kinematics function is completely determined by the class to which the configuration space of the arm belongs.
This classification result makes the aforesaid design problem tractable, making it sufficient to design a pair of inverse kinematics for each class of configuration spaces (three of them in total).

The motivation for this work comes from a more extensive problem of motion planning for the end effector of a robot arm, in which the ability to continuously sample one configuration from each connected component of the constrained configuration spaces of the arm enables us to dramatically reduce the dimensionality of the space in which the planning has to be performed, without sacrificing completeness guarantees. We start the paper with the general motivation, but address only the problem of sampling such configurations when there is no obstacle in the environment -- a problem that in itself is non-trivial. Incorporating obstacles and a complete graph search-based planning for the end effector using this technique will be presented in a future papers.

We demonstrate the simplicity and the low complexity of the presented algorithm through a Javascript + HTML5 based implementation available at \url{http://hans.math.upenn.edu/~subhrabh/nowiki/robot_arm_JS-HTML5/arm.html}.

\end{abstract}


\vspace{0.3in}
\noindent\textbf{Acknowledgements:} 
We would like to thank Prof. Vijay Kumar, Department of Mechanical Engineering and Applied Mechanics, University of Pennsylvania, and Prof. Robert Ghrist, Department of Mathematics, University of Pennsylvania, for their valuable suggestions and insightful discussions during the course of writing this paper.

\newpage
\section{Introduction and Related Work}

In context of studying the topology of the configuration space of 
polygons (\emph{i.e.}, a robot arm with a fixed end-effector position) 
some remarkable results were proposed in the PhD thesis of B. Jaggi~\cite{jaggi:thesis:paper:92}. Later on, some of those results were rediscovered and enhanced in \cite{Kapovich94onthe} and also reported in \cite{Milgram04thegeometry}. 
While these beautiful mathematical results are worth exploring in their own rights, in this paper we use them as a first step towards solving a very practical path planning problem for the end effector of a planar robot arm. 
To that end we discover a classification of configuration spaces of planar robot arms (Proposition~\ref{prop:sys-paths} in this paper) based on how the 
forward kinematics function for the arm passes through certain critical values. 
This insight then serves as a key to solving the design problem proposed in this paper -- construction of pairs of continuous inverse kinematics that give one unique configuration in each connected component of the constrained configuration spaces.

\subsection{Motivation}

\subsubsection{General Motivation}

This paper is motivated by two key needs:
\begin{itemize}

 \item[1.] \emph{The need for an inverse kinematics to be continuous -- That is, as we change the end effector position continuously, the joint angles computed using the IK algorithm should also change continuously and should not have abrupt jumps:}  
 In literature, for planar and spatial arms with a few segments, this has often been achieved by explicit trigonometric formulae developed from the geometry of the specific arm~\cite{Paul:manipulator:81,Murray1994}. However, very often such formulations are limited to the specific arm that the IK is designed for, becomes increasingly complex with the increase in the number of segments, and does suffer from isolated singularities and/or discontinuities.
 For more complex arms, traditionally numerical gradient-decent type algorithms have been used~\cite{Gupta:IK:85,Lumelsky:IK:84,Goldenberg:IK:85,Wang:CCD:91}. However, the problem very often being nonlinear and non-convex, guarantees of completeness or even continuity are difficult to achieve. Furthermore, numerical techniques are often computationally expensive.
 A mixed numerical and analytic technique has been used in \cite{Trinkle:complete:04} for computing inverse kinematics in context of path planning for the end effector of an arm with $2$ segments and point obstacles in the environment. 
 A closed form solutions to the inverse kinematics problem has recently been hinted in \cite{Cajar:IK:07} using a \emph{triangulation} approach. But the problem being addressed for a spatial arm, some numerical techniques were also adopted.

 The IK algorithm that we propose and use in this paper is an analytically computable one of $O(n^2)$ complexity ($n$ being the number of segments of the arm), akin to the triangulation approach of \cite{Cajar:IK:07}. However our additional requirement of identifying and passing through certain critical points in the configuration space, as will be described in the next point, has necessitated a more careful and formal construction of the IK algorithm. In particular, our algorithm has a recursive or incremental flavor to it, wherein we construct the inverse kinematics by breaking up the arm into smaller components.

 \item[2.] \emph{Design of tuples of continuous inverse kinematics which, for a given end effector position, return only one unique configuration in every connected component of the constrained configuration space:}
 The motivation for wanting to achieve this will be clear in the next subsection. This problem is challenging because of the plethora of geometric properties of the configuration space of a robot arm that one is faced with as the number of segments as well as the lengths of the segments are changed. In \cite{jaggi:thesis:paper:92,Kapovich94onthe,Milgram04thegeometry} it was proved that there are only two possible types as far as the connectivity of the constrained configuration space (\emph{i.e.}, the configuration space with end effector position fixed; equivalently the moduli configuration space of a polygon with $(n+1)$ sides) is concerned: Either it has a single connected component, or it a disjoint union of two $(n-2)$-dimensional tori. 
 Due to a fundamental theorem from Morse theory~\cite{Nicolaescu2007morse,milnor1963morse}, a continuous transition from a state in one of these types to one in the other type, that results in the change in the topology (connectivity) of the constrained configuration space, will mean that the system will pass through a critical value of the Morse function which maps configurations to the base length.
 We call the corresponding critical points \emph{vital critical points}. In order to design the pairs of IK with the desired properties, we need to ensure that the system passes through the vital critical points at the vital critical values of the base length.

 The way we handle this design problem is to first identify the possible vital critical points for a given system, and how the system passes through the vital critical points as the base length is changed continuously. One surprising result that we derive (Proposition~\ref{prop:sys-paths}) is that there are only three such possible types of paths passing through the different vital critical points despite the large variety of robot arms that one can come up with. This essentially gives us a classification of configuration spaces of planar robot arms, and thus simplifies our problem of design of the IK pairs significantly --- we only need to consider the three possible classes, and design three corresponding pairs of IKs.

\end{itemize}


\subsubsection{Motivation from a Path Planning Problem} \label{sec:motivation}

Although the aforesaid problem objectives are interesting in their own rights, in this section we elaborate on a more practical motivation for the problem that initiated this line of research in the first place. While this motivation relies largely on the presence of obstacles in the environment, we will also illustrate the non-triviality of the proposed solution even in environments without obstacle. This paper is solely dedicated to the later, while the more general case with obstacles will be considered in a forthcoming paper as discussed in Section~\ref{sec:discussion}.

\begin{figure}[h]
\centering
      \includegraphics[width=0.75\textwidth, trim=130 40 250 120, clip=true]{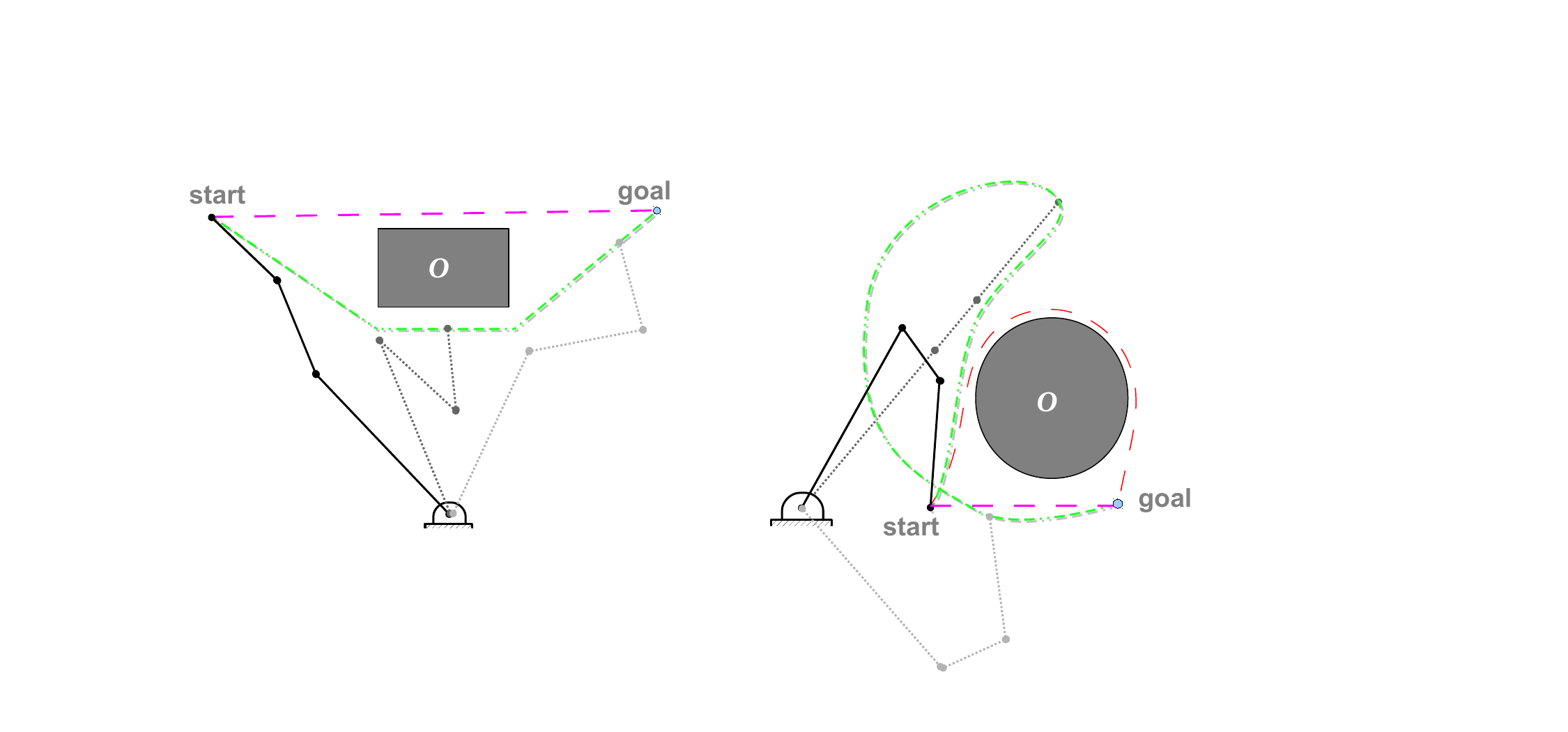}
\caption{The end effector of the robot arm may not follow the trajectories marked in magenta starting from the initial configuration as shown. In the left figure a trajectory in a different homotopy class (shown in green) solves the problem since it can be followed by the end effector. But not so in the case of the right one (the trajectory in red is in a different homotopy class, but still cannot be followed). In this case, a different trajectory (shown in green, which is in the same homotopy class as the magenta one), but passing through a critical configuration, can be followed. This issue is arguably endemic to environments with obstacles. However the proper solution to this problem involves sampling one configuration in each connected component of the constrained configuration spaces for different end-effector positions.
And achieving that solution, even in complete absence of obstacles, involves non-trivial techniques. This paper, the first in a series, deals with solving the problem of continuously sampling one configuration from the constrained configuration spaces for each end-effector position, assuming that there is no obstacle in the environment.}  \label{fig:motivation}
\end{figure}

The practical motivation behind this line of research comes from the recent advances in robotics, computer graphics as well as health sciences relating to kinematics of human motion and tools consting of linkages~\cite{buss20033d,kumar2008kinesthetic}. The problem of inverse kinematics has been actively studied for several decades due to the many promising benefits of these applications.
Often the subject of interest that leads to the inverse kinematics problem is the end-effector of the robot arm, rather than the arm itself. There may be applications in which one may need to plan shortest feasible paths for the end effector from an initial position to a final position (for example, in minimally invasive assisted laparoscopic  surgery using a robotic manipulator~\cite{payne:08:laparoscopic}, where the objective is to reach a tissue of interest inside the body).
Being reduced to a typical path planning problem, a roboticist will be inclined to discretize the configuration space, place a vertex on every feasible cell in the discretization, and establish edges between neighboring vertices, thus creating a graph in which a search-based algorithm such as A* \cite{Hart-Astar} or D* \cite{Stentz-dstar} may be employed. Such an approach has guarantees on completeness and optimality as opposed to continuous gradient-decent like approaches which suffer from issues of local minima and suboptimality~\cite{KhoslaVolpeICRA88,KimKhosla02,ref:Conner03}.
However it is to be remembered that the configuration space of the arm is the entire space of the joint angles, which, for an arm with $n$ segments, is a $n$-dimensional space (a $n$-torus, $\mathbb{T}^n$, with punctures caused by presence of obstacles or infeasible arm configurations). Discretization and graph construction in such a high dimensional configuration space creates an extremely large graph with high degree (number of emanating edges) per vertex. While such graphs have indeed been used in solving the path planning problem for planar robotic arms~\cite{Hansen-anytimeastarjair}, typically a full-blown search being highly expensive, a compromise has to be made in the optimality of the search result or the completeness of the algorithm using randomized techniques like the ARA* algorithm.

One may attempt to simplify the problem by instead trying to plan a trajectory in the $2$-dimensional space of the end-effector, $(\mathbb{R}^2 - \mathcal{O})$ (where $\mathcal{O}$ is the set of obstacles), and hence try to make the arm move (using, for example, a feedback controller) such that the end effector follows the planned trajectory. However very often the arm will fail to follow such a planned trajectory, even though every point on the trajectory is reachable (see Figure~\ref{fig:motivation} and Section 3.6.1 of \cite{PhDthesis:12}).

One may try to resolve this issue by planning multiple trajectories in the different homotopy classes in $(\mathbb{R}^2-\mathcal{O})$ (the multiple homotopy classes being created due to the presence of obstacle) as was discussed in \cite{PhDthesis:12,planning:AURO:12}. However a close investigation reveals that this approach is not enough (for example the case on the right in Figure~\ref{fig:motivation}).

\begin{figure*}[h]
\centering
      \includegraphics[width=0.75\textwidth, trim=0 0 0 0, clip=true]{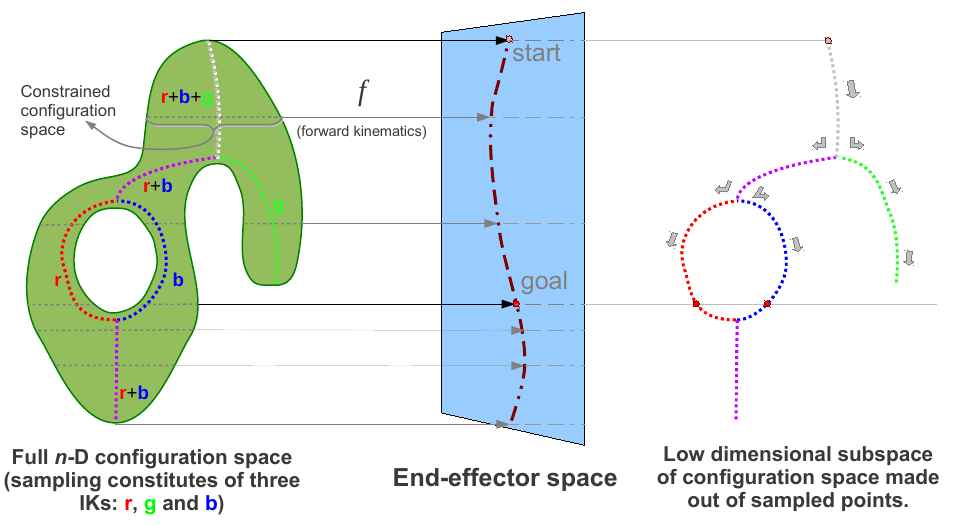}
\caption{(Informal/schematic diagram for illustration) The idea behind creating a faithful low-dimensional representation of the configuration space by sampling one point from each connected component of every constrained configuration space. The sampled subspace can be constructed as union of the image of multiple inverses kinematics (maps from end-effector space to configuration space) passing through certain critical points in the configuration space.}  \label{fig:sampling}
\end{figure*}

\vspace{0.1in}
\noindent We proposed the correct solution approach to be the following:
\begin{quote}
 Given an end-effector position, we \uline{sample exactly one configuration from every connected component of the constrained configuration space} (recall that the \emph{constrained configuration space} is the configuration space of the arm with the end effector fixed, which, in presence of obstacles, will \emph{exclude} configurations in which a segment of the arm intersects an obstacle). 
 As the end effector is moved in a continuous fashion, the sampled points from the configuration space will obviously change, but it should change continuously as well. 
 This leads to the notion of a `proximity' among the different sampled configurations as the end effector position is changed. A graph can hence be created --- \uline{vertices correspond to sampled configurations from constrained configuration spaces for every discrete end-effector positions, and edges are established between \emph{adjacent} configurations} (that is, the corresponding sampled configurations for the \emph{adjacent} end-effector positions).
\end{quote}
Since for every end effector position we sample a finite number of distinct configurations (one for each connected component of the constrained configuration space), the set of all the sampled configuration is a space with the same dimensionality as the end effector space (but with a different topology, akin to the topology of a \emph{covering space}~\cite{Hatcher:AlgTop}). Thus we have reduced the dimensionality of the problem significantly from the one encountered when planning in the entrire $n$-dimensional configuration space of the arm.
It is however worth noting that within a connected component of a constrained configuration space every configuration can be reached from the other keeping the end effector fixed. 
Thus sampling one configuration from each connected component gives us a \uline{\emph{faithful} low-dimensional representation of the $n$-dimensional joint space}. This is illustrated in Figure~\ref{fig:sampling}.

Due to the above, it is in fact easy to observe that there exists a path in the full configuration space of the arm, whose image under the end-effector map (\emph{i.e.}, the map that takes in a configuration and returns the corresponding end-effector position) is a given trajectory, if and only if there also exists a path in the subspace consisting of the sampled points whose image under the end-effector map is the same trajectory. 
Thus we end up reducing the dimensionality of the problem without compromising the completeness guarantees.


\vspace{0.1in}
\noindent \textbf{\emph{The non-triviality even in absence of obstacles in the environment:}} 
Although the illustrations and the rationales described above may indicate that obstacles are the only culprits in making a constrained configuration space have disconnected components and hence the need for the continuous sampling as described, that in fact is not true. Due to the results of \cite{jaggi:thesis:paper:92,Kapovich94onthe}, even when there is no obstacle in the environment, a constrained configuration space may have one or two connected components.
In this paper we design a way of continuously sampling one configuration from each connected component of the constrained configuration spaces as described, and as the first step, we assume that there is no obstacle in the workspace. We moreover assume that there is no joint limit and the the segments of the arm lie in different planes so that there is no self-collision of the arm. 
As we will see, this problem in itself is highly non-trivial. This is the first paper in a planned series, where the eventual objective will be to extend the techniques developed in this paper to environments with obstacles and arms with joint limits, and hence be able to execute graph search algorithms to do quick and efficient planning for robot arms with completeness guarantees.

\subsection{Organization of the Paper}

We start with some well-known results and simple observations in Section~\ref{sec:preliminaries}.
Following this, the main content of the paper is organized into two sections: \emph{i.} Section~\ref{IK-alg} describes a general continuous inverse kinematics map that involves choice of certain parameters, and, \emph{ii.} the main contributions of the paper appear in Section~\ref{sec:sampling-main}, where we describe how to choose the said parameters for designing the IK pairs with the desired property (one unique configuration per connected component of the constrained configuration space).
The later is made possible by the main theoretical contribution of the paper -- the classification result of Proposition~\ref{prop:sys-paths} -- where we identify the possible types of paths that a system can take through the \emph{vital critical values} of the base-length function. The main technical contribution of the paper is summarized in Proposition~\ref{prop:complete-IK}.


%
%
%

\vspace{0.5in}

\section{Preliminaries} \label{sec:preliminaries}

\begin{figure*}[h]
\centering
      \includegraphics[width=0.5\textwidth, trim=120 150 120 120, clip=true]{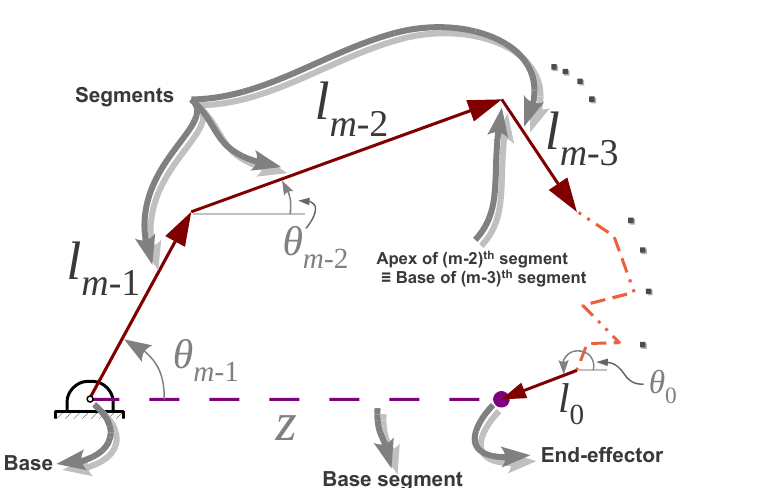}
\caption{A general planar arm at a configuration that is in its restricted configuration space, $\mathfrak{R} \subset \mathbb{T}^m$.}  \label{fig:arm-n}
\end{figure*}

The configuration space of a $n$-segmented planar robot arm is the $n$-torus, $\mathbb{T}^n$, which we coordinatize using the angles that the segments of the arm make with the positive $X$ axis (see Figure~\ref{fig:arm-n}). Thus, $[\theta_{n-1},\theta_{n-2},\cdots,\theta_0]$, with $\theta_i\in \mathbb{S}^1$, gives an unique configuration of the arm.

Throughout the paper we will assume that the end effector of the arm can be at points $q \in (\mathbb{R}^2 - \{0\})$, where $\{0\}$ is the location of the fixed base of the arm. That is, we eliminate the cases when the end effector of the robot arm coincides with the base (equivalently, the base-length, $z$, vanishes). 
Thus the configuration space that is of interest to us is $(\mathbb{T}^n - \tilde{O})$, where $\mathbb{T}^n \supset \tilde{O} = \{ [\theta_{n-1},\theta_{n-2},\cdots,\theta_0] \in \mathbb{T}^n ~\big|~ \sum_{i=0}^{n-1} l_i \sin(\theta_i) = \sum_{i=0}^{n-1} l_i \cos(\theta_i) = 0 \}$.

\subsection{The Restricted Configuration Space of a Planar Robot Arm}

The first observation that one can make is that it is sufficient to study the configuration space of the arm up to an equivalence, `$\sim$', of rotations of the entire arm about the base (see \cite{jaggi:thesis:paper:92} -- in fact the equivalence can be extended to translations and scalings of the arm as discussed in \cite{Kapovich94onthe} -- transformations that are however not relevant to a robot arm). 
In particular, this lets us design inverse kinematics assuming that the end effector will always lie on the positive $X$ axis. The corresponding set of configurations (a subspace of the complete configuration space of interest, $(\mathbb{T}^n - \tilde{O})$) will constitute the \emph{restricted configuration space}, $\mathfrak{R}$.

\begin{figure*}[h]
\centering
      \includegraphics[width=0.7\textwidth, trim=0 0 0 0, clip=true]{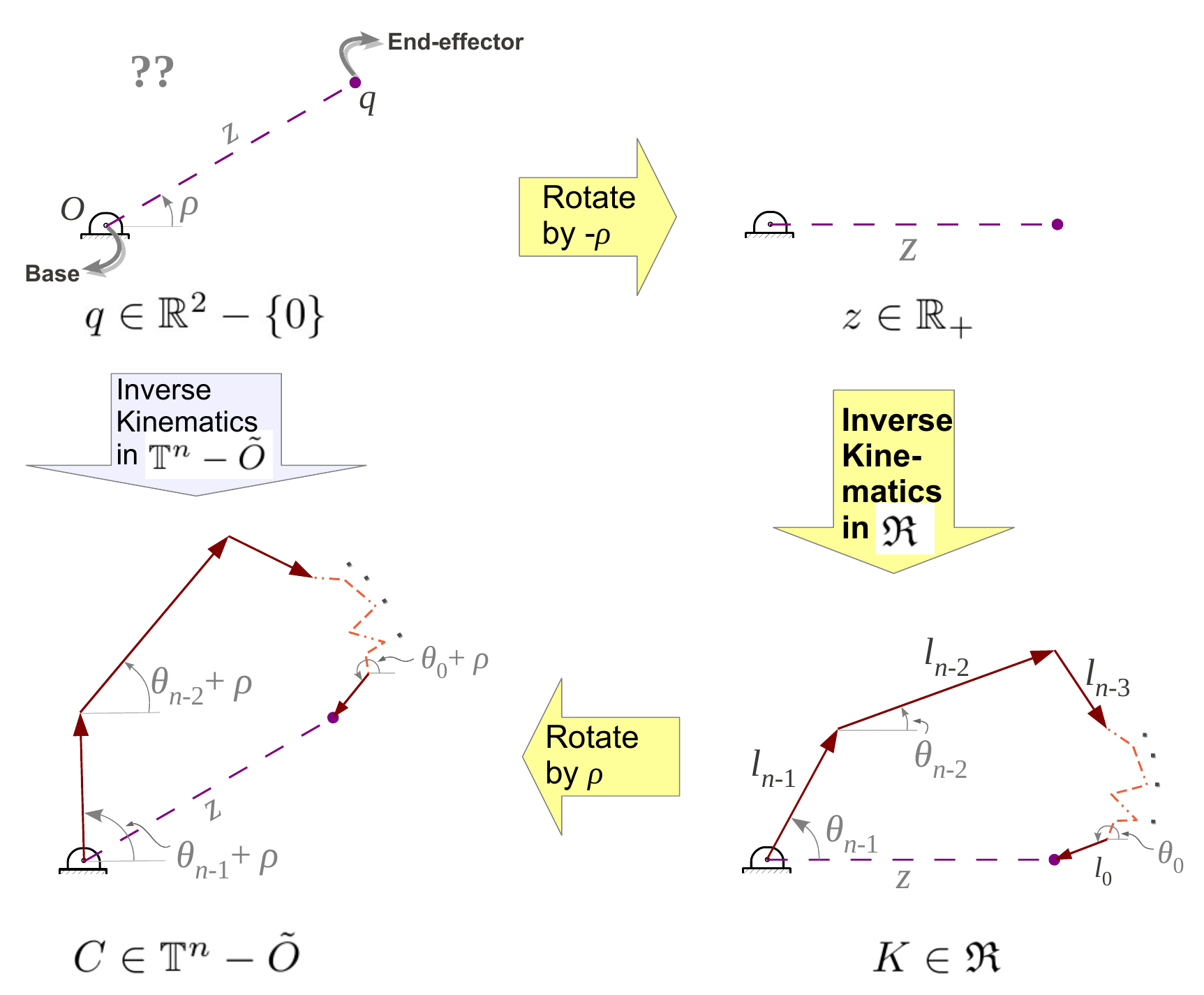}
\caption{A continuous inverse kinematics for the entire configuration space of the arm, $(\mathbb{T}^n - \tilde{O}) \rightarrow (\mathbb{R}^2 - \{0\})$, can always be constructed from a continuous inverse kinematics on the \emph{restricted} configuration space, $\mathfrak{R} \rightarrow \mathbb{R}_+$. This gives the rationale behind 
why it is sufficient to consider the restricted configuration space instead of the entire configuration space of the arm.}  \label{fig:restricted}
\end{figure*}

The sufficiency in focusing on $\mathfrak{R}$ follows from the simple observation that given any arbitrary configuration of the arm, $C \in (\mathbb{T}^n - \tilde{O})$, we can rotate the entire arm about the base by an angle, $-\rho$, ($\rho = \text{atan2}(q_y,q_x) \in \mathbb{S}^1$ being the angle that the base segment makes with the positive $X$ axis) and the resulting configuration will have the end effector on the positive $X$ axis. Furthermore, this rotation, $-\rho$, depends only on the end effector position and not on the entire configuration, $C$. Thus, given only an end effector position, $q = [q_x, q_y]^T \in \mathbb{R}^2- \{0\}$, we can compute the required rotation, following which we compute a configuration $K \in \mathfrak{R}$ (using an inverse kinematics, $\mathrm{IK}: \mathbb{R}_+ \rightarrow \mathfrak{R}$) such that the end effector is on the positive $X$ axis at a distance of $z = \|q\| = \sqrt{q_x^2 + q_y^2}$ from the base, and thus rotate $K$ by an angle of $\rho$ to get a configuration, $C\in (\mathbb{T}^n - \tilde{O})$, with end effector at $q$ (see Figure~\ref{fig:restricted}).
Furthermore, if $\mathrm{IK}: \mathbb{R}_+ \rightarrow \mathfrak{R}$ is continuous/differentiable, so will be this sequence of maps involving the rotations. 

%
%
%



Thus, in all the following discussions we will only consider configurations such that the end effector lies on a reference ray (positive $X$ axis is chosen for convenience). We write $\mathfrak{R} \subset (\mathbb{T}^n - \tilde{O})$ to denote this configuration space, and call it the \emph{restricted configuration space}. 
In terms of coordinates, if we coordinatize the $\mathbb{T}^n$ by the angles that the segments (say of lengths $l_{n-1}, l_{n-2}, \cdots, l_0$) make with the positive $X$ axis, $[\theta_{n-1},\theta_{n-2},\cdots,\theta_0]$, then $\mathfrak{R}$ denotes the subspace satisfying the condition $q_y = \sum_{j=0}^{n-1} l_j \sin(\theta_j) = 0, ~q_x = \sum_{j=0}^{n-1} l_j \cos(\theta_j) > 0$.
We will use the induced coordinatization for $\mathfrak{R}$ --- the standard coordinates on its embedding space $\mathbb{T}^n$, defined by the orientations of segments.


\begin{note}
 More formally, there exists a trivial fibration~\cite{steenrod1951topology} of the configuration space, 
$(\mathbb{T}^n -\tilde{O}) = \mathbb{S}^1 \times \mathfrak{R}$ (with projection map $\mathsf{p}: (\mathbb{T}^n -\tilde{O}) \rightarrow \mathfrak{R}$, considered as a principal $\mathbb{S}^1$ bundle), 
and a trivial fibration of the end-effector space, 
$(\mathbb{R}^2 - \{0\}) = \mathbb{S}^1 \times \mathbb{R}_+$ (with projection map $\mathsf{p}': (\mathbb{R}^2 - \{0\}) \rightarrow \mathbb{R}_+$, considered as a principal $\mathbb{S}^1$ bundle). 

The end-effector map $D: \mathfrak{R} \rightarrow \mathbb{R}_+$ is a map between the bases spaces, and a pullback of the later bundle under this map gives the former one. 
Define the resulting bundle map $\widehat{D}: (\mathbb{T}^n -\tilde{O}) \rightarrow (\mathbb{R}^2 - \{0\})$.
Thus the following relation holds (by definition of pullback bundle): $\mathsf{p}' \circ \widehat{D} = D \circ \mathsf{p}$.

The bundles being trivial, $\widehat{D}$ and $D$ also relates as $\widehat{D} = \mathrm{Id}_{\mathbb{S}^1} \times D$.
%
Thus, if there is a map $IK: \mathbb{R}_+ \rightarrow \mathfrak{R}$ such that $D \circ IK = \mathrm{Id}_{\mathbb{R}_+}$, then defining $\widehat{IK} = \mathrm{Id}_{\mathbb{S}^1} \times IK :~ (\mathbb{R}^2 - \{0\}) \rightarrow (\mathbb{T}^n -\tilde{O})$, we can observe $\widehat{D} \circ \widehat{IK} = (\mathrm{Id}_{\mathbb{S}^1} \times D) \circ (\mathrm{Id}_{\mathbb{S}^1} \times IK) = \mathrm{Id}_{\mathbb{S}^1} \times \mathrm{Id}_{\mathbb{R}_+} = \mathrm{Id}_{(\mathbb{R}^2-\{0\})}$.
There are couple of immediate consequences of the above observations:
\begin{itemize}
 \item[i.] If $IK$ is of differentiability class $C^\alpha$, so is $\widehat{IK}$.
 \item[ii.] If $z\in \mathbb{R}_+$ and $q \in (\mathbb{R}^2 - \{0\})$ relate as $\mathsf{p}'(q) = z$, then there is a homeomorphism between $D^{-1}(z) \subset \mathfrak{R}$ 
and $\widehat{D}^{-1} (q) \subset (\mathbb{T}^n - \tilde{O})$ given by $\mathsf{p} \big|_{\widehat{D}^{-1} (q)}: \widehat{D}^{-1} (q) \rightarrow D^{-1}(z)$. Thus the topology of $D^{-1}(z)$ is same as the topology of $\widehat{D}^{-1} (q)$ whenever $\mathsf{p}'(q) = z$.
\end{itemize}
Thus it will be sufficient to study the map $D$ and construct a map as $IK$ instead of studying the `hatted' versions of the same. 
We will refrain from more detailed discussion or formal proofs of the above statements in order to remain focused on the main topic of the paper.
\end{note}

\subsection{Number of Connected Components of Configuration Space for Fixed Base-Length}

We are given a robot arm, $R$, with $m$ segments and segment lengths $l_{m-1}, l_{m-2}, \cdots, l_0$ (respectively starting from the base -- Figure~\ref{fig:arm-n}). 
Suppose the length of the base segment (\emph{i.e.}, the line joining the origin and the end-effector) is $z$. 
We denote the map from the space of restricted configurations to the distance between the origin and the end effector as $D_{R}: \mathfrak{R} \rightarrow \mathbb{R}_{+}, ~ [\theta_{m-1}, \theta_{m-2}, \cdots, \theta_{0}] \mapsto z$. Due to Lemma 2 of \cite{jaggi:thesis:paper:92}, this is a Morse function. 

Now, consider the set of lengths of the sides of the closed polygon formed by the the segments of the arm and the base segment: $\{z, ~l_{m-1}, l_{m-2}, \cdots, l_0\}$ (thus, it is a polygon with $m+1$ sides). The elements of this set, when ordered according to their value, are written as $l_\ordind{m} \geq l_\ordind{m-1} \geq l_\ordind{m-2} \geq \cdots \geq l_\ordind{1} \geq l_\ordind{0}$. Thus, $l_\ordind{m} = \max \{z, ~l_{m-1}, l_{m-2}, \cdots, l_0\}$, and so on.
The following result follows from Theorem 6 of \cite{jaggi:thesis:paper:92} and Theorem 1 of \cite{Kapovich94onthe}:

\begin{theorem}[\cite{jaggi:thesis:paper:92,Kapovich94onthe}] \label{theorem:connectivity-condition}
 Keeping the end-effector of the arm fixed at a distance $z$ from the origin, and defining $l_\ordind{i}$ as above,
 \begin{itemize}
  \item[i] The constrained configuration space, $D_{R}^{-1} (z)$, is connected iff \[ l_\ordind{m} ~+~ \sum_{j=0}^{m-3} l_\ordind{j} ~~\geq~~ l_\ordind{m-1} ~+~ l_\ordind{m-2} \]
  \item[ii] The constrained configuration space, $D_{R}^{-1} (z)$, is not connected (and is a disjoint union of two $(m-2)$-dimensional tori) iff \[ l_\ordind{m} ~+~ \sum_{j=0}^{m-3} l_\ordind{j} ~~<~~ l_\ordind{m-1} ~+~ l_\ordind{m-2} \]
 \end{itemize}
\end{theorem}

\begin{figure*}[h]
\centering
      \includegraphics[width=0.5\textwidth, trim=40 50 40 30, clip=true]{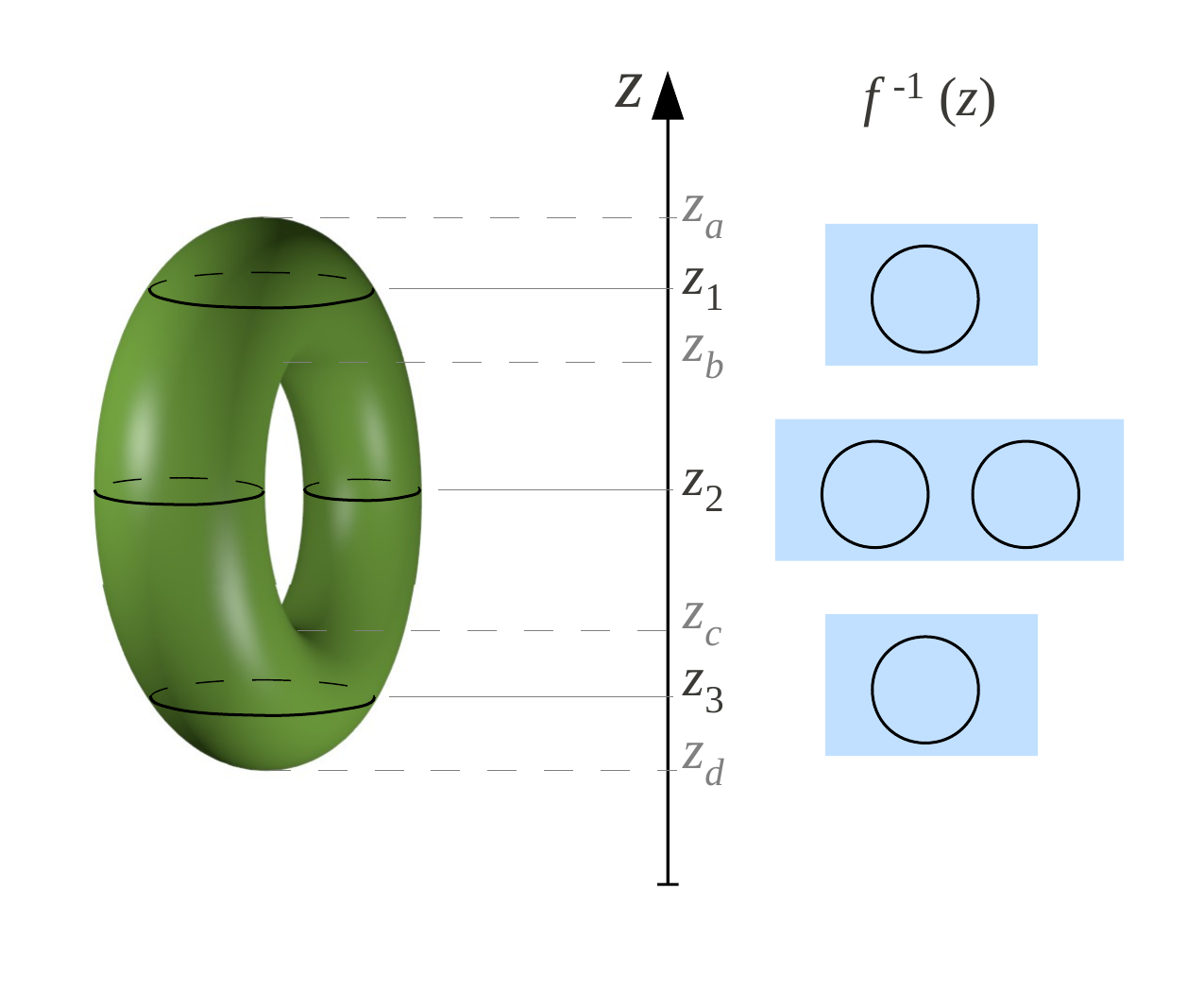}
\caption{The height function on a $2$-torus. $f^{-1}(z)$ has a single connected component when $z \in [z_d,z_c]\cup[z_b,z_a]$, and two connected components when $z\in (z_c,z_b)$.}  \label{fig:morse}
\end{figure*}

We will use a simple analogy for interpreting the above result, and hence motivate our forthcoming problem objectives. For simplicity (and for the sake of visualization), we consider the $2$-torus, $\mathbb{T}^2$, (this will be the low-dimensional analog of the restricted configuration space, $\mathfrak{R}$) with the height function, $f:\mathbb{T}^2 \rightarrow \mathbb{R}_{+}$, as the Morse function (this will be the analog of $D_{R}$). 
Given a $z\in \mathbb{R}_{+}$ (and suppose $z$ is not a critical point), as one can observe from Figure~\ref{fig:morse}, $f^{-1}(z)$ can either be a single circle or a disjoint union of two circles. The result of Theorem~\ref{theorem:connectivity-condition} essentially makes a similar assertion, but on the much larger space, $\mathfrak{R}$, and the more complicated Morse function $D_{R}$. Not only that, it also gives the range of $z$ (in form of inequalities) on which the pre-images of the Morse function are connected ($[z_d,z_c]\cup[z_b,z_a]$ in the analogous picture of Figure~\ref{fig:morse}) or disconnected ($(z_c,z_b)$ in the analogous picture of Figure~\ref{fig:morse}).

\begin{definition}[Vital Critical Points and Vital Critical Values] \label{def:vital-critical}
 Suppose $z_c$ is a critical value of the Morse function, $D_{R}$. We say $z_c$ is a \emph{vital critical value} if for an $\epsilon \rightarrow 0_+$ the number of connected components of $D_{R}^{-1}(z_c - \epsilon)$ or $D_{R}^{-1}(z_c + \epsilon)$ 
 is different from the number of connected components of $D_{R}^{-1}(z_c)$ (and neither of them is zero).
 Thus, due to Theorem~\ref{theorem:connectivity-condition}, for $z_c$ being a vital critical value, it must satisfy
 \[ l_\ordind{m} ~+~ \sum_{j=0}^{m-3} l_\ordind{j} ~~=~~ l_\ordind{m-1} ~+~ l_\ordind{m-2} \]
where, $l_\ordind{m} \geq l_\ordind{m-1} \geq l_\ordind{m-2} \geq \cdots \geq l_\ordind{1} \geq l_\ordind{0}$ are the elements of the set $\{z_c, l_{n-1}, l_{n-2},\cdots, l_0\}$ after rearrangement in order of their numeric values. 
The corresponding critical point in $\mathfrak{R}$ is a \emph{vital critical point}.
\end{definition}

\subsection{Equivalence of Configurations under Reordering} \label{sec:reorder-equivalence}

 Suppose the lengths of the segments of a robot arm, $R$, in sequence, starting from the base, are $r_{m-1}, r_{m-2}, \cdots, r_0$. 
Say its restricted configuration space is $\mathfrak{R}$, which is coordinatized by the standard coordinates (orientations of individual segments) in its embedding space, $\mathbb{T}^m$.
 
 There is the obvious base-length map, $D_{R} : [\theta_{m-1}, \theta_{m-2}, \cdots, \theta_0] (\in \mathfrak{R}) \mapsto z \in \mathbb{R}_+$, that gives a base-length for a given set of angles from the restricted configuration space, which is a Morse function as discussed earlier.
%
 An inverse kinematics, $\mathrm{IK}: \mathbb{R}_+ \rightarrow \mathfrak{R}, ~ z \mapsto [\theta_{m-1}, \theta_{m-2}, \cdots, \theta_0]$ (where $z$ is a given base-length), is an `inverse' of the base-length function
in the sense 
that $D_{R} \circ \mathrm{IK} = \text{Id}$ (the identity map on $\mathbb{R}_+$).
 
\vspace{0.05in}
 In this section we will review the properties that remain invariant upon permutation of the segments of the arm. That is, we will consider a arm, $R'$, with segment lengths $r_{\sigma(m-1)}, r_{\sigma(m-2)}, \cdots, r_{\sigma(0)}$ (respectively, starting from the base), where $\sigma$ is a permutation of the ordered set $[m-1, m-2, \cdots, 0]$.
 Let the restricted configuration space of this arm be $\mathfrak{R}'$, and the corresponding base-length map $D_{R'}$.

\begin{lemma}
 Using the usual coordinates for $\mathfrak{R}$ and $\mathfrak{R}'$ due to their respective embeddings in copies of $\mathbb{T}^m$, coordinatized by their respective segment orientations, if $[\theta_{m-1}, \theta_{m-2}, \cdots, \theta_0] \in \mathfrak{R}$, then $[ \theta_{\sigma(m-1)}, \theta_{\sigma(m-2)}, \cdots, \theta_{\sigma(0)} ] \in \mathfrak{R}'$. 
 
 Furthermore, this permutation map, $\sigma: \mathfrak{R} \rightarrow \mathfrak{R}', ~[ \theta_{m-1}, \theta_{m-2}, \cdots, \theta_0 ] \mapsto [ \theta_{\sigma(m-1)}, \theta_{\sigma(m-2)}, \cdots, \theta_{\sigma(0)} ]$, is a diffeomorphism and preserves baselength (\emph{i.e.} $D_{R}([\theta_*]) = D_{R'}\circ \sigma([\theta_*])$).
\end{lemma}
\begin{quoteproof}
 We observe that the end effector position, $[x_e, y_e]$, can be written as
{\small \begin{eqnarray*}
 x_e & = & \sum_{j=0}^{m-1} r_j \cos(\theta_j) ~=~ \sum_{j=0}^{m-1} r_{\sigma(j)} \cos(\theta_{\sigma(j)}) \\
 y_e & = & \sum_{j=0}^{m-1} r_j \sin(\theta_j) ~=~ \sum_{j=0}^{m-1} r_{\sigma(j)} \sin(\theta_{\sigma(j)}) 
\end{eqnarray*} }
This immediately implies that if $[\theta_*]\in\mathfrak{R}$, then $[\theta'_*] = \sigma([\theta_*])\in\mathfrak{R}'$, and that the base-length is preserved. 
That the permutation map is a diffeomorphism is trivial.
\end{quoteproof}

\begin{figure*}[h]
\centering
      \includegraphics[width=0.75\textwidth, trim=0 0 0 0, clip=true]{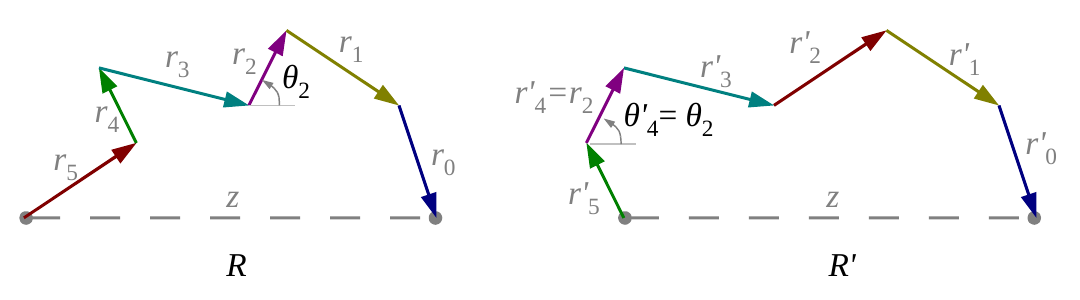}
\caption{Two arms $R$ and $R'$, with the segments permuted. Any computations on one of the arms (say, angles computed using an inverse kinematics) can be used, via a permutation, for the other arm.}  \label{fig:permute}
\end{figure*}
%

The interesting interpretation of the above lemma is that although $\mathfrak{R}$ and $\mathfrak{R}'$ are the restricted configuration spaces of two different arms, as far as forward and inverse kinematics are concerned, it does not matter which arm we use for computation --- results from one can be used for the other, with the angles only permuted (illustrated in Figure~\ref{fig:permute}). 
The following corollaries are just reiteration of this fact.

\begin{cor} 
 Suppose we are given an inverse kinematics, $IK: \mathbb{R}_+ \rightarrow \mathfrak{R}, ~z \mapsto [\theta_{m-1}, \theta_{m-2}, \cdots, \theta_0]$. 
 Then for a arm, $R'$, with segment lengths $r_{\sigma(m-1)}, r_{\sigma(m-2)}, \cdots, r_{\sigma(0)}$, where $\sigma$ is a permutation of the set ordered $[m-1, m-2, \cdots, 0]$, the map $IK' = \sigma \circ IK$ is a valid inverse kinematics that gives configurations in the restricted configuration space of $R'$.
 Furthermore, $IK'$ will have the same continuity, differentiability and other smoothness properties as $IK$.
\end{cor}

\begin{cor}
 If $IK_1, IK_2: \mathbb{R}_+ \rightarrow \mathfrak{R}$ are two inverse kinematics such that for a given $z\in\mathbb{R}_+$, $IK_1(z)$ and $IK_2(z)$ are configurations in different connected components of $D_{R}^{-1}(z) \subset \mathfrak{R}$, then $\sigma \circ IK_1(z)$ and $\sigma \circ IK_2(z)$ are also configurations in different connected components of ${D}_{R'}^{-1}(z) \subset \mathfrak{R}'$.
\end{cor}


\subsection{Identifying that Configurations are in Disconnected Components}

\begin{prop} \label{prop:disconnected-configs}
 Suppose we are given a $(m+1)$-sided polygon ($m$-segmented robot arm, with a fixed base-length) with side lengths $l_{m}, l_{m-1}, l_{m-2}, \cdots, l_{0}$.
 As before, we sort this set of $(m+1)$ values to obtain $l_\ordind{m} \geq l_\ordind{m-1} \geq l_\ordind{m-2} \geq \cdots \geq l_\ordind{1} \geq l_\ordind{0}$.
 Assume $l_\ordind{m} + \sum_{j=0}^{m-3} l_\ordind{j} < l_\ordind{m-1} + l_\ordind{m-2}$ (i.e., there are two disconnected components of the configuration space for this polygon due to Theorem~\ref{theorem:connectivity-condition}).
 For feasibility of the polygon, we also have $|l_\ordind{m-1} - l_\ordind{m-2}| \leq l_\ordind{m} + \sum_{j=0}^{m-3} l_\ordind{j}$ (which is obvious in this case. See Assumption 1 of \cite{jaggi:thesis:paper:92} for a more general conditions for feasibility).

 Consider two configurations, $K$ and $K'$, of the polygon. If in configuration $K$ the angle between the segment of length $l_\ordind{m-1}$ and the segment of length $l_\ordind{m-2}$ is $\psi \in (0,\pi) \subset \mathbb{S}^1$, and if the angle between the same segments in configuration $K'$ is $\psi' \in (-\pi,0) \subset \mathbb{S}^1$, then $K$ and $K'$ are in different connected components of the configuration space of the polygon.
\end{prop}
\begin{quoteproof}
Due to the discussion in Section~\ref{sec:reorder-equivalence}, it is sufficient to consider a polygon with ordered segment lengths, $l_\ordind{m}, l_\ordind{m-1}, l_\ordind{m-2}, \cdots, l_\ordind{0}$, and look at the angle between the segments of length $l_\ordind{m-1}$ and $l_\ordind{m-2}$.

 We start with an example (a symmetric one, for easier understanding) showing two such configurations, $K$ and $K'$, as illustrated in the figure below
\begin{center}
  \includegraphics[width=0.6\textwidth, trim=50 120 50 150, clip=true]{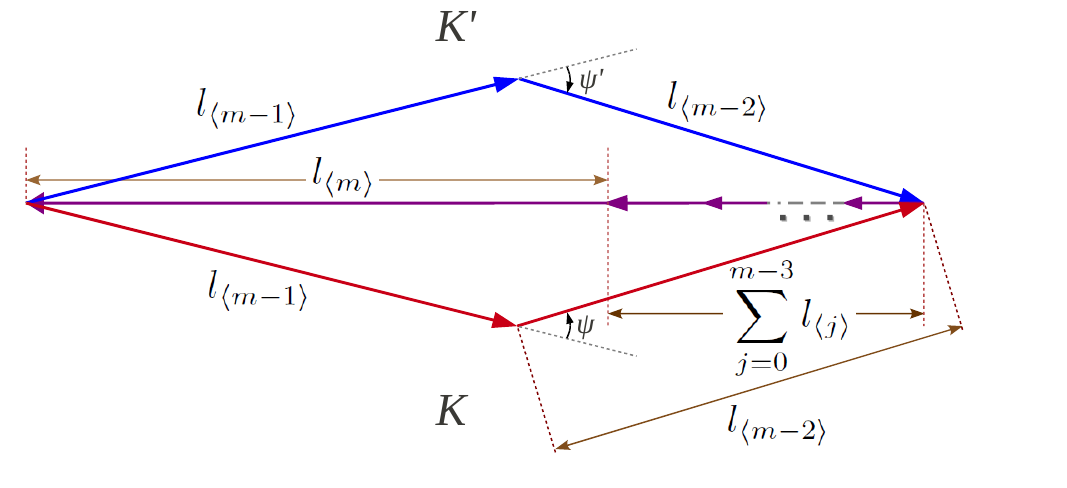}
\end{center}
The proof will consist of geometrically showing the impossibility of constructing a continuous path from one configuration to the other. Without loss of generality we keep the segment of length $l_\ordind{m}$ fixed in its position along the positive $X$ axis as shown (any configuration that is not so, can be rotated and translated as a whole to bring segment of length $l_\ordind{m}$ to that configuration).

Since $\psi \in (0,\pi) \subset \mathbb{S}^1$ and $\psi' \in (-\pi,0) \subset \mathbb{S}^1$, if we have a continuous path from one configuration to the other, 
there should at least be one point on that path (a configuration) such that 
the angle between the segments of length $l_\ordind{m-1}$ and $l_\ordind{m-2}$ 
is either $0$ or $\pi$.

\begin{center}
  \includegraphics[width=0.6\textwidth, trim=0 50 50 50, clip=true]{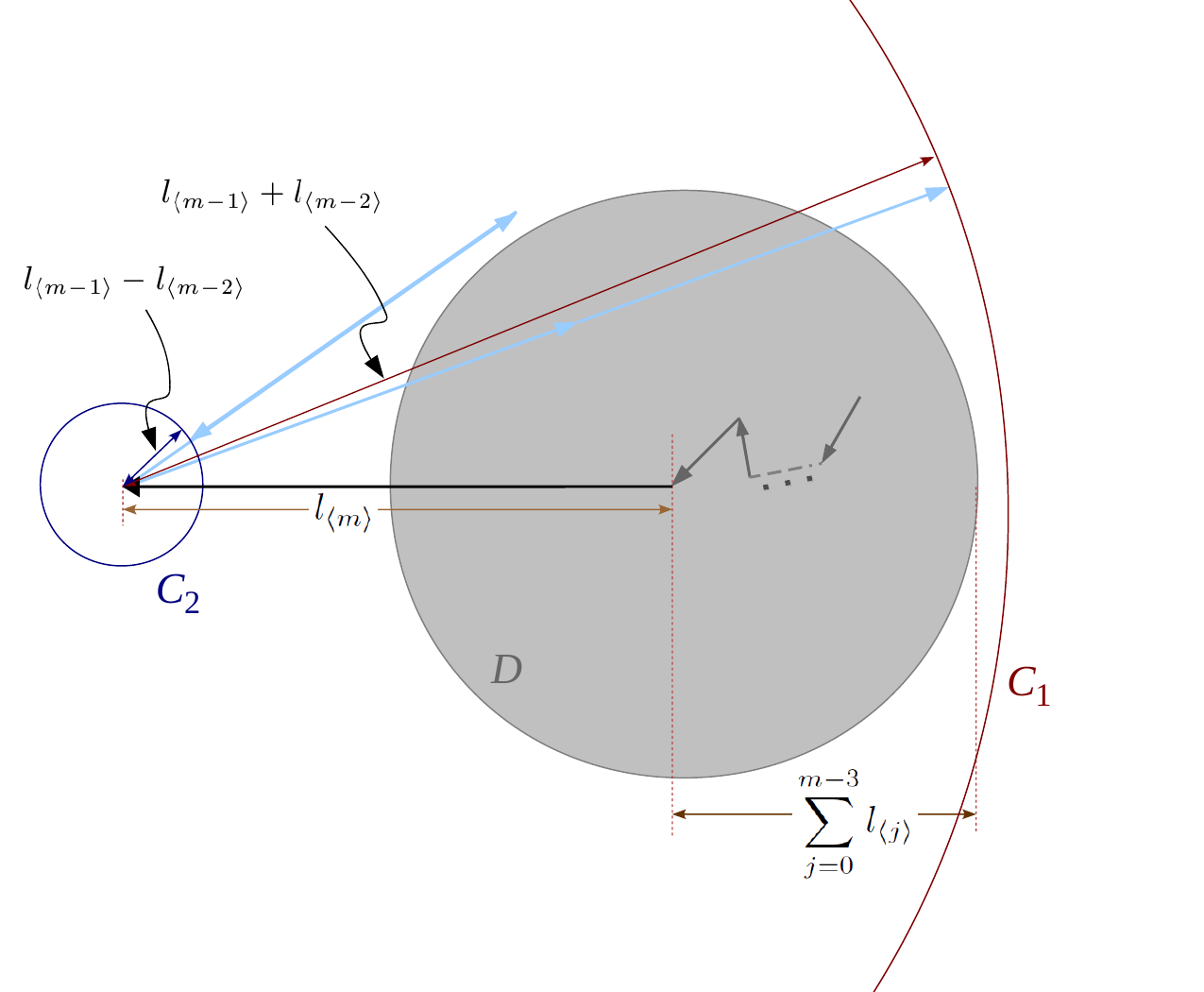}
\end{center}
However, in such a configuration the segments of length $l_\ordind{m-1}$ and $l_\ordind{m-2}$ will line up, and thus the apex of the segment of length $l_\ordind{m-2}$ will lie on a circle of radius $l_\ordind{m-1} + l_\ordind{m-2}$ or $l_\ordind{m-1} - l_\ordind{m-2}$ centered at the apex of the segment of length $l_\ordind{m}$ (circles $C_1$ and $C_2$ respectively in the above figure). 
On the other hand, the base of the segment of length $l_\ordind{m-3}$ (which should coincide with apex of the segment of length $l_\ordind{m-2}$) can stay within a disk, $D$, of radius $\sum_{j=0}^{m-3} l_\ordind{j}$ centered at the base of the segment of length $l_\ordind{m}$. 
But clearly, \\
$l_\ordind{m} + \sum_{j=0}^{m-3} l_\ordind{j} ~<~ l_\ordind{m-1} + l_\ordind{m-2} ~~\Rightarrow~~ C_1 \cap D = \emptyset$,\\
and, \\
$\left\{ \begin{array}{l}
 l_\ordind{m} + \sum_{j=0}^{m-3} l_\ordind{j} < l_\ordind{m-1} + l_\ordind{m-2} \\
 \text{and, }~~ l_\ordind{m} ~\geq~ l_\ordind{m-1}
  \end{array} \right.
$
$\Rightarrow~ 
\left\{ \begin{array}{l}
 \sum_{j=0}^{m-3} l_\ordind{j} < l_\ordind{m-2} \\
 \text{and, }~~ l_\ordind{m} ~\geq~ l_\ordind{m-1}
  \end{array} \right.
$ 
$\Rightarrow l_\ordind{m} - \sum_{j=0}^{m-3} l_\ordind{j} > l_\ordind{m-1} - l_\ordind{m-2}$
$~~\Rightarrow~~~ C_2 \cap D = \emptyset$.\\
Thus, such a configuration in which the angle between $l_\ordind{m-1}$ and $l_\ordind{m-2}$ is $0$ or $\pi$ is impossible to attain. Hence the said path between the configurations cannot exist.
\end{quoteproof}



\section{A General Continuous Inverse Kinematics Algorithm} \label{IK-alg}


We are given a robot arm, $R$, (or a part of the arm consisting only of some segments) with segment lengths $l_{m-1}, l_{m-2}, \cdots, l_0$ (not necessarily ordered). In this section we will describe a general algorithm (of computational complexity $O(m^2)$) for computing a continuous inverse kinematics 
$\mathsf{IK}: \mathbb{R}_{+} \rightarrow \mathfrak{R}$.
 Once again we use the coordinatization of $\mathbb{T}^m$ to describe points on $\mathfrak{R}$, so that $\mathsf{IK}: z \mapsto [\theta_{m-1}, \theta_{m-2}, \cdots, \theta_0] \in \mathfrak{R} \subset (\mathbb{T}^m - \tilde{O})$.

\begin{figure*}[h]
\centering
\subfigure[A general robot arm configuration.]{
      \label{fig:arm-recursive}
      \includegraphics[width=0.4\textwidth, trim=160 150 155 140, clip=true]{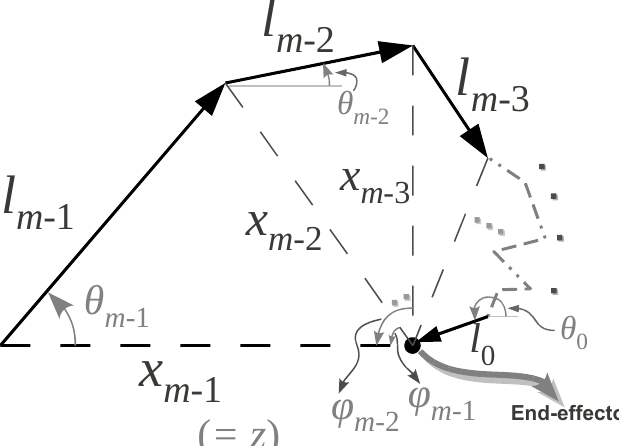}
}
\subfigure[The functions $\Theta^{+/-}_{l_p},\Phi^{+/-}_{l_p}: \mathbb{R}_{+}^2 \rightarrow \mathbb{S}^1$]{
      \label{fig:a-triangle}
      \includegraphics[width=0.36\textwidth, trim=160 170 160 140, clip=true]{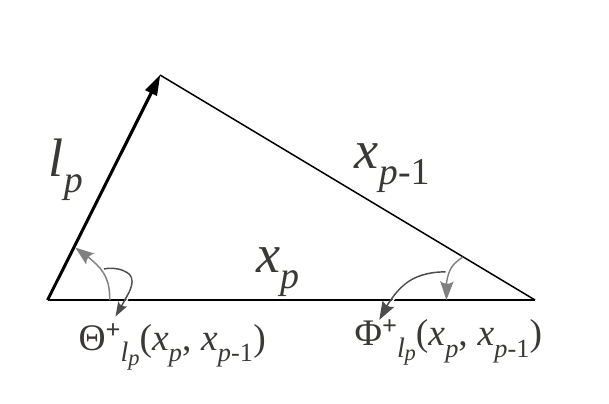}
      \includegraphics[width=0.36\textwidth, trim=160 170 160 140, clip=true]{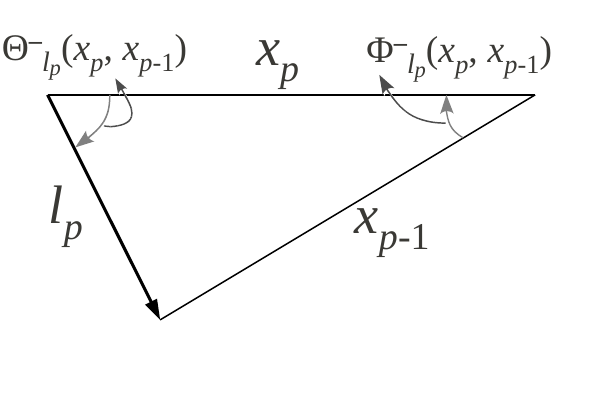}
}
\caption{.}  \label{fig:arm-recursive-algorithm}
\end{figure*}

Consider the arm configuration described in Figure~\ref{fig:arm-recursive}.
The configuration is completely described by the orientations, $\theta_p\in\mathbb{S}^1,~p=0,1,\cdots,m-1$, that the $p^{th}$ segment makes with the positive $X$ axis.
As indicated in the figure, $x_p$ is the length of the line segment connecting the base of the $p^{th}$ segment with the end effector.

Figure~\ref{fig:a-triangle} shows the two possible triangles subtended by the $p^{th}$ segment at the end effector such that the two other sides of the triangle are $x_p$ and $x_{p-1}$ (the order being consistent with the direction in which the segment points). The angles $\Theta^{+}_{l_p}(x_p, x_{p-1})$ and $\Phi^{+}_{l_p}(x_p, x_{p-1})$ are simple trigonometric functions of $l_p, x_p, x_{p-1} \in \mathbb{R}_{+}$ and return values in $[0,\pi] \subset \mathbb{S}^1$, whereas $\Theta^{-}_{\bullet}(\bullet,\bullet) = -\Theta^{+}_{\bullet}(\bullet,\bullet)$ and $\Phi^{-}_{\bullet}(\bullet,\bullet) = -\Phi^{+}_{\bullet}(\bullet,\bullet)$ return values in $\{\pi\}\cup(-\pi,0] \subset \mathbb{S}^1$.
All these functions are continuous, except when any side of the triangle is of length zero. By hypothesis $l_p>0$. Hence, in order to ensure that we don't encounter singularities or discontinuities in designing the inverse kinematics, we need to steer clear of $x_p=0$ or $x_{p-1}=0$.

\begin{note} \label{note:plus-minus-equal}~
 \begin{itemize}
  \item[i.]  If $x_p = | x_{p-1} - l_p |$, we have $\Theta^{+}_{l_p}(x_p, x_{p-1}) = \Theta^{-}_{l_p}(x_p, x_{p-1}) = \pi \in \mathbb{S}^1$ and $\Phi^{+}_{l_p}(x_p, x_{p-1}) = \Phi^{-}_{l_p}(x_p, x_{p-1}) = 0$.
  \item[ii.] If $x_p = x_{p-1} + l_p$, we have $\Theta^{+}_{l_p}(x_p, x_{p-1}) = \Theta^{-}_{l_p}(x_p, x_{p-1}) = 0$ and $\Phi^{+}_{l_p}(x_p, x_{p-1}) = \Phi^{-}_{l_p}(x_p, x_{p-1}) = 0$.
  \item[iii.] Consistent with the definition of $x_*$, we have $x_0 = l_0$. For the triangle subtended by segment of length $l_0$, clearly the angle between the segment and the line joining its base with the end effector (of length $x_0$) is $0$. Thus, for convenience we define $\Theta^{\pm}_{l_0}(x_0, x_{-1}) = 0$ (although $x_{-1}$ is not defined).
 \end{itemize}
\end{note}
\begin{note} \label{note:arm-angles}
Given a feasible set of values for $x_{m-1}, x_{m-2}, \cdots, x_1, x_0 = l_0$, and a choice of $m-1$ signs, $s_{m-1}, s_{m-2}, \cdots, s_1$ (such that $s_j$ is either $\mathquote{+}$ or $\mathquote{-}$),
\begin{itemize}
  \item[i.] A valid configuration of the arm such that the base of the segment of length $l_p$ is a distance of $x_p$ from the end effector, is described by the angles $\theta_q = \Theta^{s_q}_{l_q}(x_q, x_{q-1}) - \phi_{q+1} = \Theta^{s_q}_{l_q}(x_q, x_{q-1}) - \sum_{k=q+1}^{m-1} \Phi^{s_k}_{l_k}(x_k, x_{k-1}), ~~\forall~q=m\!-\!1,m\!-\!2,\cdots,0$ made by segment of length $l_q$ with the positive $X$ axis (note that for $q=0$, the sign $s_0$ is irrelevant since $\Theta^{\pm}_{l_0}(x_0, x_{-1}) = 0$).
  \item[ii.] The angle between the $q^{th}$ segment (of length $l_q$) and $(q-1)^{th}$ segment (of length $l_{q-1}$) is given by $\theta_q - \theta_{q-1} = \Theta^{s_q}_{l_q}(x_q, x_{q-1}) + \Phi^{s_q}_{l_q}(x_q, x_{q-1}) - \Theta^{s_{q-1}}_{l_{q-1}}(x_{q-1}, x_{q-2}), ~~q = 1,2,\cdots, m-1$
 \end{itemize}
\end{note}

%
%


\subsection{The Range Function} \label{sec:recursive-range}

We once again focus on Figure~\ref{fig:a-triangle}. If $x_{p-1}$ can assume all values in the interval $[\underline{x}_{p-1}, \overline{x}_{p-1}]$, we would like to study the range of values that $x_p$ can assume for the triangle in Figure~\ref{fig:a-triangle} to be a feasible one. 

\begin{figure*}[h]
\centering
\subfigure[Relation between $\lbrack \underline{x}_{p}, \overline{x}_{p}\rbrack$ and $\lbrack\underline{x}_{p-1}, \overline{x}_{p-1}\rbrack$, and an inverse kinematics component function, $f_{p-1}$. A \emph{sign switching set}, $S_p = \{\xi_p^1,\xi_p^2,\xi_p^3,\xi_p^4\}$, is also shown. The `$+$' and `$-$' indicate a choice of \emph{sign assignment}.]{
\label{fig:range-recursive}
      \includegraphics[width=0.5\textwidth, trim=0 0 0 30, clip=true]{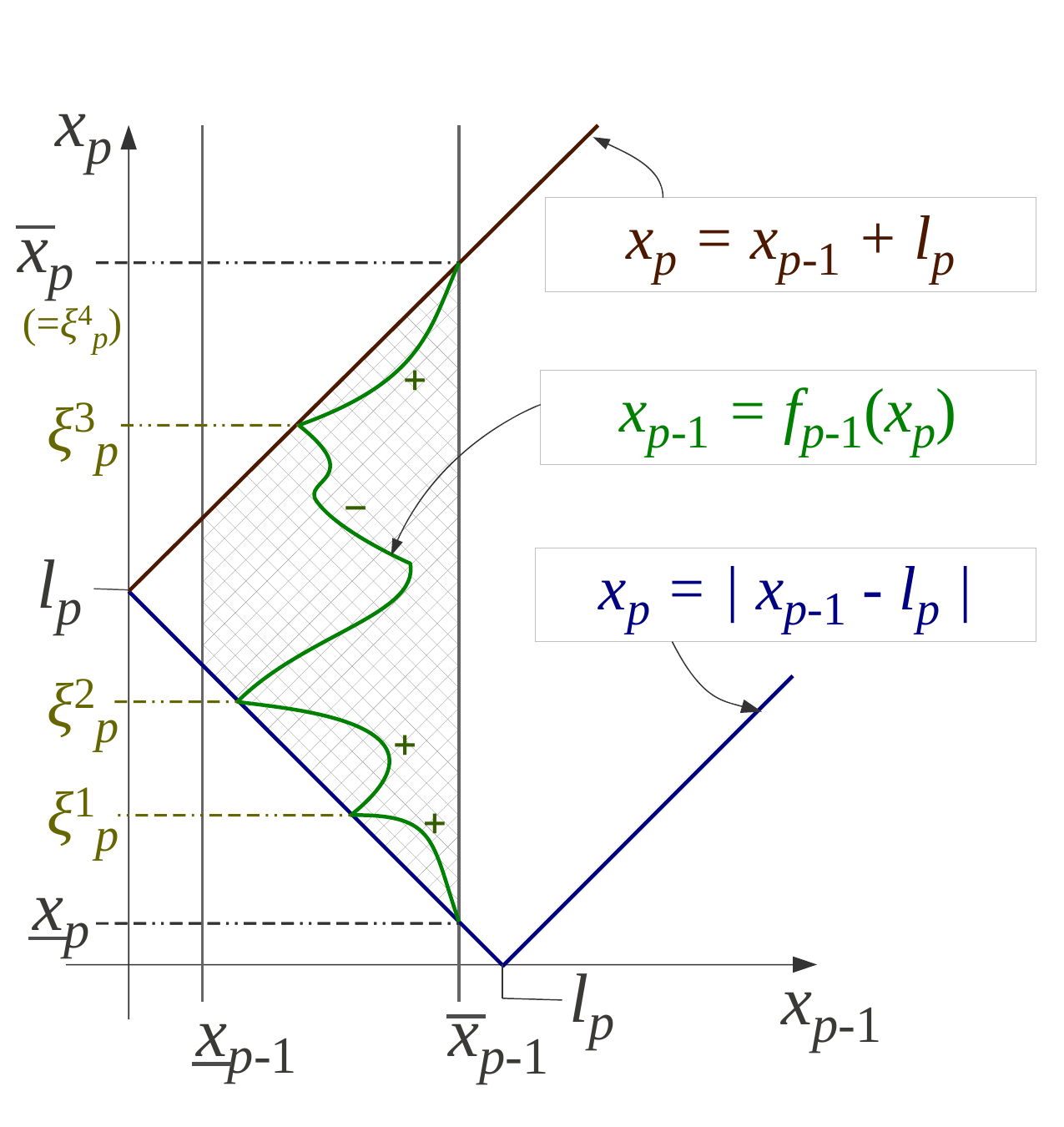}
 } \hspace{0.05in}
\subfigure[Other possible shapes for the region $|x_{p-1} - l_p| \leq x_p \leq x_{p-1} + l_p$, $x_p\in \lbrack \underline{x}_{p}, \overline{x}_{p}\rbrack$.]{
      \includegraphics[width=0.26\textwidth, trim=0 30 0 0, clip=true]{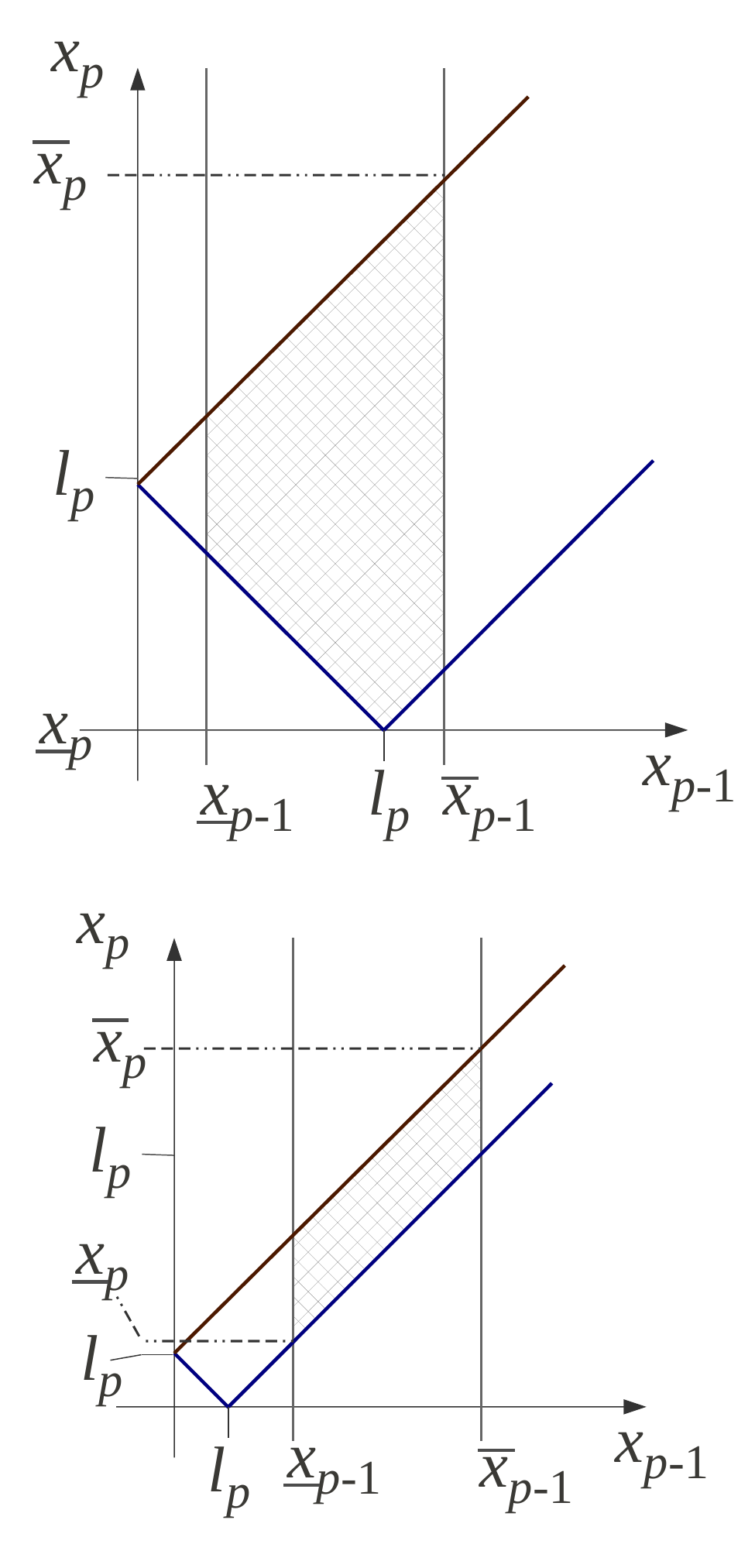}
}
\caption{The regions defined by $|x_{p-1} - l_p| \leq x_p \leq x_{p-1} + l_p$(left: $\overline{x}_p<l_p$, top right: $\underline{x}_p<l_p\leq \overline{x}_p$, and bottom right: $l_p \leq \underline{x}_p$), and an example of a IKCF-IKCSSS-IKCSA tuple.}  \label{fig:range-recursive-full}
\end{figure*}

The relation (triangle inequality) that $x_p$ and $x_{p-1}$ must satisfy is simple:
\begin{equation}
 |x_{p-1} - l_p| ~~\leq~~ x_p ~~\leq~~ x_{p-1} + l_p \label{eq:triangle-inequality}
\end{equation}
It is easy to deduce from these inequalities that $x_p$ can thus assume any value in the interval $[\underline{x}_{p}, \overline{x}_{p}]$ (see Figure~\ref{fig:range-recursive-full}) defined by,
\begin{eqnarray}
 \underline{x}_{p} ~=:~ \underline{\underline{R}}(\underline{x}_{p-1}, \overline{x}_{p-1}; l_{p}) & = & \left\{ \begin{array}{ll}
                          \underline{x}_{p-1} - l_{p}, & \text{if }~ l_{p} \leq \underline{x}_{p-1} \\
                          0, & \text{if }~ \underline{x}_{p-1} < l_{p} \leq \overline{x}_{p-1} \\
                          l_{p} - \overline{x}_{p-1}, & \text{if }~ \overline{x}_{p-1} < l_{p}
                         \end{array} \right. \nonumber \\ 
 \overline{x}_{p} ~=:~ \overline{\overline{R}}(\underline{x}_{p-1}, \overline{x}_{p-1}; l_{p}) & = & l_{p} ~+~ \overline{x}_{p-1} \label{eq:range-recursive}
\end{eqnarray}
Clearly, $\underline{x}_{0} = \overline{x}_{0} = l_0$. Thus, using the recursive relation of \eqref{eq:range-recursive} we can work out $\underline{x}_{p}$ and $\overline{x}_{p}$ for all $p=0,1,2,\cdots,m-1$ (the computational complexity being linear in $m$).

Moreover, since $\overline{x}_p$ or $\underline{x}_p$ depend only on the values in the set $\{l_0, l_1, \cdots, l_p\}$, we re-write the above recursive relations using the following simplified notations:
\begin{eqnarray}
 \underline{x}_{p} = \underline{R}(\{l_0, l_1, \cdots, l_p\}) \nonumber \\
 \overline{x}_{p} = \overline{R}(\{l_0, l_1, \cdots, l_p\})  \label{eq:range}
\end{eqnarray}

\imphead{Interpretation} $\underline{R}(\{l_0, l_1, \cdots, l_p\})$ is the minimum possible value of the base length of an arm with segment lengths $l_0, l_1, \cdots, l_p$, and $\overline{R}(\{l_0, l_1, \cdots, l_p\})$ is the maximum possible value.

\begin{prop} [The closed-form expressions for $\underline{R}$ and $\overline{R}$]~\label{prop:range-closed}
\begin{itemize}
 \item[1.] $\overline{R}(\{l_0, l_1, \cdots, l_p\}) = \sum_{j=0}^p l_j$
 \item[2.] Since, due to the discussion of Section~\ref{sec:reorder-equivalence}, the order of the elements in the set $\{l_0, l_1, \cdots, l_p\}$ does not change the value of $\underline{R}(\{l_0, l_1, \cdots, l_p\})$, without loss of generality we assume $l_p \geq l_{p-1} \geq \cdots \geq l_1 \geq l_0$.
 Then, $\underline{R}(\{l_0, l_1, \cdots, l_p\}) = \left\{ \begin{array}{ll} 
                               l_p - \sum_{j=0}^{p-1} l_j, & \text{ if } l_p > \sum_{j=0}^{p-1} l_j, \\
                               0, & \text{ otherwise.} 
                            \end{array}\right. $\newline
\end{itemize}
\end{prop}
\begin{quoteproof}
 `1.' follows trivially from the recursive expression for $\overline{x}_{p}$ in \eqref{eq:range-recursive}.

\vspace{0.05in} We prove `2.' by considering the value of $\underline{R}(\{l_0, l_1, \cdots, l_p\})$ separately for the two cases.

\noindent \textbf{Case $l_{p} > \sum_{j=0}^{p-1} l_{j}$:} Using \eqref{eq:range-recursive},
\begin{eqnarray*}
& & \underline{R}(\{l_0, l_1, \cdots, l_p\}) \\
  & = & l_p - \overline{x}_{p-1} ~~\text{[where, $\overline{x}_{p-1} = \overline{R}(\{l_{0}, l_{1}, \cdots, l_{p-1}\})$,} \\ & & \qquad\qquad\qquad\quad \text{and since $\overline{R}(\{l_{0}, l_{1}, \cdots, l_{p-1}\}) = \sum_{j=0}^{p-1} l_{j} < l_{p}$]} \\
 & = & l_{p} - \sum_{j=0}^{p-1} l_{j}
\end{eqnarray*}

\noindent \textbf{Case $l_{p} \leq \sum_{j=0}^{p-1} l_{j}$}:
This condition implies $l_{p} \leq  \overline{R}(\{l_0, l_1, \cdots, l_p\}) =: \overline{x}_{p-1}$. 
Furthermore, $\underline{x}_{p-1}$ is the minimum value of the base of an arm with segment lengths $\{l_{0}, l_{1}, \cdots, l_{p-1}\}$. It is easy to check that since the values are ordered, we have $0 \leq l_{p-1} \!-\! l_{p-2} \!+\! l_{p-3} \!-\! l_{p-4} \!+\! \cdots l_{0} \leq l_{p-1}$, and the base-length of $(l_{p-1} \!-\! l_{p-2} \!+\! l_{p-3} \!-\!\cdots)$ is achievable by the arm with segment lengths $\{l_{0}, l_{1}, \cdots, l_{p-1}\}$ (in the configuration when all the segments are aligned along a single line, but with alternating orientations). Thus for its minimum possible value we must have $\underline{x}_{p-1} \leq (l_{p-1} \!-\! l_{p-2} \!+\! l_{p-3} \!-\!\cdots) \leq l_{p-1}$. Thus we have $\underline{x}_{p-1} \leq l_{p-1} \leq \overline{x}_{p-1}$. Using \eqref{eq:range-recursive} we thus immediately have
\begin{eqnarray*}
 \underline{R}(\{l_{0}, l_{1}, \cdots, l_{p}\}) & = & 0
\end{eqnarray*}
\end{quoteproof}

\begin{cor} \label{cor:range-min} 
 \[ \underline{R}(\{l_{0}, l_{1}, \cdots, l_{p}\}) ~~\leq~~ | \pm l_{0} \pm l_{1} \pm \cdots \pm l_{p} | \]
for any assignment of `$+$' or `$-$' for the `$\pm$' signs on the right hand side of the inequality.
\end{cor}
\begin{quoteproof}
 Once again we assume, without the loss of generality, that $l_p \geq l_{p-1} \geq \cdots \geq l_1 \geq l_0$.
 The R.H.S. is always non-negative. The L.H.S. is either $0$ or  $l_{p} - \sum_{j=0}^{p-1} l_{j}$ (when $l_{p} > \sum_{j=0}^{p-1} l_{j}$). If L.H.S. is $0$ there is nothing to prove. However, if $l_{p} > \sum_{j=0}^{p-1} l_{j}$, we note that the R.H.S. can be re-written as $l_{p} \pm l_{p-1} \pm l_{p-2} \pm \cdots \pm \l_{0}$, for any assignment of `$+$' or `$-$' for the `$\pm$'s. This clearly is greater than or equal to $l_{p} - l_{p-1} - l_{p-2} - \cdots - \l_{0} \geq \underline{R}(\{l_{0}, l_{1}, \cdots, l_{p}\})$.
\end{quoteproof}

\subsection{Inverse Kinematics Components}

\begin{definition} [IKCF]
 Given $\underline{x}_{p}, \overline{x}_{p}, \underline{x}_{p-1}$ and $\overline{x}_{p-1}$ as described above, we define a $(p-1)^{th}$ ``inverse kinematics component function'' (IKCF) as a continuous function $f_{p-1}: [\underline{x}_{p}, \overline{x}_{p}] \rightarrow [\underline{x}_{p-1}, \overline{x}_{p-1}]$ such that every point on the graph of $f_{p-1}$ satisfies inequalities \eqref{eq:triangle-inequality} (\emph{i.e.}, the graph lies in the hatched region of Figure~\ref{fig:range-recursive-full}), and such that $f_{p-1}(x_p) > 0, ~\forall x_p \in [\underline{x}_{p}, \overline{x}_{p}]$.

 It is always possible to construct such a IKCF if $l_p, \overline{x}_{p} >0$, since then the region defined by inequalities \eqref{eq:triangle-inequality} is a non-empty convex polygon (Figure~\ref{fig:range-recursive-full}). 
 Such a function returns a feasible value of $x_{p-1}$ given a value of $x_p$ in its feasible range. Clearly, $f_0: [\underline{x}_{1}, \overline{x}_{1}] \rightarrow \{l_0\}$ is a constant function (since $\underline{x}_{0} = \overline{x}_{0} = l_0$).
%
\end{definition}

\imphead{Notation}
If $f_{p-1}(x_{p}) = x_{p-1}$, to denote the corresponding point on the graph of $f_{p-1}$ we will write $(x_{p-1}, x_{p})$ (instead of the standard $(x_{p}, x_{p-1})$).
\vspace{0.07in}

\begin{definition} [IKCSSS]
 For a given IKCF, $f_{p-1}: [\underline{x}_{p}, \overline{x}_{p}] \rightarrow [\underline{x}_{p-1}, \overline{x}_{p-1}]$, let
 {\small \[ \Xi(f_{p-1}) ~=~ \big\{ x_p \in [\underline{x}_{p}, \overline{x}_{p}] ~\big| \text{ Either } x_p = | f_{p-1}(x_p) - l_p| \text{ or } x_p = f_{p-1}(x_p) + l_p\text{.} \big\} \]}
 Note that this set contains at least one point, namely, $\overline{x}_p$. 
 A finite countable subset of points, $S_p = \{ \xi_p^1, \xi_p^2, \cdots, \xi_p^{h_p-1}, \xi_p^{h_p}=\overline{x}_p\} \subseteq \Xi(f_{p-1})$, such that $\overline{x}_p$ is an element of the set, is called an ``inverse kinematics component sign switching set'' (IKCSSS).
 Without loss of generality we assume $S_p$ to be an ordered set with $\xi_p^1 \leq \xi_p^2 \leq \cdots \leq \xi_p^{h_p-1} \leq \xi_p^{h_p} = \overline{x}_p$ (see Figure~\ref{fig:range-recursive}).
\end{definition}

\vspace{0.07in}
\begin{definition} [IKCSA]
 Given a IKCF, $f_{p-1}$, and a corresponding valid choice of IKCSSS, $S_p$, we define a ``inverse kinematics component sign assignment'' (IKCSA) as a map 
$\mathbf{sg}_p: S_p \rightarrow \{\mathquote{+}, \mathquote{-}\}$. 
\end{definition}

\vspace{0.07in}
\begin{definition} [IKCISF] \label{def:ikcisf}
The reason of defining the IKCSSS and a corresponding IKCSA is that now we can assign a `sign' to the intervals $[\underline{x}_p, \xi_p^1)$, $[\xi_p^1, \xi_p^2)$, $[\xi_p^2, \xi_p^3)$, $\cdots$, $[\xi_p^{h_p-1}, \overline{x}_p]$, as follows:
\begin{eqnarray} 
 \mathscr{S}_{S_p\!, \mathbf{sg}_p} & : & [\underline{x}_p, \overline{x}_p] \rightarrow \{\mathquote{+}, \mathquote{-}\} \nonumber \\
                                           & : & x_p \mapsto \left\{ \begin{array}{ll}
												\mathbf{sg}_p(\xi_p^1), & \text{ if } \underline{x}_p \leq x_p < \xi_p^{1} \\
												\mathbf{sg}_p(\xi_p^2), & \text{ if } \xi_p^{1} \leq x_p < \xi_p^{2} \\
												\vdots & \\
												\mathbf{sg}_p(\xi_p^{h_p}), & \text{ if } \xi_p^{h_p-1} \leq x_p \leq \xi_p^{h_p} = \overline{x}_p \\
												\end{array} \right. \label{eq:interval-sign-fun}
\end{eqnarray}
This function, $\mathscr{S}_{S_p\!, \mathbf{sg}_p}$, defined using a given pair of IKCSSS and IKCSA, will be refereed to as ``inverse kinematics component interval sign function'' (IKCISF), or simply the ``interval sign function''. 
This is illustrated using the $\mathquote{+}$ or $\mathquote{-}$ in Figure~\ref{fig:range-recursive}.
\end{definition}
\noindent [\emph{Note:} If $\xi_p^{j-1} = \xi_p^j$, then $\mathbf{sg}_p(\xi_p^j)$ is essentially not used in the construction of $\mathscr{S}_{S_p\!, \mathbf{sg}_p}$.]

\vspace{0.05in}
We will write $\{f_{p-1}, \mathscr{S}_{S_p\!, \mathbf{sg}_p}\}$ to indicate a choice of IKCF, and an interval sign assignment due to the choice of a corresponding IKCSSS-IKCSA pair. 
Due to the following lemma, a choice of these determines a continuous map from the range of $x_p$ to the space of configuration for the triangle with sides $l_p, x_p$ and $f_{p-1}(x_p)$ (see Figure~\ref{fig:a-triangle}).

\vspace{0.05in}
\begin{lemma} \label{lemma:trig-fun-signed}
%
 Given a IKCF-IKCSSS-IKCSA tuple, $\{f_{p-1}, \mathscr{S}_{S_p\!, \mathbf{sg}_p}\}$, 
the following functions, $\Theta_{p}, \Phi_{p}:  [\underline{x}_{p}, \overline{x}_{p}]-\{0\} \rightarrow \mathbb{S}^1$, are continuous:
\begin{eqnarray}
 \Theta_{p}(x_p) & := & \Theta^{ \mathscr{S}_{S_p\!, \mathbf{sg}_p}(x_p) }_{l_p} (x_p, f_{p-1}(x_p) ) \nonumber \\
 \Phi_{p}(x_p) & := & \Phi^{ \mathscr{S}_{S_p\!, \mathbf{sg}_p}(x_p) }_{l_p} (x_p, f_{p-1}(x_p) )
\end{eqnarray}
%
where, by $\Theta^{s}_{l_p}$ we mean $\Theta^{+}_{l_p}$ or $\Theta^{-}_{l_p}$, depending on whether $s$ is `$+$' or `$-$' (and likewise for $\Phi$).
\end{lemma}
\begin{quoteproof}
 We observe that $f_{p-1}$ is a continuous function, and that $\Theta^{+/-}_{l_p}, \Phi^{+/-}_{l_p}: \mathbb{R}_{+}^2 \rightarrow \mathbb{S}^1$ are continuous functions in their respective domains except where either of its inputs are zero. But by the construction of IKCF, $f_{p-1}(x_p)$ does not give a value of zero. Thus for all $x_p \in [\underline{x}_p,\overline{x}_p] - \{0\}$, $\Theta^{+/-}_{l_p}(x_p, f_{p-1}(x_p) )$ and $\Phi^{+/-}_{l_p}(x_p, f_{p-1}(x_p) )$ are continuous.

 Next we observe that $\Theta^{+}_{l_p}(\xi_p^i, f_{p-1}(\xi_p^i)) = \Theta^{-}_{l_p}(\xi_p^i, f_{p-1}(\xi_p^i))$ and $\Phi^{+}_{l_p}(\xi_p^i, f_{p-1}(\xi_p^i)) = \Phi^{-}_{l_p}(\xi_p, f_{p-1}(\xi_p^i))$ for all $\xi_p^i \in S_p$. 

 Thus $\Theta_{p},\Phi_{p}: \mathbb{R}_{+} \rightarrow \mathbb{S}^1$ are made up of piece-wise continuous functions that agree at the points where the constituent functions are pieced together. This concludes the proof.
\end{quoteproof}
Note that in the above Lemma, $\Phi_p$ and $\Theta_p$ are defined for $p=1,2,\cdots,m-1$. 
We extend the definition for $p=0$ by letting $\Phi_0, \Theta_0: \{l_0\} \rightarrow 0 \in \mathbb{S}^1$ (constant function that returns the identity element of the circle group, $\mathbb{S}^1$).

%

\subsection{The Inverse Kinematics Algorithm} \label{sec:recursive-IK}

Thus, given IKCF-IKCSSS-IKCSA tuples, $\{f_{p-1}, \mathscr{S}_{S_p\!, \mathbf{sg}_p}\}$, for $p=1,2,\cdots,m-1$, 
and given a base-length, $z \equiv x_{m-1} > 0$, 
we construct a continuous inverse kinematics for the entire arm as follows:
\begin{quote}

For a given value of $z \equiv x_{m-1}$ (the base-length), we can compute
\begin{equation} \label{eq:x-k}
 x_k ~=~  F_k(z) ~:=~ f_{k}\circ f_{k+1}\circ f_{k+2}\circ\cdots \circ f_{m-2}(z) 
\end{equation}
Clearly, $F_k$ are continuous $\forall k=0,1,2,\cdots,m-2$. Moreover, by the definition of IKCF, if $z>0$, then $x_k>0$.

Thus, a complete configuration for the arm is determined by (refer to Figure~\ref{fig:arm-recursive-algorithm}) the following expression for the orientation of the $q^{th}$ segment, 
\begin{eqnarray} 
\theta_q & = & \Theta_{q}( x_{q} ) ~-~ \phi_{q+1}  \nonumber \\
         & = &  \Theta_{q}( x_{q} ) ~-~ \sum_{k=q+1}^{m-1} \Phi_{k}( x_{k} )  \nonumber \\
         & = &  \Theta^{ \mathscr{S}_{S_q\!, \mathbf{sg}_q}(x_q) }_{q}( x_{q}, f_{q-1}(x_q) ) ~-~ \sum_{k=q+1}^{m-1} \Phi^{ \mathscr{S}_{S_k\!, \mathbf{sg}_k}(x_k) }_{k}( x_{k}, f_{k-1}(x_k) ) \nonumber \\ & & \qquad\qquad\qquad\in~ \mathbb{S}^1 \label{eq:theta-q}
\end{eqnarray}
where $\Phi_k, \Theta_k$ are dependent on the tuples $\{f_{k-1}, \mathscr{S}_{S_k\!, \mathbf{sg}_k}\}$ and are as defined in Lemma~\ref{lemma:trig-fun-signed}. 
Since $x_k>0$, using Lemma~\ref{lemma:trig-fun-signed} it follows that $\theta_q$ varies continuously with the base-length, $z$.

Furthermore, the fact that $[\theta_{m-1}, \theta_{m-2}, \cdots, \theta_0 ] \in \mathfrak{R} \subset \mathbb{T}^m$ follows from our very construction (see Figure~\ref{fig:arm-recursive-algorithm}). 
This can also be proved explicitly by simplifying the trigonometric expressions $x_e = \sum_{j=0}^{m-1} l_j \cos(\theta_j), ~y_e = \sum_{j=0}^{m-1} l_j \sin(\theta_j)$ and using the trigonometric formulae for $\Theta^{+/-}$ and $\Phi^{+/-}$ in the orientations of the segments described in Note~\ref{note:arm-angles}. This we omit, considering it a simple but lengthy exercise. 
\end{quote}
%

\noindent We thus have the following proposition:
\begin{prop} \label{prop:continuous-IK}
 The map from base-length, $z>0$, to the space of arm configurations described by the segment orientations as determined by \eqref{eq:theta-q} (along with \eqref{eq:x-k}) is a continuous map from $\mathbb{R}_+$ to the restricted configuration space, $\mathfrak{R}$. 

Since the map is completely determined by the tuple $\{f_{p-1}, \mathscr{S}_{S_p\!, \mathbf{sg}_p}\}, ~p=1,2,\cdots,m-1$, for brevity we write this map as
\begin{eqnarray*}
 \mathsf{IK}_{\{f_{*-1}, \mathscr{S}_{S_*\!, \mathbf{sg}_*}\}} & : & \mathbb{R}_+ \rightarrow \mathfrak{R} \\
                                                &  &  z \mapsto [\theta_{m-1}, \theta_{m-2}, \cdots, \theta_0 ]
\end{eqnarray*}
\end{prop}



\subsection{Some Particular IKCFs}

\begin{figure*}[h]
\centering
      \includegraphics[width=0.45\textwidth, trim=40 50 40 30, clip=true]{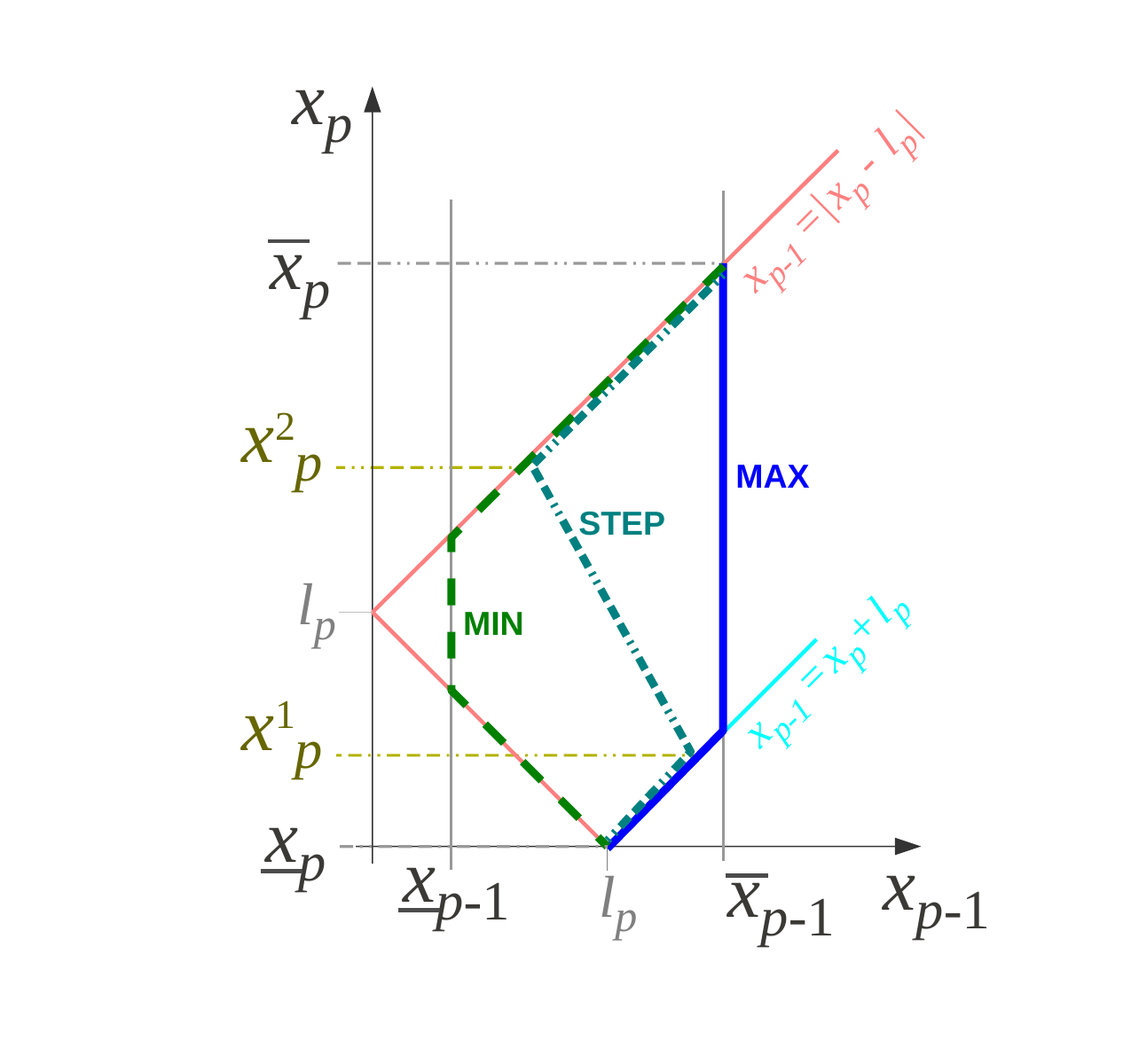}
\caption{The functions $\mathsf{MIN}$, $\mathsf{STEP}$ and $\mathsf{MAX}$ that we use to construct IKCFs.}  \label{fig:fn-minmaxstep}
\end{figure*}

In this section we will describe three types of functions that will be particularly useful for constructing IKCFs in the next sections.
First of all we are given the following parameters for defining $f_{p-1}$: The domain, $[\underline{x}_{p}, \overline{x}_{p}]$, the codomain $[\underline{x}_{p-1}, \overline{x}_{p-1}]$, and $l_p$.
Thus we define the following:
\begin{eqnarray}
 \mathsf{MIN}_{\{\underline{x}_{p-1}, l_p\}} & \!:\! & [\underline{x}_{p}, \overline{x}_{p}] \rightarrow [\underline{x}_{p-1}, \overline{x}_{p-1}], \nonumber \\
              & \!:\! & x_p ~\mapsto~ \max(\underline{x}_{p-1},~ |\underline{x}_{p-1} - l_p|) \label{eq:fun-min} \\
 & & \nonumber \\
 \mathsf{MAX}_{\{\overline{x}_{p-1}, l_p\}} & \!:\! & [\underline{x}_{p}, \overline{x}_{p}] \rightarrow [\underline{x}_{p-1}, \overline{x}_{p-1}], \nonumber \\
              & \!:\! & x_p ~\mapsto~ \min(\overline{x}_{p-1},~ \underline{x}_{p-1} + l_p) \label{eq:fun-max} 
\end{eqnarray}
Also, for given $x^1_p, x^2_p \in [\underline{x}_{p}, \overline{x}_{p}]$, with $x^1_p < x^2_p$, we define
\begin{eqnarray}
 & &\!\!\!\! \mathsf{STEP}_{\{\underline{x}_{p-1}, \overline{x}_{p-1}, l_p, x^1_p, x^2_p\}} ~:~ [\underline{x}_{p}, \overline{x}_{p}] \rightarrow [\underline{x}_{p-1}, \overline{x}_{p-1}], \nonumber \\
 & &           \qquad ~:~ x_p ~\mapsto~ \left\{ \begin{array}{ll}
                               \mathsf{MAX}_{\{\overline{x}_{p-1}, l_p\}}(x_p), &\text{if $x_{p} \leq x^1_p$} \\
                               \frac{x^2_p-x_p}{x^2_p-x^1_p} \mathsf{MAX}_{\{\overline{x}_{p-1}, l_p\}}(x_p) + \frac{x_p-x^1_p}{x^2_p-x^1_p}\mathsf{MIN}_{\{\underline{x}_{p-1}, l_p\}}(x_p) , &\text{if $x^1_p < x_{p} < x^2_p$}\\
                               \mathsf{MIN}_{\{\underline{x}_{p-1}, l_p\}}(x_p), &\text{if $x^2_p \leq x_{p} $}
                              \end{array} \right. \nonumber \\
\end{eqnarray}

Referring to Figure~\ref{fig:fn-minmaxstep}, the function $\mathsf{MIN}$ essentially returns the minimum of the possible values of $x_{p-1}$ (within the feasible region -- the \emph{hatched} region of Figure~\ref{fig:range-recursive-full}) for a given $x_p$, while $\mathsf{MAX}$ returns the maximum. $\mathsf{STEP}$ on the other hand, returns the minimum for $x_p$ greater than $x^2_p$, and maximum for $x_p$ less than $x^1_p$, while linearly interpolating in between.




\section{Sampling Exactly One Configuration from Each Connected Component of $D_{R}^{-1}(z)$} \label{sec:sampling-main}

\begin{figure*}[h]
\centering
      \includegraphics[width=0.45\textwidth, trim=40 50 40 30, clip=true]{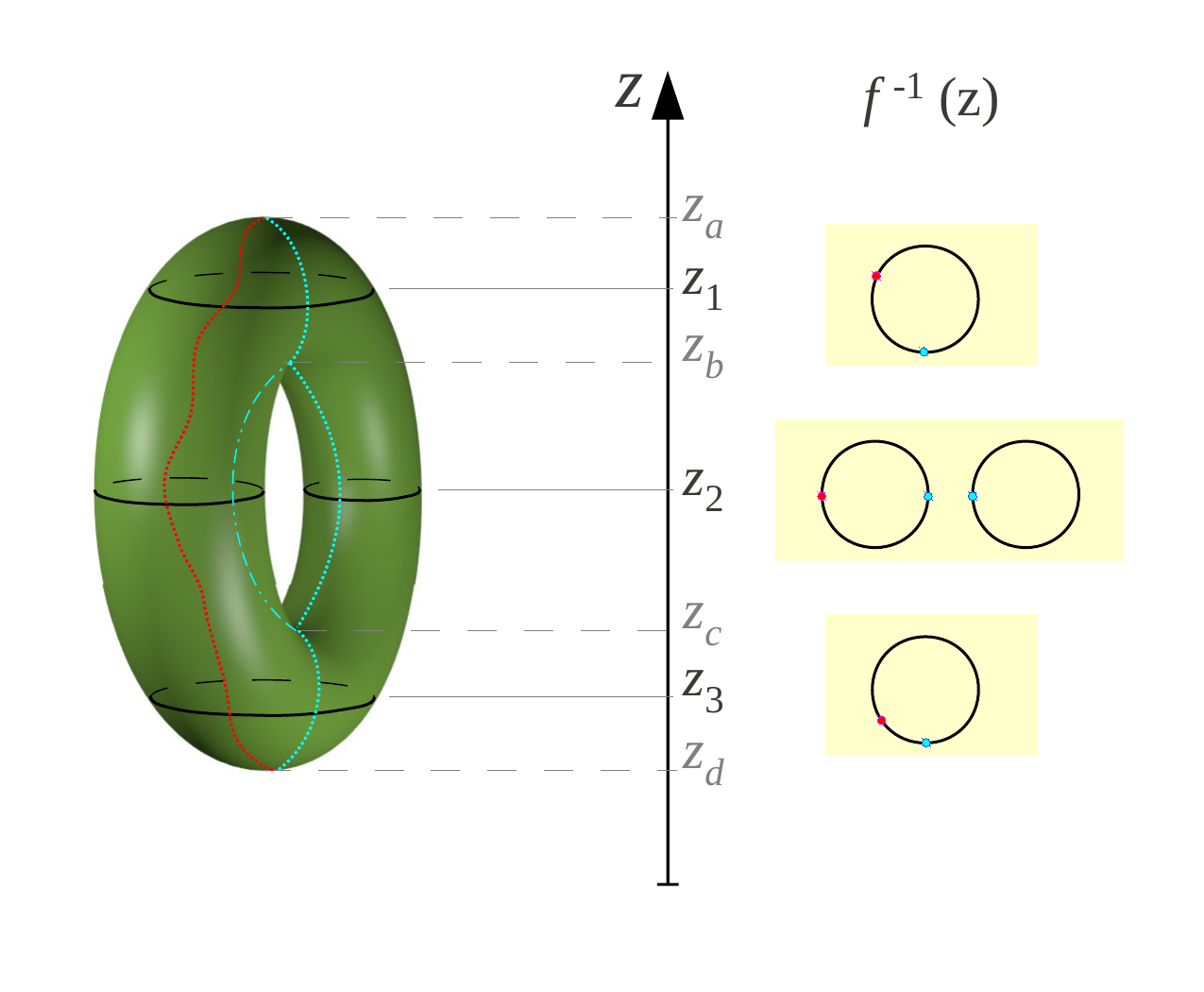}
\caption{The problem of finding a continuous global section. Red dotted curve: One configuration sampled per end-effector position. Cyan dotted curve: One configurations sampled from each connected component of $f^{-1}(z)$. Note how the later passes through the \emph{vital critical points} at the \emph{vital critical values}, $z_b, z_c$.}  \label{fig:sample-1}
\end{figure*}

We are given a robot arm, $R$. Due to the discussions in Section~\ref{sec:reorder-equivalence} it is sufficient to consider the arm with segments whose lengths have been sorted. So, without loss of generality we will assume that the segments of $R$ have lengths $l_{n-1} \geq l_{n-2} \geq \cdots \geq l_1 \geq l_0$.

The objective of this section is to design two inverse kinematics such that whenever $D_{R}^{-1}(z)$ has a single connected component the IKs agree (return the same configuration). But whenever $D_{R}^{-1}(z)$ has two components, the IKs should return two configurations, one in each of the connected components. Due to Theorem~\ref{theorem:connectivity-condition} these are the only two possibilities as far as the connectivity of $D_{R}^{-1}(z)$ is concerned --- either one connected component, or two.

It is immediately obvious that this means that we need to identify the \emph{vital critical points} in $\mathfrak{R}$, and design the inverse kinematics such that they pass through those points --- agreeing on one side of the vital critical value, while returning distinct configurations in the different connected components on other side of the critical value. This is illustrated once again in Figure~\ref{fig:sample-1} using the low-dimensional analogy with the hight function on $2$-torus. The cyan dotted curve in the figure consists of two maps, $g, \hat{g}: \mathbb{R} \rightarrow \mathbb{T}^2$, such that $g(z) = \hat{g}(z), ~\forall ~z\in [z_d, z_c] \cup [z_b, z_a]$, but when $z\in (z_c, z_b)$ we want $g(z) \neq \hat{g}(z)$ to be points in the different connected components of $f^{-1}(z)$.

\subsection{A Discrete State Representation}

As discussed, without loss of generality we consider the arm, $R$, with segments of ordered lengths, $l_{n-1} \geq l_{n-2} \geq \cdots \geq l_1 \geq l_0$.

Given a base-length, $z$, we consider the set $\{z, ~l_{n-1}, l_{n-2}, \cdots, l_0\}$ of $(n+1)$ values (the lengths of the sides of the closed polygon consisting of a fixed base-length), and arrange the elements in it in order of magnitude, so that the ordered elements are $l_\ordind{n} \geq l_\ordind{n-1} \geq l_\ordind{n-2} \geq \cdots \geq l_\ordind{1} \geq l_\ordind{0}$. Thus, $l_\ordind{n} = \max \{z, ~l_{n-1}, l_{n-2}, \cdots, l_0\}$, $l_\ordind{n-1}$ is the value from the set just less than or equal to $l_\ordind{n}$, and so on (with $l_\ordind{0} = \min \{z, ~l_{n-1}, l_{n-2}, \cdots, l_0\}$).

In the following discussions we will assume $n \geq 3$. The cases for $n \leq 2$ (\emph{i.e.} two or lesser segments) are too trivial, and computations for those easily extend from the computations for $n \geq 3$. We will however discuss such special cases separately in Note~\ref{note:n-2}.

\vspace{0.1in}

We partition the entire restricted configuration space, $\mathfrak{R}$, into some open sets (which we refer to as \emph{state blocks}) and some closed sets (which we refer to as \emph{state transitions}). In the terminology of \cite{Kapovich94onthe}, these are analogous to the ``\emph{chambers}'' and ``\emph{walls}'' respectively.

\subsubsubsection{The state blocks}
\begin{enumerate}[label=\roman*.]
 \item $\state{>}{l_{n-1}}{\infty}$: These are the states (or configurations) for which the base-length, $z = l_\ordind{n}$ (equivalently, $z\in(l_{n-1}, \infty)$), and $l_\ordind{n} + \sum_{j=0}^{n-3} l_\ordind{j} > l_\ordind{n-1} + l_\ordind{n-2}$ (equivalently, $z + \sum_{j=0}^{n-3} l_{j} > l_{n-1} + l_{n-2}$). \label{state:i}
 \item $\state{<}{l_{n-1}}{\infty}$: The states for which $z = l_\ordind{n}$, and $l_\ordind{n} + \sum_{j=0}^{n-3} l_\ordind{j} < l_\ordind{n-1} + l_\ordind{n-2}$ (equivalently, $z \in (l_{n-1},\infty)$ and $z + \sum_{j=0}^{n-3} l_{j} < l_{n-1} + l_{n-2}$). \label{state:ii}
 \item $\state{>}{l_{n-3}}{l_{n-1}}$: The states for which either $z = l_\ordind{n-1}$ or $z = l_\ordind{n-2}$, and $l_\ordind{n} + \sum_{j=0}^{n-3} l_\ordind{j} > l_\ordind{n-1} + l_\ordind{n-2}$ (equivalently, $z\in (l_{n-3},l_{n-1})$ and $l_{n-1} + \sum_{j=0}^{n-3} l_{j} > z + l_{n-2}$). \label{state:iii}
 \item $\state{<}{l_{n-3}}{l_{n-1}}$: The states for which either $z = l_\ordind{n-1}$ or $z = l_\ordind{n-2}$, and $l_\ordind{n} + \sum_{j=0}^{n-3} l_\ordind{j} < l_\ordind{n-1} + l_\ordind{n-2}$ (equivalently, $z\in (l_{n-3},l_{n-1})$ and $l_{n-1} + \sum_{j=0}^{n-3} l_{j} < z + l_{n-2}$). \label{state:iv}
 \item $\state{>}{0}{l_{n-3}}$: The states for which $z = l_\ordind{i}$ for some $i\in \{n-3, n-4, \cdots, 0\}$, and $l_\ordind{n} + \sum_{j=0}^{n-3} l_\ordind{j} > l_\ordind{n-1} + l_\ordind{n-2}$ (equivalently, $z\in (0,l_{n-3})$ and $l_{n-1} + z + \sum_{j=0}^{n-4} l_{j} > l_{n-2} + l_{n-3}$). \label{state:v}
 \item $\state{<}{0}{l_{n-3}}$: The states for which $z = l_\ordind{i}$ for some $i\in \{n-3, n-4, \cdots, 0\}$, and $l_\ordind{n} + \sum_{j=0}^{n-3} l_\ordind{j} < l_\ordind{n-1} + l_\ordind{n-2}$ (equivalently, $z\in (0,l_{n-3})$ and $l_{n-1} + z + \sum_{j=0}^{n-4} l_{j} < l_{n-2} + l_{n-3}$). \label{state:vi}
\end{enumerate}
However it is worth noting that for a given system, not all states are attainable. 
For example, if $l_{n-1} = l_{n-2} = l_{n-3} =: l$, then the states $\state{>}{l_{n-3}}{l_{n-1}}$ and $\state{<}{l_{n-3}}{l_{n-1}}$ would be empty (\emph{i.e.}, there does not exist a valid configuration satisfying the conditions of the state).

\subsubsubsection{The state transitions} 
 It is easy to note that in general the possible state transitions are the ones where either the inequality sign of the state is changed or the range of $z$ has changed from one interval to an adjacent interval. In some instances though, some of the states may be empty and can be skipped. Thus we have only the following basic state transitions (along with their inverses):
 \begin{figure*}[h]
 \centering
    \includegraphics[width=0.75\textwidth, trim=0 0 0 0, clip=true]{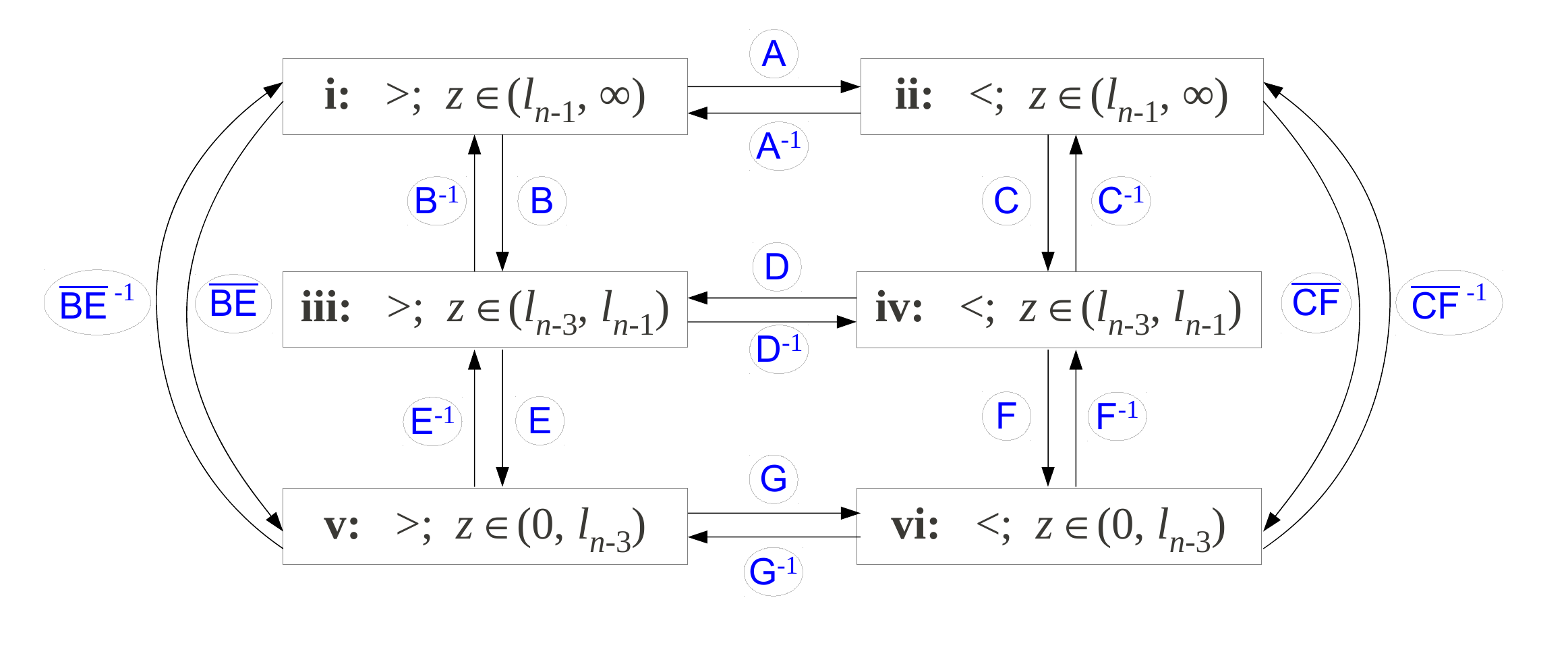}
 \caption{The possible \emph{state blocks} and \emph{state transitions}.}  \label{fig:state_transitions_full}
 \end{figure*}
\begin{enumerate}
 \item[$\zA$] and $\zA^{-1}$: This is the transitions $\state{>}{l_{n-1}}{\infty} \transitionarrow \state{<}{l_{n-1}}{\infty}$ and its inverse respectively. Thus the value of $z$ during this transition, $z_\zA$, satisfies $z_\zA + \sum_{j=0}^{n-3} l_{j} = l_{n-1} + l_{n-2}$. 

 \item[$\zB$] and $\zB^{-1}$: Transition $\state{>}{l_{n-1}}{\infty} \transitionarrow \state{>}{l_{n-3}}{l_{n-1}}$, and its inverse, can take place at $z = z_\zB := l_{n-1}$.  
 \item[$\zC$] and $\zC^{-1}$: Transition $\state{<}{l_{n-1}}{\infty} \transitionarrow \state{<}{l_{n-3}}{l_{n-1}}$ and its inverse, once again, can take place at $z = z_\zC := l_{n-1}$. 

 \item[$\zD$] and $\zD^{-1}$: This transition is $\state{>}{l_{n-3}}{l_{n-1}} \transitionarrow \state{<}{l_{n-3}}{l_{n-1}}$ and its inverse. This takes place at a value of $z = z_\zD$ such that $l_{n-1} + \sum_{j=0}^{n-3} l_{j} = z_\zD + l_{n-2}$.

 \item[$\zE$] and $\zE^{-1}$: Transition $\state{>}{l_{n-3}}{l_{n-1}} \transitionarrow \state{>}{0}{l_{n-3}}$, and its inverse, can take place at $z = z_\zE := l_{n-3}$.  
 \item[$\zF$] and $\zF^{-1}$: Transition $\state{<}{l_{n-3}}{l_{n-1}} \transitionarrow \state{<}{0}{l_{n-3}}$ and its inverse, once again, can take place at $z = z_\zF := l_{n-3}$. 

 \item[$\zG$] and $\zG^{-1}$: This transition is $\state{>}{0}{l_{n-3}} \transitionarrow \state{<}{0}{l_{n-3}}$ and its inverse. This takes place at a value of $z = z_\zG$ such that $l_{n-1} + z_\zG + \sum_{j=0}^{n-4} l_{j} =  l_{n-2} + l_{n-3}$.

 \item[$\overline{\zB\zE}$] and $\overline{\zB\zE}^{-1}$: The transition $\state{>}{l_{n-1}}{\infty} \transitionarrow \state{>}{0}{l_{n-3}}$. Possible only when $l_{n-1} = l_{n-2} = l_{n-3}$, since then the state $\state{>}{l_{n-3}}{l_{n-1}}$ will be empty. The value of $z$ at which this transition happens is thus $z_{\overline{\zB\zE}} = z_{\overline{\zB\zE}^{-1}} = l_{n-1} = l_{n-2} = l_{n-3}$. 
 \item[$\overline{\zC\zF}$] and $\overline{\zC\zF}^{-1}$: Potentially possible when $l_{n-1} = l_{n-2} = l_{n-3}$. The value of $z$ at this transition is $z_{\overline{\zC\zF}} = z_{\overline{\zC\zF}^{-1}} = l_{n-1} = l_{n-2} = l_{n-3}$.
\end{enumerate}

Note that given a state, the corresponding 
\emph{state block} or if it is at a \emph{state transition} is completely determined by the value of the base-length, $z$, alone, and not on the exact configuration of the arm.

\subsubsubsection{Impossible state transitions}
\nopagebreak

One can consider transitions which go diagonally from a state block on the left column to a state block on the right column in Figure~\ref{fig:state_transitions_full}. However some easy checks prove that such transitions are impossible.

\imphead{Impossibility of the short diagonal transitions} 
The transitions like $\state{>}{l_{n-1}}{\infty} \transitionarrow \state{<}{l_{n-3}}{l_{n-1}}$ (or their inverses) are not possible.
This is because at the transition $\state{>}{l_{n-1}}{\infty} \transitionarrow \state{<}{l_{n-3}}{l_{n-1}}$ (say which happens at $z=z_T$) we need to have, \textbf{a.} $z_T = l_{n-1}$, and for an $\epsilon \rightarrow 0_{+}$ we should also have \textbf{b.} $(z_T + \epsilon) + \sum_{j=0}^{n-3} l_{j} > l_{n-1} + l_{n-2}$ (\emph{i.e.}, close to the transition, but in state block $\state{>}{l_{n-1}}{\infty}$), and \textbf{c.} $l_{n-1} + \sum_{j=0}^{n-3} l_{j} < (z_T - \epsilon) + l_{n-2}$ (\emph{i.e.}, close to the transition, but in state block $\state{<}{l_{n-3}}{l_{n-1}}$). Combining these, one is immediately led to a contradictory statement: $l_{n-2} > \sum_{j=0}^{n-3} l_{j} + \epsilon > l_{n-2}$.

Similar arguments reveal that other transitions 
like $\state{<}{l_{n-1}}{\infty} \transitionarrow \state{>}{l_{n-3}}{l_{n-1}}$, $\state{>}{l_{n-3}}{l_{n-1}} \transitionarrow \state{<}{0}{l_{n-3}}$ and $\state{<}{l_{n-3}}{l_{n-1}} \transitionarrow \state{>}{0}{l_{n-3}}$ (and their inverses) are not possible either.

\imphead{Impossibility of the long diagonal transitions} 
In the case that the states block $\state{>}{l_{n-3}}{l_{n-1}}$ and $\state{<}{l_{n-3}}{l_{n-1}}$ are non-existent (\emph{i.e.}, $l_{n-1} = l_{n-2} = l_{n-3}$), one may consider direct transitions from $z\in (l_{n-1}, \infty)$ to $z \in (0, l_{n-3})$. While transitions like $\overline{\zB\zE}$ and $\overline{\zC\zF}$ are indeed shown in Figure~\ref{fig:state_transitions_full}, the long diagonals are in fact not possible.

To see this, we consider the transition $\state{>}{l_{n-1}}{\infty} \transitionarrow \state{<}{0}{l_{n-3}}$. Once again, suppose the transition takes place at $z = z'_T = l_{n-1} (= l_{n-2} = l_{n-3})$. Looking on either side of $z'_T$ as before, for an $\epsilon \rightarrow 0_{+}$ we should have, \textbf{a.} $(z_T + \epsilon) + \sum_{j=0}^{n-3} l_{j} > l_{n-1} + l_{n-2}$, and \textbf{b.} $l_{n-1} + (z_T + \epsilon) + \sum_{j=0}^{n-4} l_{j} < l_{n-2} + l_{n-3}$. This once again leads to a contradictory statement: $\sum_{j=0}^{n-4} l_{j} + \epsilon ~> 0 > \sum_{j=0}^{n-4} l_{j} + \epsilon$.

In a similar way we can prove that $\state{<}{l_{n-1}}{\infty} \transitionarrow \state{>}{0}{l_{n-3}}$ (or its inverse) is not possible.


\subsection{Paths Through State Transition Graph}

In this section we further refine the list of possible transitions. 
In the following proposition we consider only the state transitions as we decrease $z$ continuously starting at its maximum possible value ($z_\zS = \sum_{i=0}^{n-1} l_i$) towards zero. The state block at $z=z_\zS$ is clearly $\state{>}{l_{n-1}}{\infty}$. The objective of the proposition is to identify the state block and transitions that the system passes through as $z$ is decreased continuously, so as to enable us construct the continuous inverse kinematics as desired. 
Once again, we emphasize that without loss of generality (due to the discussions in Section~\ref{sec:reorder-equivalence}), we assume $l_{n-1} \geq l_{n-2} \geq \cdots \geq l_1 \geq l_0$.

The following proposition essentially splits the space of possible configuration spaces of robot arms into three classes -- $I$, $II$ and $III$ -- such that the type of path that the system (robot arm) follows through the vital critical values is completely determined by the class of the configuration space.
This, in the section that follows, makes it sufficient to construct only three pairs of continuous inverse kinematics -- one for each class -- such that they pass through the vital critical points, thus enabling us to achieve our design objective.

\begin{prop} \label{prop:sys-paths}
If we continuously decrease $z$, starting at $\infty$ (state block $\state{>}{l_{n-1}}{\infty}$), then 
%
%
there are only three possible classes of paths (through state transitions) that a system can take. The path classes are the ones shown in Figure~\ref{fig:state_transitions}, and the respective conditions under which a system takes one of them are:
\begin{itemize}
 \item[I.] $\zB\zE\zG$ (and $\overline{\zB\zE}\zG$ -- i.e., the possibility that $\zB$ and $\zE$ are simultaneous when $\state{>}{l_{n-3}}{l_{n-1}}$ is empty): Taken when $l_{n-2} \leq \sum_{i=0}^{n-3} l_i$.
 \item[II.] $\zA\zC\zD\zE\zG$: Taken when $l_{n-2} > \sum_{i=0}^{n-3} l_i$, and, either $n\neq 3$ or $l_{n-1} \neq l_{n-2}$.
 \item[III.] $\zA\zC\zF$: Taken when $l_{n-2} > \sum_{i=0}^{n-3} l_i$, $n=3$ and $l_{n-1} = l_{n-2}$.

 \begin{figure*}[h]
 \centering
    \includegraphics[width=0.65\textwidth, trim=0 0 0 0, clip=true]{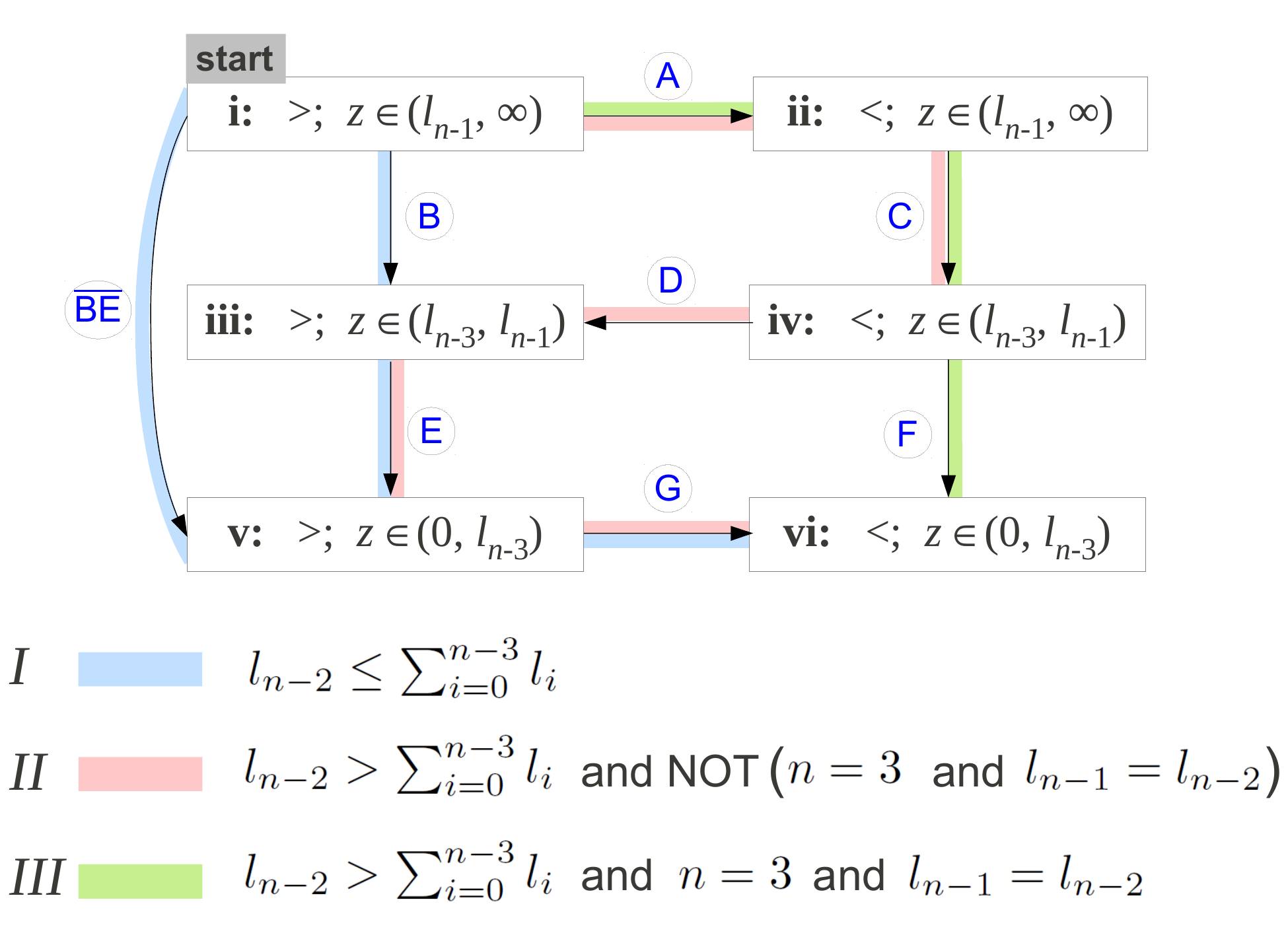}
 \caption{Figure for Proposition~\ref{prop:sys-paths}: The graph shows possible state transitions as $z$ is decreased continuously. This gives rise to three classes of paths (marked in different colors). The conditions under which the system takes these path classes are described in the legend.}  \label{fig:state_transitions}
 \end{figure*}

\end{itemize}
[Notes:
 \textbf{i)} Depending on the minimum value that $z$ can assume, the system may not traverse a path in a class in entirety as $z$ is decreased continuously. For example, if $\underline{R}(\{l_{n-1}, l_{n-2}, \cdots, l_0\}) \in (l_{n-3}, l_{n-1})$ (see \eqref{eq:range}) and $l_{n-2} \leq \sum_{i=0}^{n-3} l_i$, then the system will be able to make the transition $\zB$ but $z$ will reach its minimum possible value fore it can possibly make the transition $\zE$. 
\textbf{ii)} The path class $\zA\zC\zF$ involves an equality condition for it to be taken, making the systems that would follow a path in that class a set of measure zero in the space of all systems (a system being defined by $n$ and the set of lengths $\{l_{n-1}, l_{n-2},\cdots, l_0\}$).]
\end{prop}

\begin{quoteproof}
 The main task at hand is to prune the edges of the diagram shown in Figure~\ref{fig:state_transitions_full} and to deduce the conditions under which the state transitions take place.

\imphead{Possible transitions from $\state{>}{l_{n-1}}{\infty}$} First, it is obvious that at the highest possible value of $z$, \emph{i.e.} $z=z_\zS:=\overline{R}(\{l_{n-1}, l_{n-2}, \cdots, l_0\})=\sum_{j=0}^{n-1} l_j$, the corresponding state block is $\state{>}{l_{n-1}}{\infty}$. As we keep decreasing $z$ from this state, there are three possible transitions that can happen:
 \begin{itemize}
  \item[$\zB$:] Assume \inlinebox{$l_{n-2} \leq \sum_{i=0}^{n-3} l_i$} and \inlinebox{$l_{n-1} \neq l_{n-3}$}. Thus the state block $\state{>}{l_{n-3}}{l_{n-1}}$ is non-empty.
 
  The state blocks on either side of the transition `$\zA$' are described by $z + \sum_{j=0}^{n-3} l_{j} \lessgtr l_{n-1} + l_{n-2}$. 
  Thus the value of $z$ at which the transition $\zA$ would occur is $z_\zA = l_{n-1} + l_{n-2} - \sum_{j=0}^{n-3} l_{j}$.
  However, $\zB$ occurs at $z_\zB = l_{n-1}$. Thus, as we decrease $z$, the transition $\zB$ will occur first (instead of $\zA$) if $z_\zB > z_\zA ~\Rightarrow l_{n-2} < \sum_{i=0}^{n-3} l_i$.

  If $z_\zB = z_\zA$, then $l_{n-2} = \sum_{i=0}^{n-3} l_i$. We consider a value of $z$ just less than the transition value of $z_\zB = z_\zA = l_{n-1}$, that is, we consider $z = \lim_{\epsilon \rightarrow 0_+} (l_{n-1} - \epsilon)$. In this state block obviously $l_\ordind{n} \!=\! l_{n-1},~ l_\ordind{n-1} \!=\! z,~ l_\ordind{n-2} \!=\! l_{n-2},~ l_\ordind{n-3} \!=\! l_{n-3}, ~\text{etc}$. Now plugging these values into the inequality $l_\ordind{n} + \sum_{j=0}^{n-3} l_\ordind{j} \lessgtr l_\ordind{n-1} + l_\ordind{n-2}$ we readily verify that the state block at $z = \lim_{\epsilon \rightarrow 0_+} (l_{n-1} - \epsilon)$ in this case is indeed $\state{>}{l_{n-3}}{l_{n-1}}$. Thus transition $\zB$ (instead of $\zA$) occurs even in the case when $l_{n-2} = \sum_{i=0}^{n-3} l_i$.

  \item[$\overline{\zB\zE}$:] Assume \inlinebox{$l_{n-1} = l_{n-3}$}. Then we should also have $l_{n-1} = l_{n-2} = l_{n-3}$ (since the segment lengths are arranged in order by construction). Thus it follows that \inlinebox{$l_{n-2} \leq \sum_{i=0}^{n-3} l_i$} (equality holds only if $n = 3$. Recall we consider cases $n \geq 3$ only). 

  The transition $\overline{\zB\zE}$ would occur at $z_{\overline{\zB\zE}} = l_{n-1} (= l_{n-2} = l_{n-3})$. Thus $l_{n-2} \leq \sum_{i=0}^{n-3} l_i ~\Rightarrow~ l_{n-1} + l_{n-2} - \sum_{i=0}^{n-3} l_i \leq l_{n-1} ~\Rightarrow~ z_{\overline{\zB\zE}} \geq z_\zA$, thus implying that the value of $z$ for the transition $\overline{\zB\zE}$ is reached before (or simultaneously with) the value for transition $\zA$ is reached, as we decrease $z$ continuously.

  It can be verified once again (as was done for the case of $\zB$ above) that in the case of $z_{\overline{\zB\zE}} = z_\zA$ (possible only when $n=3$), the system enters the state block $\state{>}{l_{0}}{l_{n-3}}$ after making the transition through the value of $z = l_{n-1} (= l_{n-2} = l_{n-3})$.

  \item[$\zA$:] The only remaining case is that of \inlinebox{$l_{n-2} > \sum_{i=0}^{n-3} l_i$}. Due to the discussion above, this implies \inlinebox{$l_{n-1} \neq l_{n-3}$} (since we proved $l_{n-1} = l_{n-3} ~\Rightarrow~ l_{n-2} \leq \sum_{i=0}^{n-3} l_i$. Thus, $l_{n-1} \neq l_{n-3} ~\Leftarrow~ l_{n-2} > \sum_{i=0}^{n-3} l_i$). In this case $z_\zA > z_\zB$, thus implying that as we decreases the value of $z$, it reaches $\zA$ and transitions into state block $\state{<}{l_{n-1}}{\infty}$.
 \end{itemize}

\imphead{Possible transitions from $\state{>}{l_{n-3}}{l_{n-1}}$} As illustrated in the original state transition diagram, Figure~\ref{fig:state_transitions_full}, there are three possible transitions from this state block, namely, to the earlier state block $\state{>}{l_{n-1}}{\infty}$ (transition $\zB^{-1}$), or to $\state{<}{l_{n-3}}{l_{n-1}}$ (transition $\zD^{-1}$) or to $\state{>}{0}{l_{n-3}}$ (transition $\zE$).
 
 Transition $\zB^{-1}$ is clearly not possible if we are decreasing $z$.
 Furthermore, in state block $\state{>}{l_{n-3}}{l_{n-1}}$, the inequality that holds is $l_{n-1} + \sum_{j=0}^{n-3} l_{j} > z + l_{n-2}$. Thus, if $z$ is decreased further, and we still have $z \in (l_{n-3}, l_{n-1})$ (thus $z$ remains on the R.H.S. of the inequality), it is not possible for the inequality sign to flip.
 Thus, transition $\zD^{-1}$ is also not possible.

 However it is easy to observe that there is nothing obvious preventing transition $\zE$ when $z$ is decreased, thus making it the only possible transition.

\imphead{Possible transitions from $\state{>}{0}{l_{n-3}}$} Once again we have three possibilities: $\overline{\zB\zE}^{-1}$, $\zE^{-1}$ or $\zG$ (Figure~\ref{fig:state_transitions_full}). The first two are clearly not possible if we are decreasing $z$ continuously. Thus $\zG$ is the only possible transition from this state block.

\imphead{Possible transitions from $\state{<}{0}{l_{n-3}}$} The inequality that holds at this transition is $l_{n-1} + z + \sum_{j=0}^{n-4} l_{j} < l_{n-2} + l_{n-3}$. Moreover $z$ is not among the top three values in the ordered set of lengths, $l_\ordind{*}$ (\emph{i.e.}, $z = l_\ordind{i}$ for some $i < n-2$). Thus decreasing $z$ further will keep satisfying the inequality $l_\ordind{n} + \sum_{j=0}^{n-3} l_\ordind{j} < l_\ordind{n-1} + l_\ordind{n-2}$, keeping $z$ to the L.H.S. of it. Thus the system will remain in state block $\state{<}{0}{l_{n-3}}$ (\emph{i.e.} terminal state block, from which no transition is possible).

\imphead{Possible transitions from $\state{<}{l_{n-1}}{\infty}$} The possibilities are $\zA^{-1}$, $\zC$ and $\overline{\zC\zF}$.
 In the state block $\state{<}{l_{n-1}}{\infty}$, the inequality that holds is $z + \sum_{j=0}^{n-3} l_{j} < l_{n-1} + l_{n-2}$. Thus once again it is obvious that the transition $\zA^{-1}$ is not possible if $z$ is decreased continuously.

 For $\overline{\zC\zF}$ to happen, we need to have $l_{n-1} = l_{n-2} = l_{n-3}$. However we have seen earlier that if that is true, then, starting at state block $\state{>}{l_{n-1}}{\infty}$, as we continuously decrease $z$, the system will transition to $\state{>}{0}{l_{n-3}}$ (\emph{i.e.}, will make transition $\overline{\zB\zE}$), following which there is no possibility of transitioning to $\state{<}{l_{n-1}}{\infty}$ in the first place. So the transition $\overline{\zC\zF}$ is never possible.

 Thus the only possible transition from this state block, as $z$ is decreased continuously, is the transition $\zC$ into the state block $\state{<}{l_{n-3}}{l_{n-1}}$.

\imphead{Possible transitions from $\state{<}{l_{n-3}}{l_{n-1}}$} The possible transitions from this state block are $\zC^{-1}$, $\zD$ and $\zF$.
 Since we are only decreasing $z$ continuously, the possibility of $\zC^{-1}$ happening is eliminated immediately. The remaining two possibilities, and the conditions for them to be realized, are described below:
 \begin{itemize}
  \item[$\zF$:] If possible, say transition $\zF$ happens. Clearly $\zF$ takes place at $z = z_{\zF} := l_{n-3}$. 
 So right after the transition has taken place into the state block $\state{<}{0}{l_{n-3}}$, \emph{i.e.} at $z' = \lim_{\epsilon\rightarrow 0_+} (z_{\zF} - \epsilon)$, the condition that needs to be satisfied is $l_{n-1} + z' + \sum_{j=0}^{n-4} l_{j} < l_{n-2} + l_{n-3} ~\Rightarrow~ l_{n-1} + \sum_{j=0}^{n-4} l_{j} < l_{n-2} + \epsilon, ~\epsilon\rightarrow 0_+ $.
%
 Since by construction $l_{n-1} \geq l_{n-2}$ and by hypothesis $l_j > 0, ~\forall j$, this is possible only if \inlinebox{$n=3$} and \inlinebox{$l_{n-1}=l_{n-2}$}.

 \item[$\zD$:] This transition potentially happens at $z = z_\zD := l_{n-1} + \sum_{j=0}^{n-3} l_j - l_{n-2}$ Thus the only possible transition is $\zD$.
 If either \inlinebox{$n\neq 3$} or \inlinebox{$l_{n-1} > l_{n-2}$} then it is clear that $z_\zD > l_{n-3} (= z_\zF)$. Thus transition $\zD$ will happen before transition $\zF$ can potentially take place. 
\end{itemize}



\vspace{0.1in}
The above observations on possible transitions are combined to obtain the proposed result.
\end{quoteproof}

\begin{note}[Extending to system with $n=2$]
\label{note:n-2}
In this case there is only one possible path class that the system can take: $\zA\zC\zF$ (and the minimum possible value of $z$ will be reached before or at the transition $\zF$). 
For such a system, starting at state block $\state{>}{l_{n-1}}{\infty}$ (the state block actually being empty -- \emph{i.e.}, there does not exist a valid configuration satisfying the conditions of the state block), 
as we decrease $z$, transition $\zA$ occurs when $z = l_{n-1} + l_{n-2}$ (maximum possible value of $z$). Following that 
clearly $l_\ordind{n} < l_\ordind{n-1} + l_\ordind{n-2}$ always holds true due to triangle inequality. 
Assuming $l_{n-3} = 0$ (so that the state block $\state{<}{l_{n-3}}{l_{n-1}} \equiv \state{<}{0}{l_{n-1}}$ is meaningful), the system can also make the transition $\zC$.
For obvious reasons the transition $\zD$ is never possible for such a system (in order to satisfy triangle inequality). Moreover the state block $\state{<}{0}{l_{n-3}} \equiv \state{<}{0}{0}$ is empty. But transition $\zF$ can happen when $l_{n-1}=l_{n-2}$, since then $z$ can attain a value of $z_\zF = 0$.
\end{note}

\begin{figure*}[h]
 \centering
    \includegraphics[width=0.75\textwidth, trim=0 0 0 0, clip=true]{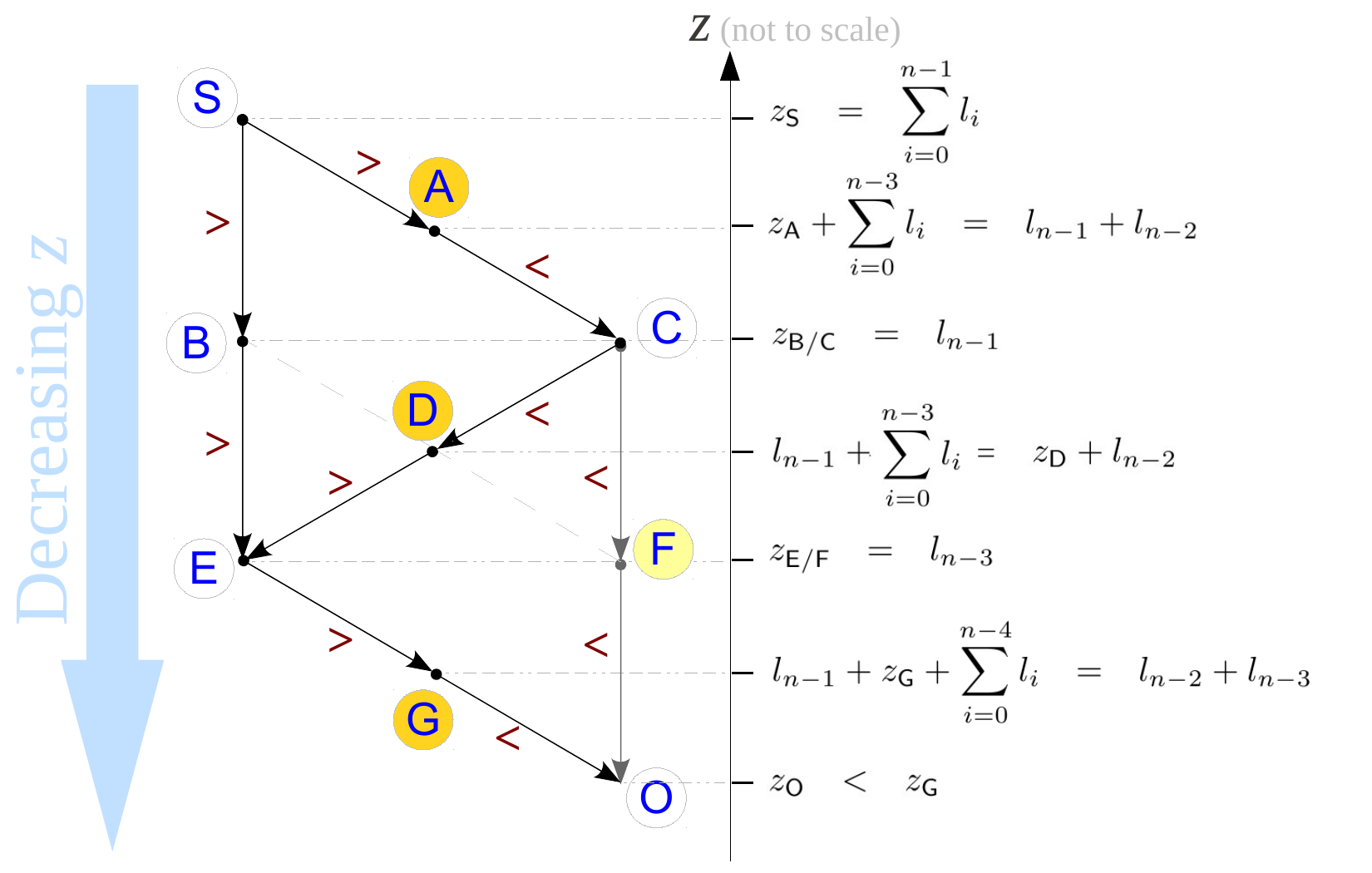}
 \caption{Graph with edges representing state block, and vertices the state transitions as $z$ decreases continuously (dual to Figure~\ref{fig:state_transitions}). The transitions in yellow are the ones corresponding to \emph{vital critical} points. 
 The transition $\zF$ is reached only in path class $\lite{\zS\,}\zA\zC\zF\lite{\,\zO}$, in which case $z_\zF$ is a vital critical value. 
          }  \label{fig:state_graph}
\end{figure*}

For easier visualization, we re-draw Figure~\ref{fig:state_transitions} as a graph with the state transitions represented by graph vertices and the state blocks by edges. This, as shown in Figure~\ref{fig:state_graph}, is a more natural (and potentially intuitive) representation than before since in reality the state transitions happen at single values of $z$ (\emph{i.e.}, points), while the state block themselves constitute extended ranges of $z$ (\emph{i.e.}, segments).
In this new figure we insert two terminal vertices: \emph{i.} A vertex, $\zS$, where the state is $z=z_\zS:=\overline{R}(\{l_{n-1}, l_{n-2}, \cdots, l_0\})=\sum_{i=0}^{n-1} l_i$, and, \emph{ii.} a vertex, $\zO$, where $z=z_\zO :=0$. 

Every system will start at $\zS$ at the highest possible value of $z=z_\zS$ (which is attainable), and due to Proposition~\ref{prop:sys-paths}, will take a path in class $\lite{\zS\,}\zB\zE\zG\lite{\,\zO}$ ($\lite{\zS\,}\overline{\zB\zE}\zG\lite{\,\zO}$ when the edge $\zB\!-\!\zE$ is of length zero), $\lite{\zS\,}\zA\zC\zD\zE\zG\lite{\,\zO}$ or $\lite{\zS\,}\zA\zC\zF\lite{\,\zO}$. However, as discussed earlier, the paths may not be traversed in their entirety since the minimum possible value of $z$ (which is given by $\underline{R}(\{l_{n-1}, l_{n-2}, \cdots, l_0\})$) may prevent it from doing so.
In fact, since we restrict $z$ to positive values (see Proposition~\ref{prop:continuous-IK}), $\zO$ will not be exactly reachable in any system. 

The \emph{vital critical points} are the transitions in which $l_\ordind{n} + \sum_{j=0}^{n-3} l_\ordind{j} = l_\ordind{n-1} + l_\ordind{n-2}$, and are marked by dark yellow in Figure~\ref{fig:state_graph}.

It is obvious that the vast majority of systems will follow paths in classes $\lite{\zS\,}\zB\zE\zG\lite{\,\zO}$ or $\lite{\zS\,}\zA\zC\zD\zE\zG\lite{\,\zO}$. The equality conditions for following path class $\lite{\zS\,}\zA\zC\zF\lite{\,\zO}$ ($n=3$ and $l_{n-1} = l_{n-2}$) makes the systems following this path class a set of measure zero in the space of all possible systems. 
However, when a system does satisfy $n=3$ and $l_{n-1} = l_{n-2}$, it is easy to check that $z_\zF$ is indeed a vital critical value, which coincides with $z_\zD$ and $z_\zG$. (Some analysis reveals that in this case the restricted configuration space, $\mathfrak{R}$, is not a manifold near the vital critical configuration $\zF$. We will refrain from further discussion on this observation in order to avoid deviating from the main topic of the paper.)



\subsection{Inverse Kinematics Design for the Three Classes of Paths}

In this section we will exploit the fact, as proven in the previous section, that there are essentially three distinct classes of paths that the system can take as the base-length, $z$ is decreased. Hence we will design one inverse kinematics strategy (using the continuous IK algorithm described in Section~\ref{IK-alg}) for each of these paths such that that it passes through the \emph{vital critical points}, enabling us to continuously sample one configuration from each connected component of the configuration space for different end effector positions.

Due to Proposition~\ref{prop:continuous-IK}, for the arm with segment lengths $\{l_{n-1}, l_{n-2}, \cdots, l_0\}$, this implies that we need to design the IKCFs, $f_{p-1}: [\underline{x}_p, \overline{x}_p] \rightarrow [\underline{x}_{p-1}, \overline{x}_{p-1}]$ (where, $\underline{x}_p = \underline{R}(\{l_0, l_1, \cdots,l_{p}\})$ and $\overline{x}_p = \overline{R}(\{l_0, l_1, \cdots,l_{p}\})$), and the corresponding IKSSS \& IKCSA pairs, $S_p$ \& $\mathbf{sg}_p$, for all $p=1,2,\cdots,n-1$.

\vspace{0.05in}
We start by noting that in Figure~\ref{fig:state_graph} the transitions colored in yellow (transitions $\zA$, $\zD$ and $\zG$) are the vital critical points, 
since the system is in a state block $\state{>}{\alpha}{\beta}$ on one side of the transition (thus, by Theorem~\ref{theorem:connectivity-condition} there is one connected component of $D^{-1}_{n+1}(z), ~\forall z\in (\alpha,\beta)$), and is in a state block $\state{<}{\mu}{\nu}$ on the other side (thus, by Theorem~\ref{theorem:connectivity-condition}, $D^{-1}_{n+1}(z)$ consists of two disconnected components for all $z\in (\mu,\nu)$).

\begin{figure*}[h]
\centering
      \includegraphics[width=0.9\textwidth, trim=50 190 50 140, clip=true]{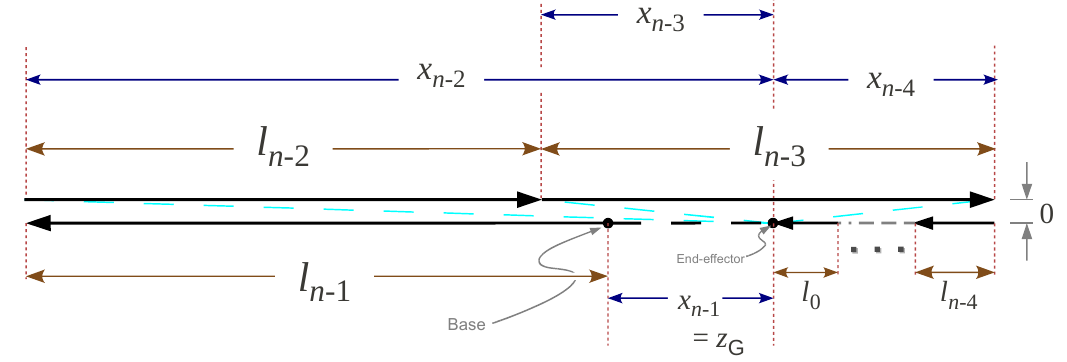}
\caption{Critical configuration at transition $\zG$, where $l_{n-2} + l_{n-3} = l_{n-1} + z_\zG +\sum_{j=0}^{n-4} l_j$. The vertical separation between the segments (marked by `$0$') is used only for easy visualization.}  \label{fig:crit-config-g}
\end{figure*}

\subsubsubsection{I: Inverse Kinematics for the path class $\lite{\zS\,}\zB\zE\zG\lite{\,\zO}$ (and $\lite{\zS\,}\overline{\zB\zE}\zG\lite{\,\zO}$)}
This class of paths is taken when $l_{n-2} \leq \sum_{j=0}^{n-3} l_j$.
On such a path the only vital critical point is the transition $\zG$.
Thus the critical design criteria is that when $z = z_\zG = l_{n-2} + l_{n-3} - l_{n-1} -\sum_{j=0}^{n-4} l_j$ (\emph{i.e.} $z$ has the \emph{vital critical value}), the configuration should be at the corresponding vital critical point (the critical configuration). 
This critical configuration is one in which all the segments are aligned along a single line as illustrated in Figure~\ref{fig:crit-config-g} (for a formal proof of this, refer to the gradient and Hessian matrix computed in the proof of Lemma 2 of \cite{jaggi:thesis:paper:92}). 
The $x_j$ marked in the figure are same as those in earlier Figure~\ref{fig:arm-recursive}, but shown in this particular critical configuration.

In order to design the complete inverse kinematics passing through the vital critical point, we first compute $\overline{x}_p$ and $\underline{x}_p$ for all $p=0,1,2,\cdots,n-1$ using Equation~\eqref{eq:range-recursive} or Proposition~\ref{prop:range-closed}.
Noting that the set $\{l_{n-1},l_{n-2},\cdots,0\}$ is already in ascending ordered, we make a few observations about $\underline{x}_p$ for $p=0,1,2,\cdots, n-1$ using Proposition~\ref{prop:range-closed} and Corollary~\ref{cor:range-min}:
\begin{itemize}
 \item[i.] $\underline{x}_{n-2} ~=~ 0$, since in this path class $l_{n-2} \leq \sum_{j=0}^{n-3} l_j$.
 \item[ii.] If $z_\zG > 0$, then $\underline{x}_{n-3} ~=~ l_{n-3} - \sum_{j=0}^{n-4} l_j > 0$ (since $z_\zG > 0, l_{n-1} \geq l_{n-2} ~\Rightarrow~ l_{n-3} > \sum_{j=0}^{n-4} l_j$, and using Proposition~\ref{prop:range-closed}).
\end{itemize}

%

\noindent If $z_\zG = l_{n-2} + l_{n-3} - l_{n-1} -\sum_{j=0}^{n-4} l_j \in [\underline{x}_{n-1}, \overline{x}_{n-1}]$, implying that the vital critical point $\zG$ is attainable, 
a configuration as shown in Figure~\ref{fig:crit-config-g} will be attainable.

Thus we need to design the IKCF-IKCSSS-IKCSA tuples, $\{f_{p-1}, \mathscr{S}_{S_p\!, \mathbf{sg}_p}\}, ~p=1,2,\cdots,n-1$.
For that, the criteria that needs to be satisfies is that when $z \equiv x_{n-1} = z_\zG$ (we assume $z_\zG \geq \underline{x}_{n-1} > 0$) we should have the values of $x_{p-1}, ~p=1,2,\cdots$ as follows (refer to Figure~\ref{fig:crit-config-g}):
\begin{eqnarray}
  x^\zG_{n-2} & = & z_\zG + l_{n-1} ~=~ l_{n-2} + l_{n-3} -\sum_{j=0}^{n-4} l_j ~~~(> l_{n-1} \geq l_{n-2}), \nonumber \\ 
  x^\zG_{n-3} & = & | x^\zG_{n-2} - l_{n-2} | ~=~ l_{n-3} - \sum_{j=0}^{n-4} l_j ~~~(=~ \underline{x}_{n-3}), \nonumber\\ 
  x^\zG_{k-1} & = & | x^\zG_{k} - l_{k} | ~=~ \sum_{j=0}^{k-1} l_j ~~~(=~ \overline{x}_{k-1}) ~, \qquad\forall k \leq n-3 \label{eq:path-I-crit-values-G}
\end{eqnarray}
where we use superscripts to indicate that these are the values at the transition $\zG$. Also for convenience we define $x^\zG_{n-1} = z_\zG$.

\begin{figure*}
\centering
\subfigure[$f^{I}_{n-2}$.]{
      \label{fig:BEG-f_n-2}
      \includegraphics[width=0.4\textwidth, trim=80 20 80 20, clip=true]{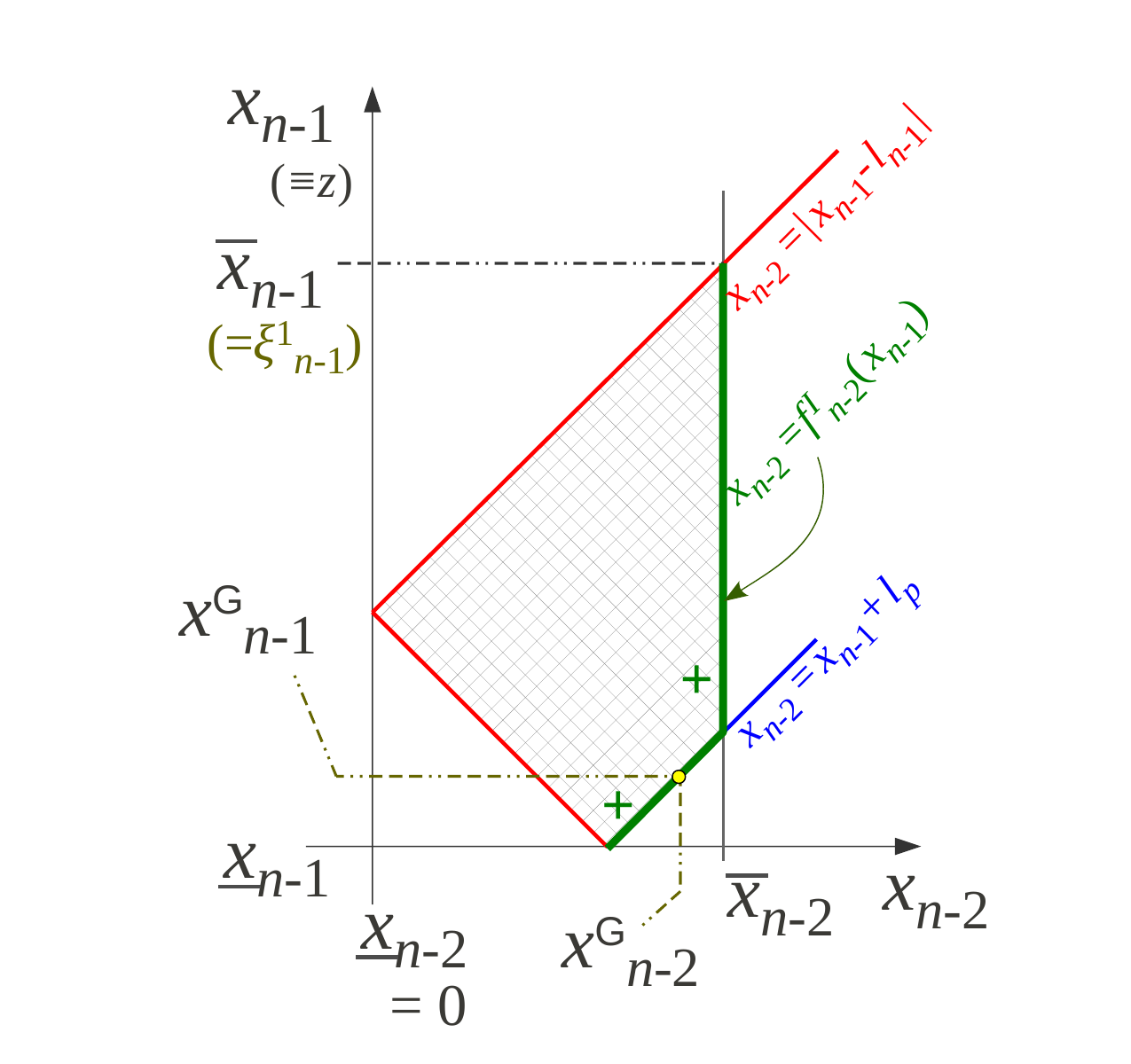}
} \hspace{-0.15in}
\subfigure[$f^{I}_{n-3}$.]{
      \label{fig:BEG-f_n-3}
      \includegraphics[width=0.4\textwidth, trim=80 20 80 20, clip=true]{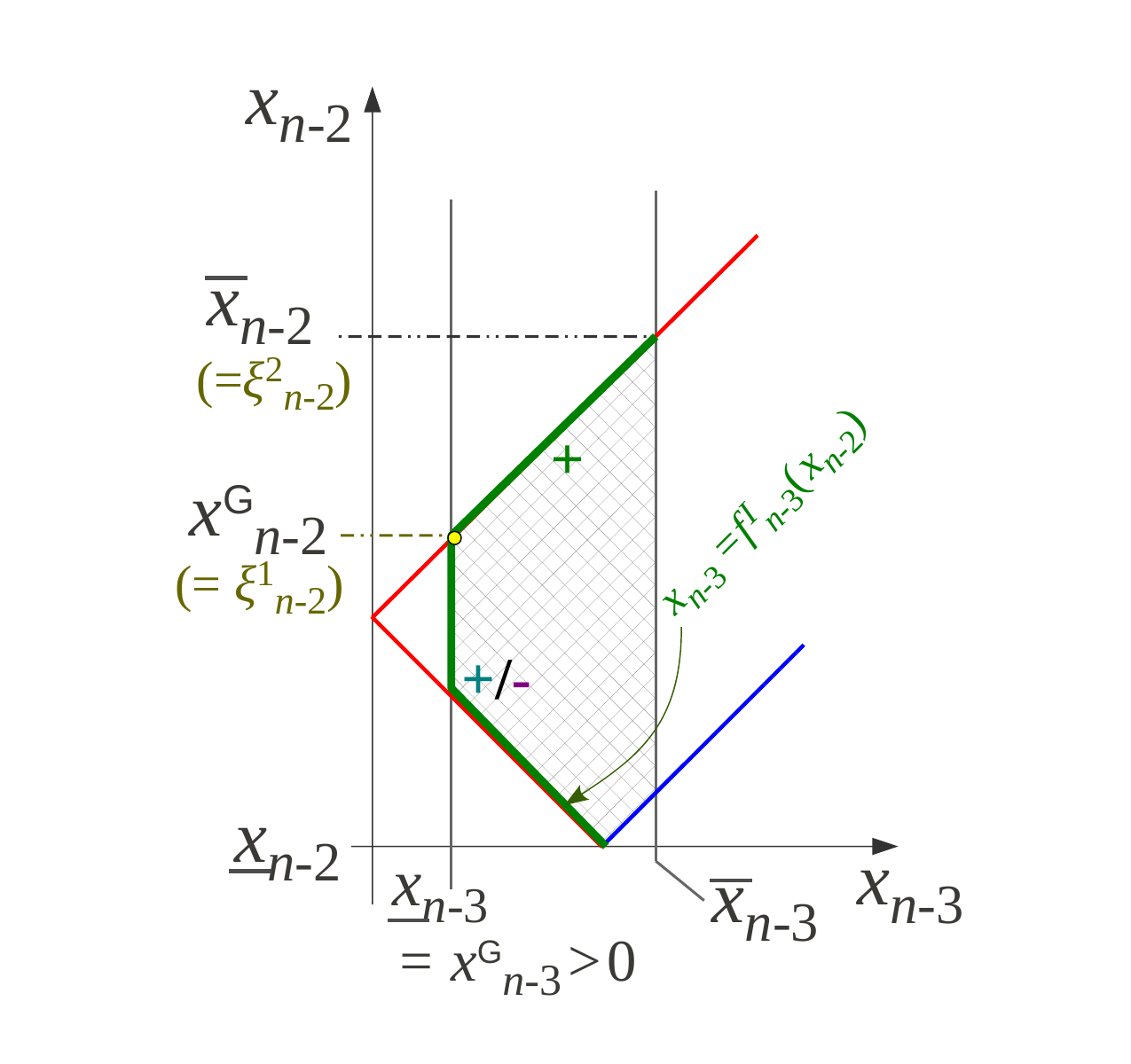}
} \hspace{-0.15in}
\subfigure[$f^{I}_{n-4}$.]{
      \label{fig:BEG-f_n-4}
      \includegraphics[width=0.4\textwidth, trim=80 20 80 20, clip=true]{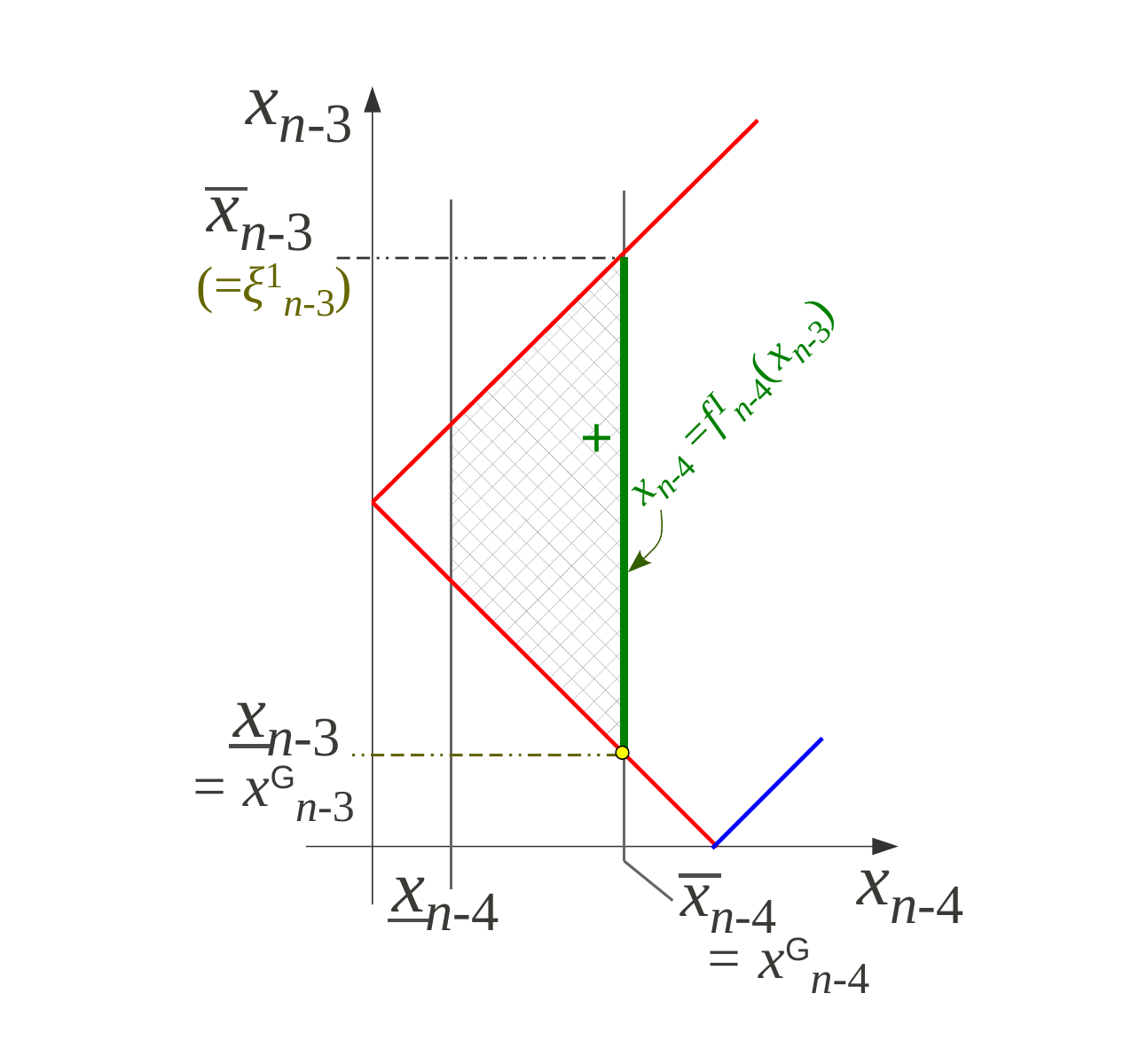}
} \hspace{-0.15in}
\subfigure[$f^{I}_{p-1}, ~p\leq n-4$.]{
      \label{fig:BEG-f_p-1}
      \includegraphics[width=0.4\textwidth, trim=80 20 80 20, clip=true]{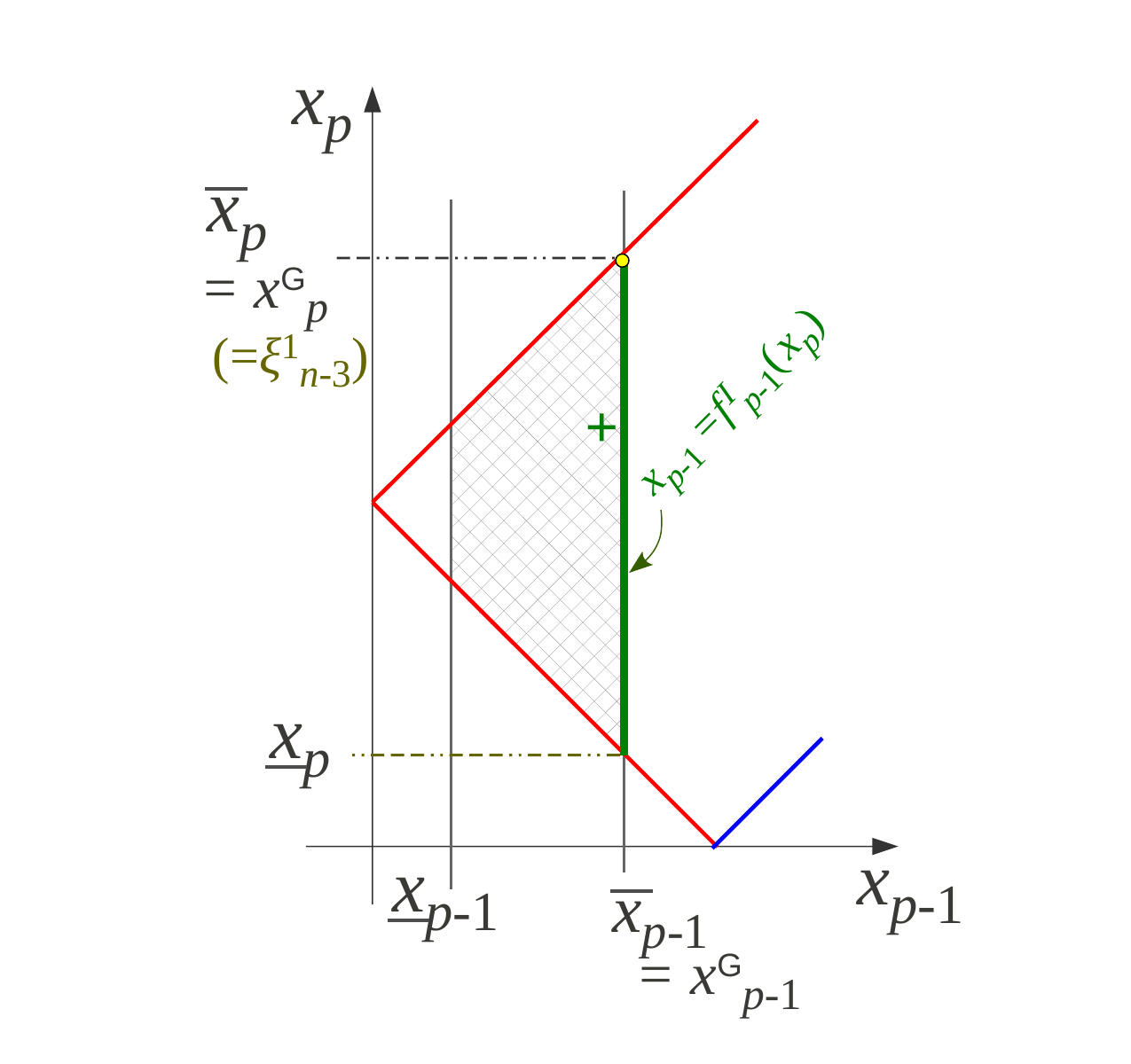}
}
\caption{Design (schematic representation) of IKCF-IKCSSS-IKCSA tuples for path class $\lite{\zS\,}\zB\zE\zG\lite{\,\zO}$. 
Note how the graphs for $f^{I}_{p-1}$ is \emph{right hugging} for $p=n-1$, while \emph{left hugging} for $p\leq n-2$, to satisfy the requirement $f_{p-1}(x^\zG_{p}) = x^\zG_{p-1}$ as determined by \eqref{eq:path-I-crit-values-G}.
(Note that in drawing these schematics we used the facts that $x^\zG_{n-2} > l_{n-2}$ and $x^\zG_{p} > l_{p}$.)}  \label{fig:BEG-all}
\end{figure*}

\begin{note}[Feasibility of $f_{p-1}$] \label{note:IKCF-feasibility}
In order to construct a $f_{p-1}$ that passes through the the point $(x^\zG_{p-1}, x^\zG_{p})$ (\emph{i.e.}, $f_{p-1}(x^\zG_{p}) = x^\zG_{p-1}$), we first need to make sure that this point lies in the corresponding feasible region in $[\underline{x}_{p-1},\overline{x}_{p-1}]\times[\underline{x}_{p},\overline{x}_{p}]$ (\emph{hatched} region in Figure~\ref{fig:range-recursive-full}). 
Clearly the points $(x^\zG_{p-1}, x^\zG_{p})$, as suggested, satisfy the conditions $|x^\zG_{p-1} - l_p| ~~\leq~~ x^\zG_p ~~\leq~~ x^\zG_{p-1} + l_p, \quad\forall p$. 
Moreover, it is easy to observe that if $z_\zG>0$, then $x^\zG_{p-1}$ has a value that can be written as $| \pm l_{p-1} \pm l_{p-2} \pm \cdots \pm l_0 |$. 
Thus due to Corollary~\ref{cor:range-min}, $\underline{x}_{p-1} \leq x^\zG_{p-1}$. 
Also for the same reason, the fact that $x^\zG_{p-1} \leq \overline{x}_{p-1} = \sum_{j=0}^{p-1} l_j$ is obvious.
%
%
%
\end{note}
The above discussions indicate the position of the points $(x^\zG_{p-1}, x^\zG_{p})$ in $[\underline{x}_{p-1},\overline{x}_{p-1}] \times [\underline{x}_{p},\overline{x}_{p}]$ as indicated in Figure~\ref{fig:BEG-all} (note that the figures assume $z_\zG > 0$).

\imphead{IKCFs}
Note from \eqref{eq:path-I-crit-values-G} (as well as referring to Figure~\ref{fig:BEG-all}), $x^\zG_{n-2}$ is the maximum that $x_{n-2}$ can attain (within the feasible bounds) when $x_{n-1} = x^\zG_{n-1}$. 
The value that $x_{n-3}$ attains is $x_{n-3}^\zG = \underline{x}_{n-3}$, which is the minimum possible when $x_{n-2} = x^\zG_{n-2}$ (assuming $z_\zG>0$, we have shown that $\underline{x}_{n-3} > 0$ as well).
Whereas for all $p\leq n-3$, $x^\zG_{p-1} = \overline{x}_{p-1}$ is the minimum that $x_{p-1}$ can attain when $x_{p} = x^\zG_{p}$.
Thus, we design the IKCFs, $f^{I}_{p-1}, ~~p=1,2,\cdots,n-1$ as follows: 
\begin{eqnarray}
%
f^{I}_{n-2}(x_{n-1}) & \!\!=\!\! & \mathsf{MAX}_{\{\overline{x}_{n-2}, l_{n-1}\}}(x_{n-1}) \nonumber \\
f^{I}_{n-3}(x_{n-2}) & \!\!=\!\! & \left\{ \begin{array}{l} \mathsf{MIN}_{\{\overline{x}_{n-3}, l_{n-2}\}}(x_{n-2}), ~~\text{if $z_\zG>0$} \\
                                                           \mathsf{MAX}_{\{\overline{x}_{n-3}, l_{n-2}\}}(x_{n-2}), ~~\text{if $z_\zG\leq 0$}
                                   \end{array} \right. \nonumber \\
f^{I}_{p-1}(x_{p}) & \!\!=\!\! & \mathsf{MAX}_{\{\overline{x}_{p-1}, l_{p}\}}(x_{p}), \qquad\forall~p \leq n-3
\label{eq:BEG-IKCFs} \end{eqnarray}
Such a construction of IKCFs will be have the property that if $z_\zG = x^\zG_{n-1} \in [\underline{x}_{n-1}, \overline{x}_{n-1}] - \{0\}$, then, due to the discussion in Note~\ref{note:IKCF-feasibility}, $f^{I}_{n-2}(x^\zG_{n-1}) = x^\zG_{n-2}, ~f^{I}_{n-3}(x^\zG_{n-2}) = x^\zG_{n-3}$, etc. 

\vspace{0.05in}
\imphead{IKCSSSs}
The choice of the IKCSSS are
\begin{equation}
S^{I}_p ~= \left\{ \begin{array}{l} \{x^\zG_{n-2}, ~\overline{x}_{n-2}\}, ~~\text{when } p=n-2 \\ \{\overline{x}_p\}, ~~\text{when } p=1,2,\cdots,n-3,n-1 \end{array} \right. \label{eq:BEG-IKCSSSs}
\end{equation}

\imphead{IKCSAs}
For every $z < z_\zG$ we need to have two configurations, one in each connected component of $D_{R}^{-1}(z)$. This is achievable by choosing two different sets of IKCSAs -- one for each connected component.
In particular, we choose the first set of IKCSAs as follows (recall, a IKCSA is a map $\mathbf{sg}_p: S_p \rightarrow \{\mathquote{+}, \mathquote{-}\}$):
\begin{eqnarray} \label{eq:ikcsa-I-1}
 \mathbf{sg}^{I}_{p} : 
\left\{ \begin{array}{l} 
 x^\zG_{p} \mapsto \mathquote{+}, ~~ \overline{x}_{p} \mapsto \mathquote{+}, ~~ \text{when } p = n-2 \\
 \overline{x}_{p} \mapsto \mathquote{+}, ~~ \forall ~p \neq n-2
\end{array} \right. 
\end{eqnarray}
and the second set of IKCSAs as
\begin{equation}
 \mathbf{\hat{sg}}^{I}_{p} : 
\left\{ \begin{array}{l} 
 x^\zG_{p} \mapsto \mathquote{-}, ~~ \overline{x}_{p} \mapsto \mathquote{+}, ~~ \text{when } p = n-2 \\
\overline{x}_{p} \mapsto \mathquote{+}, ~~ \forall ~p \neq n-2
\end{array} \right. \label{eq:ikcsa-I-2}
\end{equation}

\vspace{0.05in}
The above construction implies the interval sign function gives
$\mathscr{S}_{S^I_p, \mathbf{sg}^I_p}(x_p) = \mathquote{+}, ~\forall~p, x_p\in[\underline{x}_p,\overline{x}_p]$ using the first IKCSA of \eqref{eq:ikcsa-I-1}.
But using the second IKCSA of \eqref{eq:ikcsa-I-2}, $\mathscr{S}_{S^I_p, \mathbf{\hat{sg}}^I_p}(x_p) = \left\{ \begin{array}{l} \small \mathquote{-}, ~\text{if } p\!=\!n\!-\!2, x_{n\!-\!2} \in[\underline{x}_{n\!-\!2}, x^\zG_{n\!-\!2}) \\ \small \mathquote{+}, ~\text{otherwise.} \end{array}\right.$.
Thus the only thing that the two inverse kinematics, $\textsf{IK}_{ \{ f^{I}_{*-1}, \mathscr{S}_{S^{I}_{*}\!, \mathbf{sg}^{I}_{*}} \} }$ and $\textsf{IK}_{ \{ f^{I}_{*-1}, \mathscr{S}_{S^{I}_{*}\!, \mathbf{\hat{sg}}^{I}_{*}} \} }$, differ in is the interval sign function value corresponding to the IKCF $f_{n-3}$ in the interval $[\underline{x}_{n-2}, x^\zG_{n-2}]$ of the domain of $f_{n-3}$ (illustrated in Figure~\ref{fig:BEG-all}).

In Proposition~\ref{prop:complete-IK} we will prove that the two inverse kinematics, $\textsf{IK}_{ \{ f^{I}_{*-1}, \mathscr{S}_{S^{I}_{*}\!, \mathbf{sg}^{I}_{*}} \} }$ and $\textsf{IK}_{ \{ f^{I}_{*-1}, \mathscr{S}_{S^{I}_{*}\!, \mathbf{\hat{sg}}^{I}_{*}} \} }$, together give exactly one distinct configuration in each connected component of $D_{R}^{-1}(z)$.

\begin{note}[Case when $\zG$ is not reachable] \label{note:partial-path}
If $z_\zG = x^\zG_{n-1} \notin [\underline{x}_{n-1}, \overline{x}_{n-1}] - \{0\}$ 
(but conditions for following path class I is still satisfied) the system will always be in a $[>; *]$ state block, and thus any arbitrary IKCF-IKCSSS-IKCSA tuples, $\{f_{*-1}, \mathscr{S}_{S_*\!, \mathbf{sg}_*}\}$, returning only one configuration for every value of $z$ will be sufficient.
This will indeed be the case since $x^\zG_{n-1} \notin [\underline{x}_{n-1}, \overline{x}_{n-1}] ~\Rightarrow~ x^\zG_{n-2} \notin [\underline{x}_{n-2}, \overline{x}_{n-2}]$, and thus rendering $\mathscr{S}_{S^I_p\!, \mathbf{sg}_{n-2}}$ and $\mathscr{S}_{S^I_p\!, \mathbf{\hat{sg}}_{n-2}}$ identical in their entire domain $[\underline{x}_{n-2}, \overline{x}_{n-2}]$ (refer to Definition~\ref{def:ikcisf}).
Furthermore, the specialized definition of $f^{I}_{n-3}$ for $z_\zG \leq 0$ will ensure that we do not encounter a configuration with $x_{n-3} = 0$.
\end{note}

Note that the above constructions of IKCF-IKCSSS-IKCSA tuples are completely indifferent to whether or not transitions $\zB$ and $\zE$ are distinct or simultaneous (when $l_{n-1}=l_{n-3}$). The only fact that we utilized is that $\zG$ is the first (and only) \emph{vital critical point} that is encountered as we decrease $z$ starting from its highest value. This is true for both paths $\lite{\zS\,}\zB\zE\zG\lite{\,\zO}$ and $\lite{\zS\,}\overline{\zB\zE}\zG\lite{\,\zO}$.

\vspace{0.1in}
We next construct similar IKCF-IKCSSS-IKCSA tuples for inverse kinematics for the path classes $\lite{\zS\,}\zA\zC\zD\zE\zG\lite{\,\zO}$ and $\lite{\zS\,}\zA\zC\zF\lite{\,\zO}$. 

\begin{figure*}[h]
\centering
\subfigure[Critical configuration at transition $\zA$. Note that $z_\zA > l_{n-1}$ since $z_\zA = l_{n-1} + l_{n-2} - \sum_{j=0}^{n-3} l_j$ and $l_{n-2} > \sum_{j=0}^{n-3} l_j$ for path classes $\lite{\zS\,}\zA\zC\zD\zE\zG\lite{\,\zO}$ or $\lite{\zS\,}\zA\zC\zF\lite{\,\zO}$ -- the only classes where transition $\zA$ is possible.]{
      \label{fig:crit-config-a}
      \includegraphics[width=0.9\textwidth, trim=40 190 120 160, clip=true]{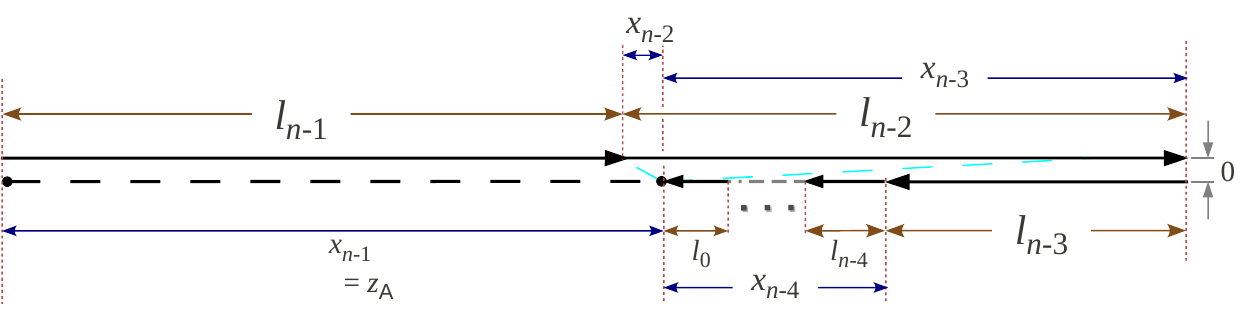}
}
\subfigure[Critical configuration at transition $\zD$.]{
      \label{fig:crit-config-d}
      \includegraphics[width=0.55\textwidth, trim=40 160 60 130, clip=true]{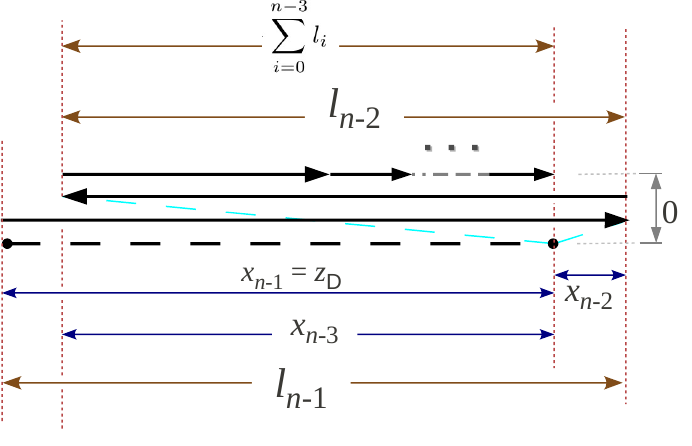}
}
\caption{Vital critical configurations $\zA$ and $\zD$.}  \label{fig:crit-configs2}
\end{figure*}

\subsubsubsection{II: Inverse Kinematics for the path class $\lite{\zS\,}\zA\zC\zD\zE\zG\lite{\,\zO}$}
For this path class to be followed we need to have $l_{n-2} > \sum_{j=0}^{n-3} l_j$. 
The vital critical values are ~$z_\zA = l_{n-1} + l_{n-2} - \sum_{j=0}^{n-3} l_j$, ~~$z_\zD = l_{n-1} - l_{n-2} + \sum_{j=0}^{n-3} l_j$, and, ~~$z_\zG = l_{n-2} + l_{n-3} - l_{n-1} -\sum_{j=0}^{n-4} l_j$. The respective vital critical configurations are shown in Figures~\ref{fig:crit-config-a}, \ref{fig:crit-config-d} and \ref{fig:crit-config-g}.

Once again, a few observations about $\underline{x}_{*}$:
\begin{itemize}
 \item[i.] $\underline{x}_{n-2} = l_{n-2} - \sum_{j=0}^{n-3} l_j > 0$ (since on this path class $l_{n-2} > \sum_{j=0}^{n-3} l_j$, and using Proposition~\ref{prop:range-closed}),
 \item[ii.] If $z_\zG > 0$, then $\underline{x}_{n-3} ~=~ l_{n-3} - \sum_{j=0}^{n-4} l_j > 0$ (same as the previous path class),
\end{itemize}

Referring to Figure~\ref{fig:crit-config-a} we have for the vital critical point $\zA$,
\begin{eqnarray}
 x^\zA_{n-1} & = & z_\zA  \nonumber \\ 
 x^\zA_{n-2} & = &  | x^\zA_{n-1} - l_{n-1} | ~=~ l_{n-2} - \sum_{j=0}^{n-3} l_j ~~~(=~ \underline{x}_{n-2}) \nonumber \\ 
 x^\zA_{k-1} & = & | x^\zA_{k} - l_{k} | ~=~ \sum_{j=0}^{k-1} l_j ~~~(=~ \overline{x}_{k-1}) 
, \quad\forall ~k \leq n-2 \label{eq:crit-pts-A-condition}
\end{eqnarray}

Likewise, referring to Figure~\ref{fig:crit-config-d} we have for the vital critical point $\zD$,
\begin{eqnarray}
 x^\zD_{n-1} & = & z_\zD  \nonumber \\
 x^\zD_{n-2} & = & | x^\zD_{n-1} - l_{n-1} | ~=~ l_{n-2} - \sum_{j=0}^{n-3} l_j ~~~(=~ \underline{x}_{n-2}) \nonumber \\ 
   x^\zD_{k-1} & = & | x^\zD_{k} - l_{k} | ~=~ \sum_{j=0}^{k-1} l_j ~~~(=~ \overline{x}_{k-1}), \quad\forall ~k \leq n-2
\end{eqnarray}

\begin{figure*}[h]
\centering
\subfigure[$f^{II}_{n-2}$.]{ \label{fig:acdeg-fn-2}
      \includegraphics[width=0.5\textwidth, trim=80 40 80 20, clip=true]{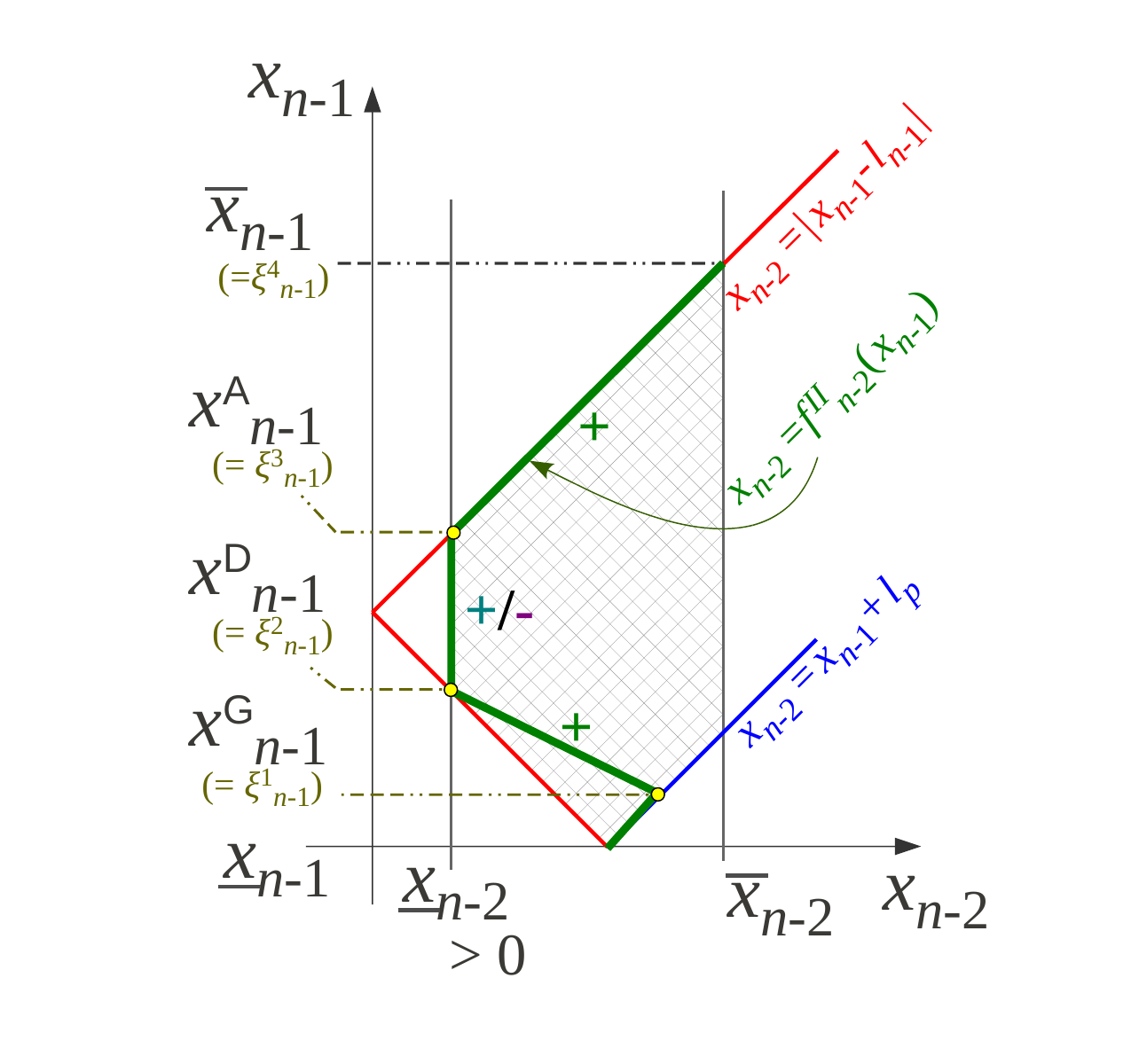}
} 
\subfigure[$f^{II}_{p-1}, ~p\leq n-2$.]{ \label{fig:acdeg-fn-3-p}
\includegraphics[width=0.22\textwidth, trim=80 0 80 120, clip=true]{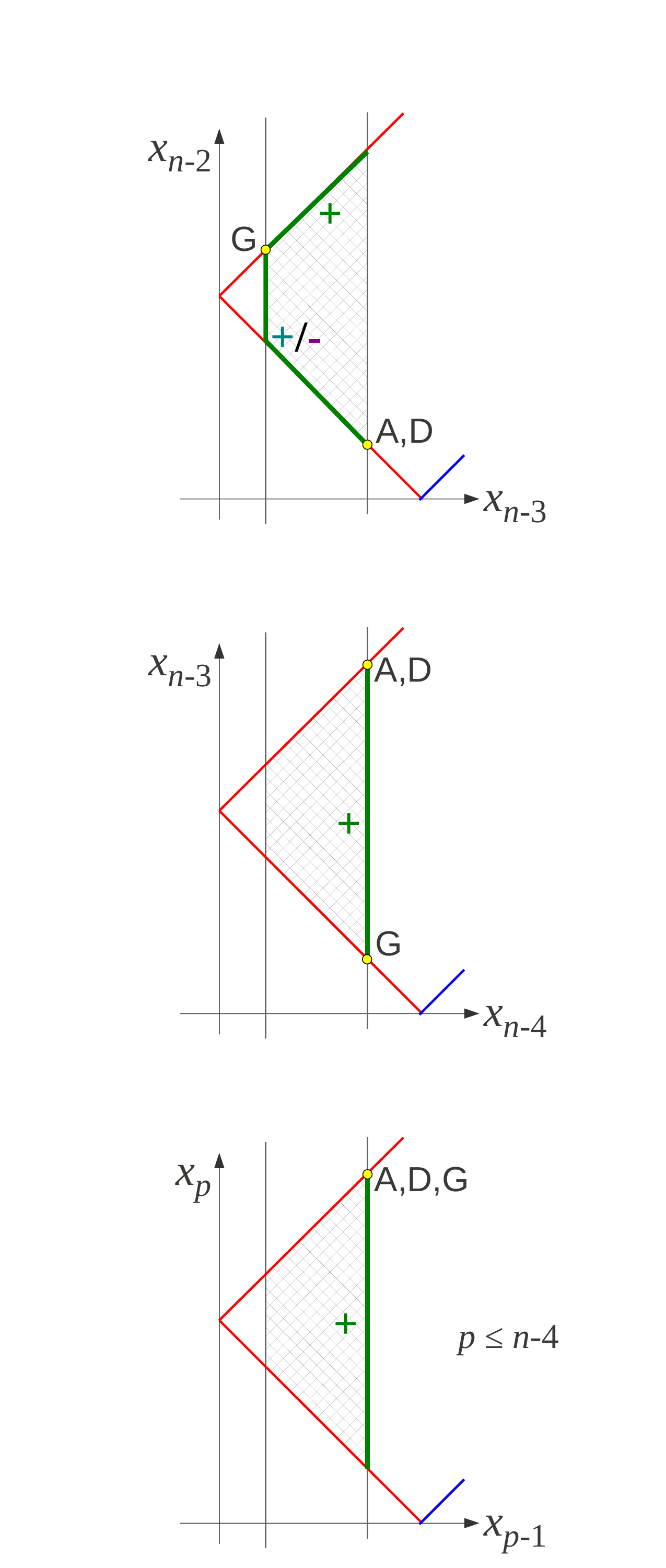}
}
\caption{Schematic representation of the IKCF, and corresponding IKCSSS and IKCSA for path class $\lite{\zS\,}\zA\zC\zD\zE\zG\lite{\,\zO}$.}  
\end{figure*}

It is interesting to note that $x^\zA_{n-2} = x^\zD_{n-2} = \underline{x}_{n-2}$. This, along with the fact that $x^\zA_{n-2} =  | x^\zA_{n-1} - l_{n-1} |$ and $x^\zD_{n-2} =  | x^\zD_{n-1} - l_{n-1} |$ indicate the positions of $(x^\zA_{n-2}, x^\zA_{n-1})$ and $(x^\zD_{n-2}, x^\zD_{n-1})$ as illustrated in Figure~\ref{fig:acdeg-fn-2}.


Using the same arguments as in Note~\ref{note:IKCF-feasibility}, we can show that if $z_\zA > 0$, then the points $(z^\zA_{p-1}, z^\zA_{p})$ lie inside the feasible region for the graph of $f^{II}_{p-1}$. Likewise for $\zD$.
 Moreover, as before, for all $p=1,2,\cdots,n-1$, when $x_p = x^\zA_{p}$, the value of $x_{p-1}$ is the minimum that it can attain inside the feasible region. The same is true for $\zD$. For $\zG$ we still need to satisfy \eqref{eq:path-I-crit-values-G}.
Furthermore, we notice that $x^\zG_{n-1} < x^\zD_{n-1} < x^\zA_{n-1}$ for this path class.
%
Thus we construct the following IKCFs (refer to Figure~\ref{fig:acdeg-fn-2}):
\begin{eqnarray}
f^{II}_{n-2}(x_{n-1}) & \!\!=\!\! & \mathsf{STEP}_{\{\underline{x}_{n-2}, \overline{x}_{n-2}, l_{n-1}, x^\zG_{n-1}, x^\zD_{n-1}\}}(x_{n-1}) \nonumber \\
f^{II}_{n-3}(x_{n-2}) & \!\!=\!\! & \left\{ \begin{array}{l} \mathsf{MIN}_{\{\overline{x}_{n-3}, l_{n-2}\}}(x_{n-2}), ~~\text{if $z_\zG>0$} \\
                                                           \mathsf{MAX}_{\{\overline{x}_{n-3}, l_{n-2}\}}(x_{n-2}), ~~\text{if $z_\zG\leq 0$}
                                   \end{array} \right. \nonumber \\
f^{II}_{p-1}(x_{p}) & \!\!=\!\! & \mathsf{MAX}_{\{\overline{x}_{p-1}, l_{p}\}}(x_{p}), \quad\forall~p \leq n-3
\label{eq:ACDEG-IKCFs}
\end{eqnarray}

\noindent And the IKCSSSs:
\begin{equation}
S^{II}_p ~= \left\{ \begin{array}{l} \{x^\zD_{n-1}, x^\zA_{n-1}, \overline{x}_{n-1} \}, ~~\text{when } p=n-1 \\ 
                                     \{x^\zG_{n-2}, \overline{x}_{n-2} \}, ~~\text{when } p=n-2 \\ 
                                     \{\overline{x}_p\}, ~~\text{when } p=1,2,\cdots,n-3
                    \end{array} \right. \label{eq:ACDEG-IKCSSSs}
\end{equation}

\noindent And two sets of IKCSAs, $\mathbf{sg}^{II}_*$ and $\mathbf{\hat{sg}}^{II}_*$ (one each for obtaining one configuration from each connected component of $D_{n+1}^{-1}(z)$):
\begin{equation}
 \mathbf{sg}^{II}_p: 
\left\{ \begin{array}{l} 
  x^\zD_{p} \mapsto \mathquote{+}, ~~x^\zA_{p} \mapsto \mathquote{+}, ~~\overline{x}_{p} \mapsto \mathquote{+}, ~~\text{when } p=n-1 \\
  x^\zG_{p} \mapsto \mathquote{+}, ~~\overline{x}_{p} \mapsto \mathquote{+}, ~~\text{when } p=n-2 \\
  \overline{x}_{p} \mapsto \mathquote{+}, ~~\forall p \in \{1,2,\cdots,n-3\}
\end{array} \right.
\end{equation}
and
\begin{equation}
 \mathbf{\hat{sg}}^{II}_p: 
\left\{ \begin{array}{l} 
  x^\zD_{p} \mapsto \mathquote{+}, ~~x^\zA_{p} \mapsto \mathquote{-}, ~~\overline{x}_{p} \mapsto \mathquote{+}, ~~\text{when } p=n-1 \\
  x^\zG_{p} \mapsto \mathquote{-}, ~~\overline{x}_{p} \mapsto \mathquote{+}, ~~\text{when } p=n-2 \\
  \overline{x}_{p} \mapsto \mathquote{+}, ~~\forall p \in \{1,2,\cdots,n-3\}
\end{array} \right.
\end{equation}

The above construction of IKCSSSs and IKCSAs give interval sign functions \\
$\mathscr{S}_{S^I_p, \mathbf{sg}^I_p}(x_p) = \mathquote{+}, ~\forall ~x_p\in[\underline{x}_p,\overline{x}_p], p\in\{1,2,\cdots,n-2,n-1\}$, and, \\
{\small $\mathscr{S}_{S^I_p, \mathbf{\hat{sg}}^I_p}(x_p) = \left\{ \begin{array}{l} \small \mathquote{-}, ~\text{if }~ p\!=\!n\!-\!1, x_{n\!-\!1} \in[x^\zG_{n\!-\!1}, x^\zD_{n\!-\!1}) ~\text{ or }~ p\!=\!n\!-\!2, x_{n\!-\!2} \in[\underline{x}_{n\!-\!2}, x^\zG_{n\!-\!2}) \\ \small \mathquote{+}, ~\text{otherwise.} \end{array}\right.$}.

Similar to the argument in Note~\ref{note:partial-path}, if any of $z_\zA, z_\zD$ or $z_\zG$ is not attainable, even then the inverse kinematics due to the above designed IKCF-IKCSSS-IKCSA tuples, $\mathsf{IK}_{\{ f^{II}_{*-1}, \mathscr{S}_{S^{II}_{*}\!, \mathbf{sg}^{II}_{*}} \}}$ and $\mathsf{IK}_{\{ f^{II}_{*-1}, \mathscr{S}_{S^{II}_{*}\!, \mathbf{\hat{sg}}^{II}_{*}} \}}$, will be sufficient. 

\subsubsubsection{III: Inverse Kinematics for the path class $\lite{\zS\,}\zA\zC\zF\lite{\,\zO}$}
In this path class we once again have $l_{n-2} > \sum_{j=0}^{n-3} l_j$. Thus we still have $\underline{x}_{n-2} = l_{n-2} - \sum_{j=0}^{n-3} l_j > 0$. 
Also, in this case $\underline{x}_{n-3} = \overline{x}_{n-3} = l_{n-3}$.

\begin{figure*}[h]
\centering
      \includegraphics[width=0.65\textwidth, trim=40 140 60 130, clip=true]{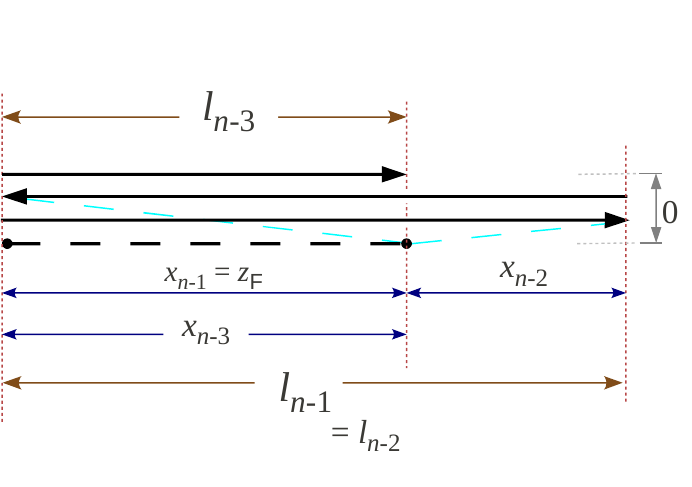}
\caption{Critical configuration at transition $\zF$, where $l_{n-2} + l_{n-3} = l_{n-1} + z_\zF$.}  \label{fig:crit-config-f}
\end{figure*}

The vital critical points are $\zA$ and $\zF$. 
Treatment of the $\zA$ is done as before. The vital critical configuration at $z = z_\zF = l_{n-1} - l_{n-2} + l_{n-3}$ is shown in Figure~\ref{fig:crit-config-f}. This vital critical value is different from the others in the sense that the system is in a state block of the form $[<; *]$ on either side of $z_\zF$ (\emph{i.e.} for both $z_\zF + \epsilon$ and $z_\zF - \epsilon$).

Writing down the values of $x^\zF_*$,
\begin{eqnarray}
 x^\zF_{n-1} & = & l_{n-2} + l_{n-3} - l_{n-1} ~~(= l_{n-3}) \nonumber \\
 x^\zF_{n-2} & = & |x^\zF_{n-1} - l_{n-1}| ~~(= l_{n-2} - l_{n-3} = \underline{x}_{n-2}) \nonumber \\
 x^\zF_{n-3} & = & |x^\zF_{n-2} - l_{n-2}| ~~(= l_{n-3} = \overline{x}_{n-3})
\end{eqnarray}

Keeping the above and \eqref{eq:crit-pts-A-condition} in view, similar to before we construct the IKCFs as follows 
\begin{eqnarray}
f^{III}_{n-2}(x_{n-1}) & \!\!=\!\! & \mathsf{MIN}_{\{\overline{x}_{n-2}, l_{n-1}\}}(x_{n-1}) \nonumber \\
f^{III}_{p-1}(x_{p}) & \!\!=\!\! & \mathsf{MAX}_{\{\overline{x}_{p-1}, l_{p}\}}(x_{p}), ~~p = n-2, n-3, \cdots, 1
\end{eqnarray}
The IKCSSS-IKCSA constructions follow similar to that of path class $\lite{\zS\,}\zA\zC\zD\zE\zG\lite{\,\zO}$.

\noindent The IKCSSSs:
\begin{equation}
S^{III}_p ~= \left\{ \begin{array}{l} \{x^\zF_{n-1}, x^\zA_{n-1}, \overline{x}_{n-1} \}, ~~\text{when } p=n-1 \\ 
                                     \{x^\zF_{n-2}, \overline{x}_{n-2} \}, ~~\text{when } p=n-2 \\ 
                                     \{\overline{x}_p\}, ~~\text{when } p=1,2,\cdots,n-3
                    \end{array} \right.
\end{equation}

\noindent And two sets of IKCSAs:
\begin{equation}
 \mathbf{sg}^{III}_p: 
\left\{ \begin{array}{l} 
  x^\zF_{p} \mapsto \mathquote{+}, ~~x^\zA_{p} \mapsto \mathquote{+}, ~~\overline{x}_{p} \mapsto \mathquote{+}, ~~\text{when } p=n-1 \\
  x^\zF_{p} \mapsto \mathquote{+}, ~~\overline{x}_{p} \mapsto \mathquote{+}, ~~\text{when } p=n-2 \\
  \overline{x}_{p} \mapsto \mathquote{+}, ~~\forall p \in \{1,2,\cdots,n-3\}
\end{array} \right.
\end{equation}
and
\begin{equation}
 \mathbf{\hat{sg}}^{III}_p: 
\left\{ \begin{array}{l} 
  x^\zF_{p} \mapsto \mathquote{+}, ~~x^\zA_{p} \mapsto \mathquote{-}, ~~\overline{x}_{p} \mapsto \mathquote{+}, ~~\text{when } p=n-1 \\
  x^\zF_{p} \mapsto \mathquote{-}, ~~\overline{x}_{p} \mapsto \mathquote{+}, ~~\text{when } p=n-2 \\
  \overline{x}_{p} \mapsto \mathquote{+}, ~~\forall p \in \{1,2,\cdots,n-3\}
\end{array} \right.
\end{equation}


\vspace{0.05in}
\begin{prop} \label{prop:complete-IK}
 The two continuous inverse kinematics thus designed for each path class,
 \[ \textsf{IK}_{ \{ f^{\#}_{*-1}, \mathscr{S}_{S^{\#}_{*}\!, \mathbf{sg}^{\#}_{*}} \} }~ \text{ and } ~\textsf{IK}_{ \{ f^{\#}_{*-1}, \mathscr{S}_{S^{\#}_{*}\!, \mathbf{\hat{sg}}^{\#}_{*}} \} }\] 
 [where,
 $\# \equiv I$ for path class $\lite{\zS\,}\zB\zE\zG\lite{\,\zO}$ (and $\lite{\zS\,}\overline{\zB\zE}\zG\lite{\,\zO}$) -- used when $l_{n-2} \leq \sum_{j=0}^{n-3} l_j$; 
~$\# \equiv II$ for path class $\lite{\zS\,}\zA\zC\zD\zE\zG\lite{\,\zO}$ -- used when $l_{n-2} > \sum_{j=0}^{n-3} l_j$ and, either $n\neq 3$ or $l_{n-2} \neq l_{n-3}$; 
~$\# \equiv III$ for path class $\lite{\zS\,}\zA\zC\zF\lite{\,\zO}$ -- used when $l_{n-2} > \sum_{j=0}^{n-3} l_j, n=3$ and $l_{n-2} = l_{n-3}$],\\
%
have the properties that
\begin{itemize}
 \item[i.] They pass through the vital critical points on the path at the respective vital critical values of base segment length, $z$,
 \item[ii.] They agree when the length of base segment, $z$, is such that $D_{R}^{-1}(z)$ has a single connected component, and,
 \item[iii.] Give configurations in the two different connected components of $D_{R}^{-1}(z)$ whenever $D_{R}^{-1}(z)$ has two connected components. 
\end{itemize}
\end{prop}
\begin{quoteproof}
 We first note that the two inverse kinematics for any path class, $\textsf{IK}_{ \{ f^{\#}_{*-1}, \mathscr{S}_{S^{\#}_{*}\!, \mathbf{sg}^{\#}_{*}} \} }$ and $\textsf{IK}_{ \{ f^{\#}_{*-1}, \mathscr{S}_{S^{\#}_{*}\!, \mathbf{\hat{sg}}^{\#}_{*}} \} }$, only differ in the choice of the IKCSAs.
 
 \emph{i.} is due to construction of the IKCFs, $f^{\#}_{*-1}$.

 \emph{ii.} is easy to observe since $\mathbf{sg}^{\#}_{p}$ and $\mathbf{\hat{sg}}^{\#}_{p}$ agree on $\overline{x}_p$ and $x^\zD_p$ (the later is relevant only for path $\lite{\zS\,}\zA\zC\zD\zE\zG\lite{\,\zO}$). Thus the sign functions, $\mathscr{S}_{S_p\!,\mathbf{sg}_p}$ and $\mathscr{S}_{S_p\!,\mathbf{\hat{sg}}_p}$, and hence the two inverse kinematics agree in all 
state blocks of the form $[>; *]$, \emph{i.e.}, $z$ such that $D_{R}^{-1}(z)$ has a single connected component 
(see the construction of functions $\Theta_p$ and $\Phi_p$ in Lemma~\ref{lemma:trig-fun-signed} and the inverse kinematics algorithm of Section~\ref{sec:recursive-IK}).

\vspace{0.05in}
 To prove \emph{iii.} it is sufficient to pick one particular base-length, $z$, in every interval in which the state block is of the form $[<; *]$ (\emph{i.e.}, $[\underline{x}_{n-1}, z_\zG)$ for path classes $\lite{\zS\,}\zB\zE\zG\lite{\,\zO}$ \& $\lite{\zS\,}\zA\zC\zD\zE\zG\lite{\,\zO}$; $(z_\zD, z_\zA)$ for path class $\lite{\zS\,}\zA\zC\zD\zE\zG\lite{\,\zO}$; and, $[\underline{x}_{n-1}, z_\zF)$ \& $(z_\zF, z_\zA)$ for path class $\lite{\zS\,}\zA\zC\zF\lite{\,\zO}$), and show that the two inverse kinematics return two distinct configurations, one in each connected component of $D_{R}^{-1}(z)$.
 Then, due to the continuity of either of the inverse kinematics, they will return one configuration in each of the connected components throughout the interval.
%
 Thus, we look at the two configurations returned by the two inverse kinematics when the base-length is $z_\zG - \epsilon$ (for path classes $\lite{\zS\,}\zB\zE\zG\lite{\,\zO}$ and $\lite{\zS\,}\zA\zC\zD\zE\zG\lite{\,\zO}$), $z_\zA - \epsilon$ (for path classes $\lite{\zS\,}\zA\zC\zD\zE\zG\lite{\,\zO}$ and $\lite{\zS\,}\zA\zC\zF\lite{\,\zO}$) and $z_\zF - \epsilon$ (for path class $\lite{\zS\,}\zA\zC\zF\lite{\,\zO}$), with $\epsilon\rightarrow 0_+$. 


We illustrate the proof for $\lim_{\epsilon\rightarrow 0_+} z_\zG - \epsilon$.
The system at $z = z_\zG - \epsilon$ is in state block $\state{<}{0}{l_{n-3}}$ (see Figure~\ref{fig:state_transitions}). Due to the construction of $f^{I}_{n-2}, f^{I}_{n-4}, f^{I}_{n-5}, \cdots, f^{I}_{0}$ or $f^{II}_{n-2}, f^{II}_{n-4}, f^{II}_{n-5}, \cdots, f^{II}_{0}$ it is easy to observe that (see Figure~\ref{fig:BEG-all},~\ref{fig:acdeg-fn-2}) when $x_{n-1} = z_\zG \!-\! \epsilon = x_{n-1}^\zG \!-\! \epsilon$ for some small positive $\epsilon$, we still have $x_{n-2} = x_{n-1} + l_{n-1}$, ~$x_{n-3} = \underline{x}_{n-3}$ and $x_p = \overline{x}_p, ~p=n-4, \cdots, 1, 0$.
Thus the diagram below is the schematic of the configuration when $z = z_\zG - \epsilon$.

\begin{center}
 \includegraphics[width=0.6\textwidth, trim=50 120 50 150, clip=true]{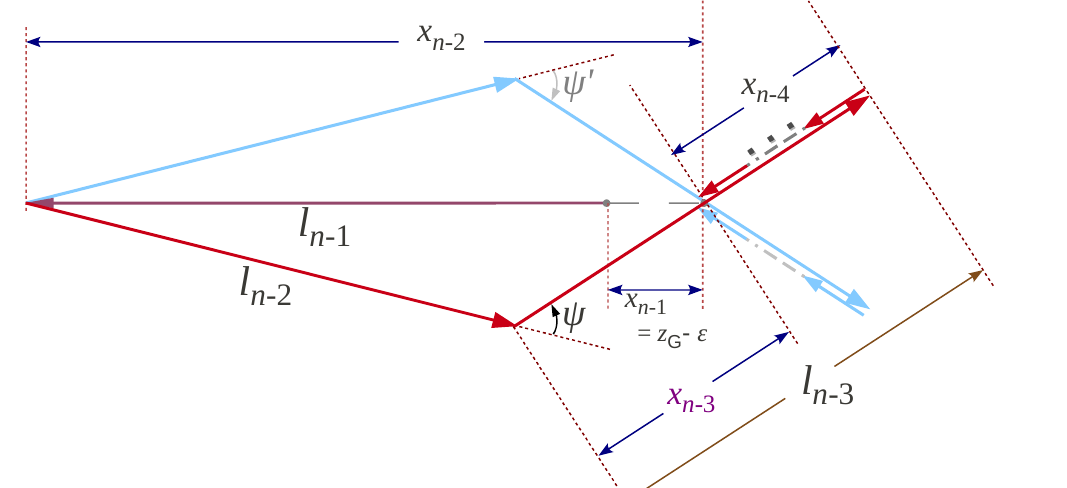}
\end{center}
%
At such values of $x_{*}$, 
we have $x_{n-1} = x_{n-2} + l_{n-1}$ and $x_p = |x_{p-1} - l_p|, ~p\in\{n-3,n-4,\cdots,1\}$. Thus in any configuration with these values of $x_*$ we have $\Theta^{+}_{l_p}(x_p,x_{p-1}) = \Theta^{-}_{l_p}(x_p,x_{p-1}) = \left\{ \begin{array}{l} 0, ~\text{for $p=n-1$} \\ \pi , ~\text{for $p\in\{n-3,n-4,\cdots,1\}$} \end{array} \right.$, and, $\Phi^{+}_{l_p}(x_p,x_{p-1}) = \Phi^{-}_{l_p}(x_p,x_{p-1}) = 0$ for $p=n-1,n-3,n-4,\cdots,1$, as illustrated above (see Note~\ref{note:plus-minus-equal}). 

However, we note that $\mathscr{S}_{S^\#_{n-2}, \mathbf{sg}^\#_{n-2}}(x^\zG_{n-2}) = \mathquote{+}$ and $\mathscr{S}_{S^\#_{n-2}, \mathbf{\hat{sg}}^\#_{n-2}}(x^\zG_{n-2}) = \mathquote{-}$ for $\#=I, II$. Moreover since $\underline{x}_{n-3} > 0$ (assuming $\zG$ is reachable), we have $x_{n-2} \neq x_{n-3} + l_{n-2}$ and $x_{n-2} \neq |x_{n-3} - l_{n-2}|$.
Thus, $\Theta^{+}_{l_{n-2}}(x_{n-2},x_{n-3}) = -\Theta^{-}_{l_{n-2}}(x_{n-2},x_{n-3}) \neq 0$ and $\Phi^{+}_{l_{n-2}}(x_{n-2},x_{n-3}) = -\Phi^{-}_{l_{n-2}}(x_{n-2},x_{n-3}) \neq 0$.
Using these (and Note~\ref{note:arm-angles}.ii.), one can easily verify that the angle between the segments of length $l_{n-2}$ and $l_{n-3}$ due to the configurations produced by the two inverse kinematics, $\textsf{IK}_{ \{ f^{\#}_{*-1}, \mathscr{S}_{S^{\#}_{*}\!, \mathbf{sg}^{\#}_{*}} \} }(z_\zG - \epsilon)$ and $\textsf{IK}_{ \{ f^{\#}_{*-1}, \mathscr{S}_{S^{\#}_{*}\!, \mathbf{\hat{sg}}^{\#}_{*}} \} }(z_\zG - \epsilon)$ (for $\# = I,II$), are the negative of each other (\emph{i.e.}, inverses of each other on the circle group $\mathbb{S}^1$). This is illustrated by the angle $\psi$ and $\psi'$ in the figure above.
Again, note that since $z_\zG - \epsilon < l_{n-3}$, if the set $\{z_\zG -\epsilon, l_{n-1}, l_{n-2}, \cdots, l_0\}$ is sorted such that the sorted elements are $l_\ordind{n} \geq l_\ordind{n-1} \geq \cdots \geq l_\ordind{0}$, then $l_{n-1} = l_\ordind{n}, l_{n-2} = l_\ordind{n-1}$ and $l_{n-3} = l_\ordind{n-2}$.
Thus, it follows from Proposition~\ref{prop:disconnected-configs}, that these two configurations are in disconnected components of $D_{R}^{-1}(z_\zG - \epsilon)$.
%
%
%

\vspace{0.05in}
The proof for $x_{n-1} = \lim_{\epsilon\rightarrow 0_+} z_\zA - \epsilon$ and $x_{n-1} = \lim_{\epsilon\rightarrow 0_+} z_\zF - \epsilon$ follows in a similar manner.
\end{quoteproof}

\section{Implementation and Results}

\begin{figure*}[h]
\centering
      \fbox{\includegraphics[width=1.0\textwidth, trim=0 0 0 0, clip=true]{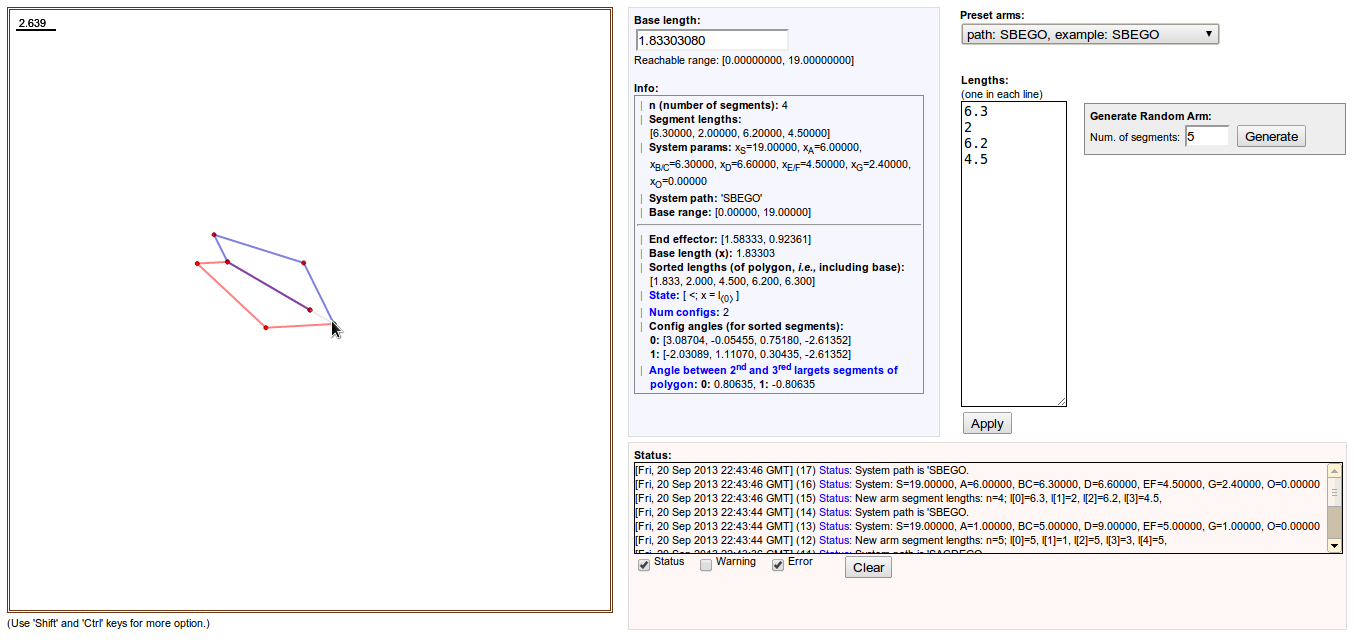}}
\caption{The Javascript and HTML5 based graphical user interface available at \url{http://hans.math.upenn.edu/~subhrabh/nowiki/robot_arm_JS-HTML5/arm.html} .}  \label{fig:JS-interface}
\end{figure*}

A simple analysis reveals that the complexity of computing the inverse kinematics 
is 
$O(n^2)$ (equation \ref{eq:theta-q} being computed for $q=0,1,\cdots,n-1$). Thus is it possible to implement the algorithm efficiently on a system-independent web-based platform with relatively low computational capabilities. In particular we chose to do a HTML5 + Javascript implementation and made it available on the web at one of the following links:
\begin{itemize} [topsep=0pt, partopsep=0pt]
 \setlength{\itemsep}{0pt} 
 \item[-] {\small \url{http://subhrajit.net/wiki/index.php?xURL=arm_JS-HTML5}}
 \item[-] {\small \url{http://hans.math.upenn.edu/~subhrabh/nowiki/robot_arm_JS-HTML5/arm.html}}
\end{itemize}
The interface is highly intuitive (Figure~\ref{fig:JS-interface}). The ``preset arms'' gives a list of preset arms that can be selected. Mouse position on the workspace on the left indicates the end effector position (holding the 'ctrl' key lets one move mouse at sub-pixel accuracy), and the program will compute the inverse kinematics to show one or two arm configurations depending on if the state block is of the form $[>; *]$ or $[<;*]$.
The `info' section displays information about the arm and the current state, and does validating checks on the number of configurations returned and whether or not the two configurations belong to different connected components.
The code is optimized to run on a Google chrome browser.

\section{Conclusion and Future Direction} \label{sec:discussion}

\begin{figure}[h]
\centering
\subfigure[In absence of obstacles there are two connected components of the constrained configuration space at $z'$.]{
      \includegraphics[width=0.45\textwidth, trim=20 20 20 20, clip=true]{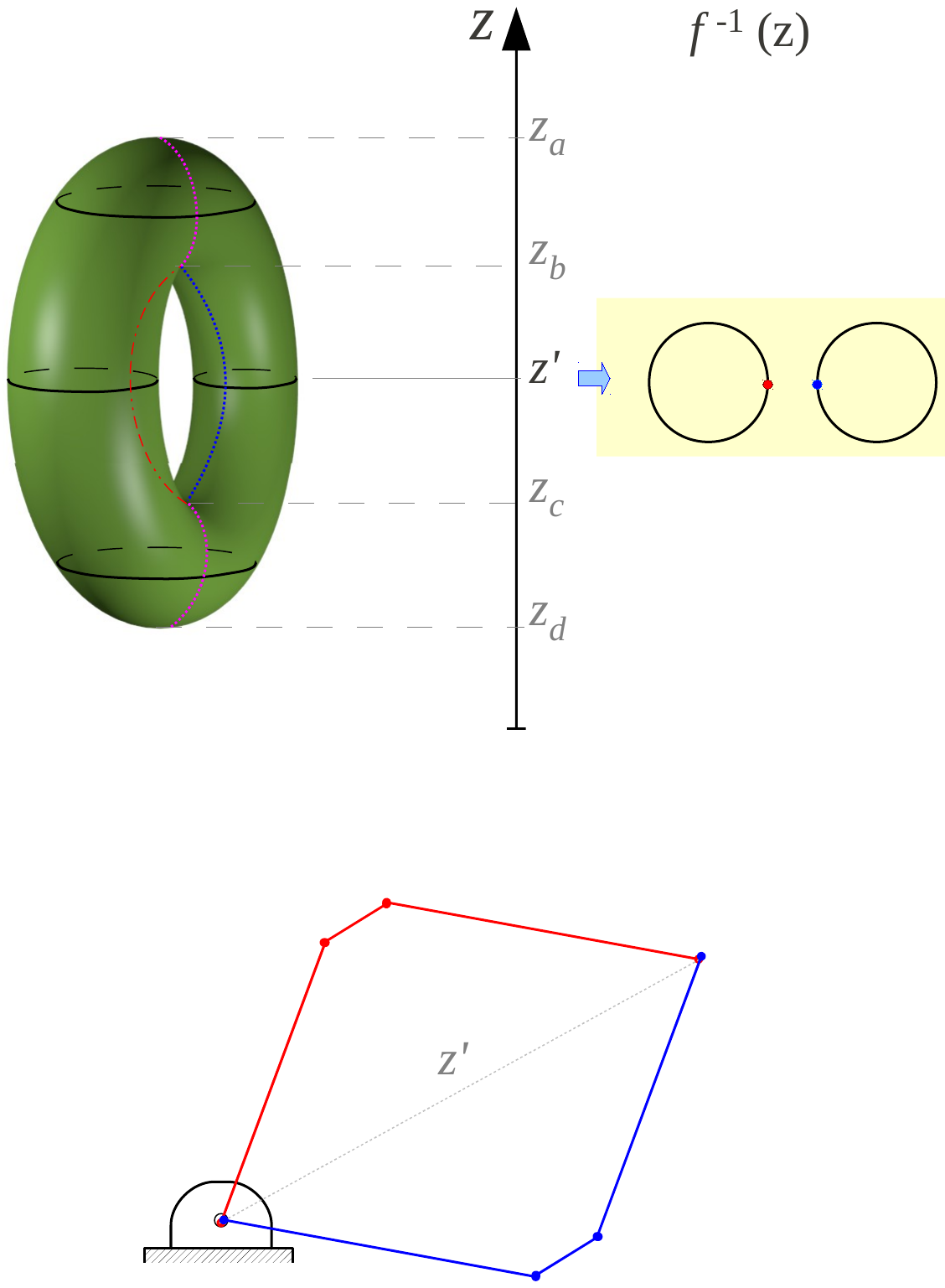}
}
\subfigure[In presence of obstacles there are three connected components of the constrained configuration space at $z'$.]{
      \includegraphics[width=0.45\textwidth, trim=20 20 20 20, clip=true]{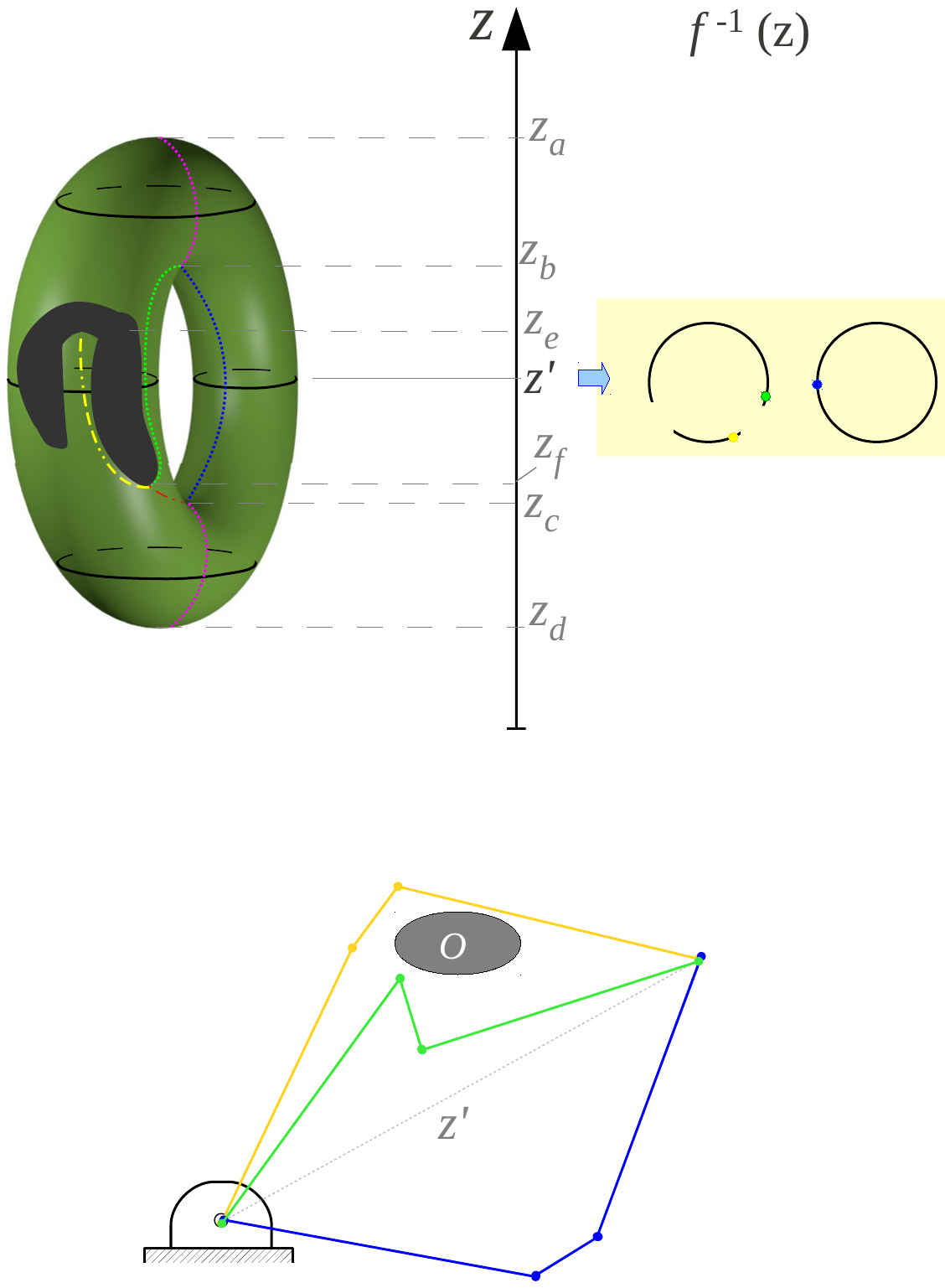}
}
\caption{An illustration (a low-dimensional analogy -- top row does not show actual configuration space) of how having obstacles in the environment can create punctures in $\mathfrak{R}$, which in turn can create more connected components in $D_{R}^{-1}(z)$, and hence new vital critical points.}  \label{fig:future}
\end{figure}

The design problem problem addressed in this paper involved sampling exactly one unique configuration for each connected component of the constrained configuration spaces, $D_{R}^{-1}(z)$, such that the sampled configurations vary continuously as $z$ is varied continuously.
We motivated the need for solving such a problem from the more practical problem of planning optimal path for the end effector of a robot arm. The said sampling approach would drastically reduce the dimensionality of the search domain in automated manipulation planning using search algorithms like A*.
The solution to this problem is made possible by a fundamental theoretical result involving classification of the configuration spaces of planar robot arms.


Moving forward, our next step will be to generalize this technique to the case when the environment contains obstacles. Identifying the \emph{vital critical points} (Difinition~\ref{def:vital-critical}) and designing inverse kinematics such that the system pass through the vital critical points will understandably be more difficult in that case, and will require numerical techniques.
Figure~\ref{fig:future} illustrates the challenges in achieving this.


\bibliographystyle{alpha}
\bibliography{IK}
\end{document}